\documentclass[10pt]{article}

\usepackage{amsmath}
\usepackage{amssymb}

\usepackage{graphicx}

\usepackage{cite}

\usepackage{color} 

\usepackage[justification=justified,singlelinecheck=off,font={small}]{caption}

\usepackage[caption=false]{subfig}
\date{\today}

\begin{document}

\begin{center}
{\Large
\textbf{On Robustness Analysis of Stochastic Biochemical Systems by Probabilistic Model Checking}
}
\\[5mm]
Lubo\v{s} Brim,
Milan \v{C}eska, 
Sven Dra\v{z}an,
David \v{S}afr\'{a}nek$^{\ast}$
\\[4mm]
Systems Biology Laboratory, Faculty of Informatics, Masaryk University, Brno, Czech Republic
\\
$\ast$ E-mail: safranek@fi.muni.cz
\end{center}
\thispagestyle{empty}

\section*{Abstract}
This report proposes a novel framework for a rigorous robustness analysis of stochastic biochemical systems. The technique is based on probabilistic model checking. We adapt the general definition of robustness introduced by Kitano to the class of stochastic systems modelled as continuous time Markov Chains in order to extensively analyse and compare robustness of biological models with uncertain parameters. The framework
utilises novel computational methods that enable to effectively evaluate the
robustness of models with respect to quantitative temporal properties and
parameters such as reaction rate constants and initial conditions.

The framework is applied to gene regulation as an example of a central biological
mechanism where intrinsic and extrinsic stochasticity plays crucial role due to
low numbers of DNA and RNA
molecules. 
Using our methods we have obtained a comprehensive and precise analysis of
stochastic dynamics under parameter uncertainty. Furthermore, we apply our
framework to compare several variants of two-component signalling networks from
the perspective of robustness with respect to intrinsic noise caused by low
populations of signalling components. We succeeded to extend previous studies
performed on deterministic models (ODE) and show that stochasticity may
significantly 
affect obtained predictions. Our case studies demonstrate that the framework can
provide deeper insight into the role of key parameters in maintaining the system
functionality and thus it significantly contributes to formal methods in
computational systems biology.

\section{Introduction} 

Robustness is one of the fundamental features of biological systems. According
to Kitano~\cite{Kitano2004} \textit{``robustness is a property that allows a
 system to maintain its functions against internal and external
 perturbations''}. To formally analyse robustness, we must thus precisely
identify model of a biological system and define formally the notions of a
system's function and its perturbations. In this paper, we propose a novel
framework for robustness analysis of stochastic biochemical systems. To this
end, inspected systems are described by means of stochastic biochemical kinetics
models, system functionality is defined by its logical properties, and system
perturbation is modelled as a change in stochastic kinetic parameters or initial conditions of the model.

Processes occurring inside living cells exhibit dynamic behaviours that can be
observed and classified as carrying out a certain function -- maintaining stable
concentrations, responding to a change of the environment, growing etc. Kinetic
models with parameters are used to formally capture cell dynamics. To observe
and analyse a dynamic behaviour on a kinetic model, all its numerical parameters
must be instantiated to a specific value.  This poses a challenge since the
precise values of all parameters (kinetic constants, initial concentrations,
environmental conditions etc.)  may not be known, may be known but without a
given accuracy of measurement or may in principle form an interval instead of
being a single value (e.g. non-homogeneous cell populations, different
structural conformations of a molecule leading to multiple kinetic rates etc.).
This implies that the behaviour of a kinetic model for a given single parametric
instantiation and its derived functionality may not provide an adequate result and
it is therefore unavoidable to take into account possible uncertainties, variance and
inhomogeneities.

The concept of robustness addresses this aspect of functional evaluation by
considering a weighted average of all behaviour across a \emph{space of
  perturbations} each altering the model parameters (hence its behaviour) in a
particular way and having a certain probability of occurrence. A general
definition of robustness was introduced by Kitano~\cite{Kitano2007}:
$$
R^{\mathcal{S}}_{\mathcal{A},\mathbf{P}} = \int_{\mathbf{P}}{\psi(p)D^{\mathcal{S}}_{\mathcal{A}}(p)dp}
$$
where $\mathcal{S}$ is the system, $\mathcal{A}$ is the function under scrutiny,
$\mathbf{P}$ is the space of all perturbations, $\psi(p)$ is the probability of
the perturbation $p\in \mathbf{P}$ and $D^{\mathcal{S}}_{\mathcal{A}}(p)$ is an
\emph{evaluation function} stating how much the function $\mathcal{A}$ is
preserved under a perturbation \textit{p} in the system $\mathcal{S}$.

For the macroscopic view as provided by the deterministic modelling framework
based on ordinary differential equations (ODEs), the concept of robustness has
been widely studied. There exist many mature analytic techniques based on static
analysis as well as dynamic numerical methods for effective robustness analysis
of ODE models. In circumstances of low molecular/cellular numbers such as in
signalling~\cite{ueda2007stochastic}, immunity reactions or gene
regulation~\cite{Gillespieetal05}, intrinsic and extrinsic noise plays an
important role and thus these processes are more faithfully modelled
stochastically. However, the existing methods and tools are not adequate for
rigorous and effective analysis of stochastic models with uncertain
parameters. In order to bridge this gap we adapt the concept of robustness to
stochastic systems.

The main challenge of the adaptation lies in the interpretation of the
evaluation function $D^{\mathcal{S}}_{\mathcal{A}}(p)$. We discuss several
definitions of the evaluation function that give us different options how to
quantify the ability of the system to preserve the inspected functionality under
a parameters perturbation. We show how absolute and relative robustness
of the stochastic systems can be captured and analysed using our
framework.

Semantics of stochastic biochemical kinetics models can be defined by \emph{Continuous
Time Markov Chains} (CTMCs) where the evolution of the probability density vector
describing the population of particular species is given by the chemical master
equation (CME)~\cite{Gillespie1977}.  A function of a system in the biological
sense is any intuitively understandable behaviour (e.g., stability of \textit{ERK}
signal effector population in high concentration observed in a given time
horizon). In order to define the robustness of a system formally we need to make
precise the intuitive and informal concept of functionality.  Our framework
builds on the formal methods where the functionality of a system is expressed
indirectly by its logical properties. This leads to a more abstract approach
emphasising the most relevant aspects of a system function and suppressing less
important technicalities. We use stochastic temporal logics, namely the bounded time
fragment of \emph{Continuous Stochastic Logic (CSL)}~\cite{Aziz1996csl} further
extended with \emph{rewards}~\cite{Kwiatkowska2006rewards} (e.g.,
$\mathsf{P}_{\geq 0.9}[ \mathsf{G}^{[t_1,t_2]} (ERK > \textit{high})]$).  To broaden
the scope of possibly captured functionality we extend CSL with a class of
\emph{post-processing functions} defined over probability density vectors. We
show that the bounded fragment of CSL with rewards and post-processing functions
can adequately capture many biological relevant behaviours that are recognisable
in finite time intervals.

Our framework is based on probabilistic model checking techniques that compute
the probability with which a given CTMC satisfies a given CSL formula. The
computation can be conducted using Monte Carlo based methods such as Gillespie's
stochastic stimulation algorithm~\cite{Gillespie1977} or numerical methods such
as uniformisation~\cite{Stewart2009uniformization}. Although Monte Carlo based
methods (often denoted as statistical model checking) can produce detailed
simulations for stochastically evolving biochemical systems, computing a
statistical description of their dynamics that is necessary for evaluating
$D^{\mathcal{S}}_{\mathcal{A}}(p)$, such as the probability density, mean, or
variance, requires a large number of individual simulations. Moreover, if
$\mathcal{A}$ describes a behaviour that occurs rarely, the evaluation of
$D^{\mathcal{S}}_{\mathcal{A}}(p)$ requires an extremely large number of
simulations to be performed to obtain sufficient accuracy. As was shown
in~\cite{munsky2006finite} in such situations numerical methods are
substantially more efficient. Since rigorous stochastic robustness analyses 
may require to compute precise probabilities of all behaviours,
we build our framework on probabilistic numerical methods.  

To analyse the robustness of the CTMC $\mathcal{C}$ with respect to the CSL
formula $\Phi$ over the space of perturbations $\mathbf{P}$, which can be discrete
but still very large or continuous and thus infinite, one needs to efficiently
compute or approximate the evaluation function $D^{\mathcal{C}}_{\Phi}$, i.e,
the values $D^{\mathcal{C}}_{\Phi}(p)$ for all $p\in \mathbf{P}$. One of the
possible approaches (recently used in~\cite{Bartocci2013}) is to effectively
sample the perturbation space $\mathbf{P}$ and use standard statistical or
numerical methods to obtain values in grid points. These values can be
afterwards interpolated linearly or polynomially. Using adaptive grid refinement
such an approach provides an arbitrary degree of precision. A disadvantage of
this method is the fact that the obtained result is an approximation not
providing any minimal and maximal upper bounds. Therefore such an approach can
neglect sharp changes or discontinuities in the landscape of the evaluation
function $D^{\mathcal{C}}_{\Phi}$.

It is worth noting that the evaluation function can be discontinuous or may
change its value rapidly on a very small perturbation interval in situations
when the given CSL formula contains nested probability operators. In particular,
this is inevitable to formulate hypothesis requiring the detailed temporal
program~\cite{Alon2004} of the biological system (e.g., temporal ordering of events). The actual
shape of the evaluation function arises from the combination of such a formula
and the particular model. Especially, high sensitivity of a model to the
perturbed parameter can intensify rapid changes in the evaluation function. An
example of a formula with a nested probability operator is mentioned in Section~\ref{sec:functionality}.

To evaluate the function $D^{\mathcal{C}}_{\Phi}$ we employ in our framework
another approach that is based upon our min-max approximation method recently
published in~\cite{CAV2013}. The method guarantees strict upper and lower
estimates of $D^{\mathcal{C}}_{\Phi}(p)$ without neglecting any sharp changes or
discontinuities. This method exploits numerical techniques for probabilistic
model checking, can provide arbitrary degree of precision and thus can be
considered as an orthogonal approach to the adaptive grid refinement.  The
framework further extends the min-max approximation to a more general class of
stochastic biochemical models (i.e., incorporation of stochastic Hill kinetics)
and a more general class of quantitative properties (i.e., including
post-processing functions) and allows us to compute the robustness of such
systems. In our framework we provide the user not only a numerical value giving
the robustness of the system but possibly also a landscape visualisation of the
evaluation function.

We demonstrate the applicability of the proposed method by means of two
biological case studies -- a model predicting dynamics of a gene regulatory
circuit controlling the $G_1/S$ phase transition in the cell cycle of mammalian
cells, and two models representing different topologies of a general
two-component signalling mechanism present in procaryotic cells. Both cases are
examples of cellular processes where stochasticity plays a crucial role
especially because of low numbers of molecules involved.

The former case study exploits the usability of the method to analyse
bistability (and its robustness) in the stochastic framework and thus provides a
stochastic analysis analogy to the study presented in~\cite{Swatetal04} under
the deterministic (ODE) setting. Robustness is employed to characterise
parameterisations of the model with respect to the tendency of the molecule
population to choose one of the possible steady states irreversibly deciding
whether the cell will or will not commit to \textit{S}-phase. The results show that
intrinsic and extrinsic noise caused by randomness in protein-DNA
binding/unbinding events and other processes controlling the chemical affinity
of involved molecules can significantly affect the cell decision. In our model,
the intrinsic noise of chemical reactions is inherently captured by stochastic
mass action kinetics whereas the extrinsic noise is considered by means of
parameter uncertainty.

The latter case study focuses on analysing the effect of intrinsic noise on the
signalling pathway functionality. In particular, two topologically different
variants of a two-component signalling pathway are exploited for different levels
of input signal and different levels of intrinsic noise appearing in
transcription of the two signalling components. The considered topologies have
been compared in the previous study presented by Steuer et
al.~\cite{steuer2011robust} where robustness has been analysed in the
setting of deterministic (ODE) models. Here the signalling mechanism is
remodelled in the stochastic setting and robustness is employed to quantify
under which circumstances the individual topologies are less amenable to
intrinsic noise of the underlying protein transcription mechanism. The results
show that the stochastic approach can uncover facts unpredictable in the
deterministic setting.

Formal analysis of complex stochastic biological systems employing both 
the numerical and the statistical methods generally suffers from extremely
high computational demands. These computational demands are even more critical 
if we need to analyse systems with uncertain parameters which is also the 
case of our framework. However, our framework has been designed in order to be adapted to high performance 
computing platforms (e.g. multi-core workstations and 
massively parallel general-purpose graphic processing units)  and also to be successfully 
combined with existing acceleration methods, see e.g.~\cite{munsky2006finite,Henzinger2009,Didier}. 
Although the acceleration is a subject of our future research (inspired by our previous results~\cite{BBS10}), we already employ the fact 
that the min-max approximation method can be efficiently parallelised. In the second case study where the analysis 
of the inspected perturbation space requires an extensive numerical computation, we utilise a high performance multi-core workstation to  achieve the acceleration. Fundamentally different approaches to overcome the time 
complexity of stochastic analyses of complex biological systems build on a moment closure 
computation and on a fluid approximation, see e.g.~\cite{Verena2013,Bortolussi2012}. These approaches 
are briefly discussed in the related work.

The main contributions of this paper can be summarised in the following way:
\begin{enumerate}
\item The adaptation of the general concept of robustness of
  Kitano~\cite{Kitano2007} to the class of stochastic systems modelled by
  CTMCs. The key step of the adaptation is a definition of the evaluation
  function that reflects the quantitative aspects of stochastic models and their
  behaviours. We discuss several definitions of the function allowing for
  different ways of capturing stochastic robustness.
\item Introduction of a novel framework based on formal methods to evaluate
  robustness of the stochastic system with respect to the functionality given by
  a stochastic temporal property and to \textit{perturbations in reaction rate
    parameters and initial conditions}. The framework significantly extends the
  min-max approximation method published in~\cite{CAV2013}, namely with the
  support for Hill kinetics and post-processing functions.
\item Demonstration of the fact that our concept of robustness can capture and
  quantify the ability of the stochastic systems to maintain their
  functionality.  We apply our framework to two biologically relevant case
  studies. Namely, it is the gene regulation of mammalian cell cycle where we
  explore the impact of stochasticity in low molecule numbers to bistability of
  a regulatory circuit controlling G1/S transition and analysis of noise
  behaviour in different topologies of two-component signalling systems. The case
  studies show that our framework provides deeper understanding of how the validity
  of an inspected hypotheses depends on reaction rate parameters and initial
  conditions.
	
 \end{enumerate}

\subsection{Related work}

The discussion on related work can be roughly divided into two parts. 
First, we summarise the existing methods
for parameter exploration and robustness analysis of stochastic models.
Second, we briefly mention the methods and tools allowing for robustness analysis of ODE models. 

In the field of stochastic models, parameter estimation methods and the concept
of robustness are not as established yet as in the case of ODE models. We have
recently published a method~\cite{CAV2013} where the CSL model checking
techniques are extended in order to systematically explore the parameters of
stochastic biochemical kinetic models. In~\cite{Verena2012} a CTMC is explored
with respect to a property formalised as a deterministic timed automaton
(DTA). It extends~\cite{Verena2011} to parameter estimation with respect to the
acceptance of the DTA. Most approaches to parameter
estimation~\cite{Reinker,Verena2011,Petzold2012} rely on approximating the
maximum likelihood. Their advantage is the possibility to analyse infinite state
spaces~\cite{Verena2011} (employing dynamic state space truncation with
numerically computed likelihood) or even models with no prior knowledge of
parameter ranges~\cite{Petzold2012} (using Monte-Carlo optimisation for
computing the likelihood). In~\cite{Verena2013} the moment closure approach is
considered to capture the distribution of highly populated species in combination
with discrete stochastic description for low populated species. The method is
able to cope with multi-modal distributions appearing in multi-stable
systems. The method introduced in~\cite{Bortolussi2012} exploits fluid (limit)
approximation techniques and in that way enables an alternative approach to CSL
model checking of stochastic models. Despite the computational efficiency, a
shared disadvantage of all the mentioned methods is that they rely on
approximations applicable only to models that include highly populated
species. This is not the case of, e.g., gene regulation dynamics.

Approaches based on Markov Chain Monte-Carlo sampling and Bayesian
inference~\cite{Wilkinson,Clarke,Cago} can be extended to sample-based
approximation of the evaluation function, but at the price of undesired
inaccuracy and high computational demands~\cite{Bernardini,Paolo}. Compared to
these methods, our method provides the upper and lower bounds of the result
which makes it more reliable and precise but at the price of higher computational
demands. The most relevant contribution to this domain has been recently
introduced by Bartocci et al.~\cite{Bartocci2013}. To our best knowledge this is
the only related work addressing robustness of stochastic biochemical
systems. The work is based on the idea to directly adapt the concept of behaviour
oriented robustness to stochastic models. Individual simulated trajectories of
the CTMC are locally analysed with respect to a formula of Signal Temporal Logic
(STL), a linear-time temporal logic interpreted on simulated time sequences. For
each simulated trajectory, the so-called satisfaction degree representing the
distance from being (un)satisfied is computed, thus resulting into a randomly
sampled distribution of the satisfaction degree. This distribution thus gives
modellers another source of information in addition to probability of formula satisfaction
(percentage of valid trajectories in the sampled set). In comparison, our method
directly (and exactly) computes the probability of formula satisfaction for a
different kind of temporal logic -- the branching-time CSL logic. This allows to
express more intricate properties that require branching time, e.g.,
multi-stability. On the other hand, our method as conceptually based on
transient analysis does not allow to compute local analysis of individual
trajectories, i.e., to obtain a satisfaction degree would require non-trivial
elaboration at the level of numerical algorithms.

In the domain of ODE models, there exist several analytic methods for effective
analysis under parameter uncertainty. They build on static analysis
(stoichiometric analysis, flux balance analysis) as well as dynamic numerical
methods (simulation, monitoring by temporal formulae, sensitivity analysis)
implemented in tools
(e.g.~\cite{Hoops2006copasi,Loew2001vcell,Fages2004biocham}). Robustness
analysis with respect to functionality specified in terms of temporal formulae
has been introduced recently~\cite{Fainekos2009robustness,Rizketal09}. There
exist two major approaches how to define and analyse
robustness. 
If only parameters of the model are perturbed, we speak of a \textit{behaviour
  oriented approach} to robustness. This approach has been explored by Fainekos
\& Pappas~\cite{Fainekos2009robustness}, further extended by A. Donz\'{e} et
al.~\cite{Donze2010robust} and implemented in the toolbox
Breach~\cite{Donze2010breach}. Another option could be to perturb the model
structure i.e. the reaction topology, as this is done in many gene knock-out
biological experiments. Such changes are in principle discrete and the problem
of robustness computation for such perturbations would reduce to solving many
individual instances of the same problem for each discrete topology. However
identifying model behaviour shared among individual perturbations can lead to
more efficient analysis~\cite{Barnat12coloredmc}. 

Yet another way to look at perturbations is from the perspective of property
uncertainty. If the system is considered fixed and all parameters exactly known,
the uncertainty then lies in the property of interest. For a specific property
such as ``The concentration of X repeatedly rises above 10 and drops below 5
within the first 20 minutes'' where all three numerical constants can be
altered, we explore how much would they have to be altered in order to affect
the property validity in the given model. This approach has been adopted for
ODEs by F. Fages et al.~\cite{Rizketal09} and implemented in the tool
BIOCHAM~\cite{Fages2004biocham}. When only parameters of the property are
perturbed, it is the case of a \textit{property oriented approach} to
robustness.

\section{Methods}
\label{sec:methods}

\subsection{Methodology Overview}
In this paper we propose a formal framework that allows to analyse the
robustness of stochastic biochemical systems with respect to a space of perturbing
parameters. The framework consists of the following objects:

\begin{itemize}
\item \emph{a finite state stochastic biochemical system given by a set of chemical
  species participating in a set of chemical reactions}
  
   Each of the reaction is
  associated with a stochastic rate function that for a fixed stochastic
  rate constant returns  the rate of the reaction. To formalise such system we use a
  population based finite state continuous time Markov Chain (CTMC), i.e, a
  state of the CTMC is given by populations of particular species and the
  evolution of the CTMC is driven by the chemical master equation (CME)~\cite{Gillespie1977,Didier}.

\item \emph{a perturbation space defined by a Cartesian product of uncertain stochastic
  rate constants given as value intervals with minimal and maximal
  bounds}
  
  Additionally, the perturbation space may also be expanded by initial
  conditions of the system (i.e, interval for the size of a population of a
  particular species) encoded in the initial state of the CTMC. The given
  stochastic system and the perturbation space induce a \emph{set of parameterised}
  CTMCs.

\item \emph{set of paths that describe the evolution of a fully instantiated
  stochastic system (i.e., in which all stochastic rate constants and the
  initial state are specified) over time}
  
  For a state of the system and a finite
  time there is a unique probability measure of all paths starting in that state
  that defines  \emph{probability distribution} over states occupied by the system
  at the given time. Each perturbation from a given perturbation space possibly
  leads to a different probability distribution.

\item \emph{stochastic temporal property interpreted over the paths and states of CTMC
  enabling to specify an a priori given quantitative hypothesis about the
  system} 
  
   We primarily focus on the \emph{bounded time fragment of Continuous
    Stochastic Logic} (CSL)~\cite{Aziz1996csl} further extended with
  rewards~\cite{Kwiatkowska2006rewards}. For most cases of biochemical
  stochastic systems the bounded time restriction is adequate since a typical
  behaviour is recognisable in finite time. Additionally, we also consider
  properties given by a class of post-processing functions defined over
  probability distributions at the given finite time.

\end{itemize}

The main goal of our framework is to analyse how the validity of an \emph{a
  priori} given hypothesis expressed as a temporal property depends on uncertain
parameters of the inspected stochastic system.  For this purpose we adapt the
general definition of robustness~\cite{Kitano2007} to the class of stochastic
systems.  While the concept of robustness is well established for deterministic
systems~\cite{Donze11rob-behaviour,Rizk09rob-property}, it has not been
adequately addressed for stochastic systems. The key difference is the fact that
evolution of a stochastic system is given by a set of paths in contrast to a
single trajectory as in the case of a deterministic system. Hence a stochastic
system at the given time is described by a probability distribution over states
of the corresponding CTMC in contrast to the single state representation of a
deterministic system. Therefore, the definition of robustness for stochastic
systems requires a more sophisticated interpretation of the evaluation function
that determines how the quantitative temporal property is preserved under a
perturbation of the system's parameters.

Similarly to Kitano, we define robustness of stochastic systems as the integral
of an evaluation function. In our case the evaluation function
$D^{\mathcal{C}_p}_{\Phi}$ for each parameter point \textit{p} from the
inspected perturbation space $\mathbf{P}$ returns the quantitative model checking
result for the respective CTMC $\mathcal{C}_p$ and the given property $\Phi$. We
show how robustness can be effectively over/under-approximated for a class of
quantitative temporal properties using new techniques for model checking of
parameterised CTMCs. Moreover, if the property can be expressed using only the
bounded time fragment of CSL with rewards (i.e., without post-processing
functions) we can extend the approach to global quantitative model checking
techniques. They enable to compute the model checking result for all states of a
CTMC with the same price as for a single state and thus to analyse the perturbation
of initial conditions in a much more effective way. Finally, we demonstrate how
robustness can capture and quantify the ability of a stochastic system to
maintain its functionality described by such class of properties.

Since the inspected perturbation space is in principle dense the set of
parameterised CTMCs to be explored is infinite. It is thus not possible to
compute the model checking result for each CTMC individually. The
straightforward approach to overcome this problem could be to sample points from
the perturbation space and use existing model checking techniques for fully
instantiated CTMCs. That way we can obtain precise model checking results in
the grid points and then interpolate them linearly or polynomially. Although an
adaptive grid refinement could provide an arbitrary degree of precision, it does
not guarantee strict lower and upper bounds. Hence such an approach could
neglect sharp changes or discontinuities of the evaluation function. Since we
want  to guarantee strict bounds of obtained results, we extend our
previously published method~\cite{CAV2013}. This method allows to compute strict
minimal and maximal bounds on the quantitative model checking results for all
CTMCs $\left\{ \mathcal{C}_p \mid p \in \mathbf{P} \right\}$ for a given
perturbation space~$\mathbf{P}$.

\subsection{Models} 

The formalism used to model a biochemical system is essential since it not only
dictates the possible behaviours that may or may not be captured, but also
determines the means of detecting them. ODEs enable the study of large ensembles
of molecules in population count and species diversity since they abstract from
the individualistic properties of each molecule such as position or its
stochastic behaviour and take as its variables only concentrations of each
species. Stochastic models such as CTMCs abstract positions of molecules but
maintain their individual reactions. Even more detailed models such as Brownian
dynamics which keep track of positions but abstract from the geometry and
orientation of each molecule could be used. However as the amount of information
about each individual molecule increases the computational complexity of proving
some property to hold over all the behaviours of a model becomes quickly
infeasible even for small models.

In our framework we focus on stochastic biochemical systems that can be
formalised as a finite state system~$\mathcal{S}$ defined by a set of \textit{N}
\emph{chemical species} in a well stirred volume with fixed size and fixed
temperature participating in \textit{M} \emph{chemical reactions}.  The number $X_i$ of
molecules of each species $S_i$ has a specific bound and each reaction is of the
form $u_1 S_1 + \ldots + u_N S_N \longrightarrow v_1 S_1 + \ldots + v_N S_N$
where $u_i, v_i \in \mathbb{N}_0$ represent \emph{stoichiometric coefficients}.

A \emph{state} of a system in time $t$ is the vector $\mathbf{X}(t) = (X_1(t),
X_2(t), \ldots, X_N(t) )$. When a single reaction with index $r \in \{1, \ldots,
M\}$ with vectors of stoichiometric coefficients $U_r$ and $V_r$ occurs the
state changes from $\mathbf{X}$ to $\mathbf{X}' = \mathbf{X} - U_r + V_r$, which
we denote as $\mathbf{X} \stackrel{r}{\rightarrow} \mathbf{X}'$. For such
reaction to happen in a state $\mathbf{X}$ all reactants have to be in
sufficient numbers and the state $\mathbf{X}'$ must reflect all species
bounds. The \emph{reachable state space} of~$\mathcal{S}$, denoted as
$\mathbb{S}$, is the set of all states reachable by a finite sequence of
reactions from \emph{an initial state} $\mathbf{X}_0$.
The set of indices of all reactions changing the state $\mathbf{X}_i$ to the
state $\mathbf{X}_j$ is denoted as $\mathsf{reac}(\mathbf{X}_i,\mathbf{X}_j) =
\{r \mid \mathbf{X}_i \stackrel{r}{\longrightarrow} \mathbf{X}_j \}
$. Henceforward the reactions will be referred directly by their indices.

According to~\cite{Gillespie1977,Didier} the behaviour of a stochastic system
$\mathcal{S}$ can be described by the CTMC $\mathcal{C} = (\mathbb{S},
\mathbf{X}_0, \mathbf{R})$ where the transition matrix
$\mathbf{R}(\mathbf{X}_i,\mathbf{X}_j)$ gives the probability of a transition
from $\mathbf{X}_i$ to $\mathbf{X}_j$. Formally, the transition matrix is
defined as:
$$\mathbf{R}(\mathbf{X}_i,\mathbf{X}_j) \stackrel{def}{=} \sum_{r \in \mathsf{reac}(\mathbf{X}_i,\mathbf{X}_j) } f_r(\mathbf{k}_r, \mathbf{X}_i)$$ 
where $f_r$ is a \emph{stochastic rate function} and $\mathbf{k}_r$ is a vector
of all numerical parameters occurring in $f_r$ such as a \emph{stochastic rate
  constant} $k_r$, stoichiometry exponents, Hill coefficients etc.

In case of mass action kinetics the stochastic rate function has the simple form
of a polynomial of reacting species populations. That is $ f_r(\mathbf{k}_r,
\mathbf{X}_i) = k_r \cdot C_{r,i}$ where $C_{r,i} \stackrel{def}{=} \prod_{l =
  1}^{N}\binom{\mathbf{X}_{i,l}}{u_l}$ corresponds to the population dependent
term such that $\mathbf{X}_{i,l}$ is the \textit{l}th component of the state
$\mathbf{X}_{i}$ and $u_l$ is the stoichiometric coefficient of the reactant
$S_l$ in reaction $r$. However, sometimes the mass action kinetics is not sufficient, 
especially, when the reactions are not elementary but are rather an abstraction
of several reactions with unknown precise dynamics (e.g. gene transcription) or
if including all elementary reactions would cause the analysis to be
computationally infeasible. In such cases dynamics are typically approximated by
Hill functions~\cite{Hill1910}, a quasi-steady-state
approximation~\cite{Madsen2012} of the law of mass conservation. For sake of
simplicity of our presentation we will further assume that for each reaction $r$
the vector $\mathbf{k}_r$ is one-dimensional and thus $\mathbf{k}_r = k_r$, the
proposed methods can however be directly used also for multi-dimensional vectors
of constants. To comply with standard notation in the area of CTMC analysis
henceforward the states $\mathbf{X}_i \in \mathbb{S}$ will be denoted as $s_i$.

The probability of a transition from state $s_i$ to $s_j$ occurring within
\textit{t} time units is $1 - e^{-\mathbf{R}(s_i,s_j) \cdot t}$, if such a
transition cannot occur then $\mathbf{R}(s_i,s_j) = 0$. The time before any
transition from $s_i$ occurs is exponentially distributed with an overall
\textit{exit rate} $E(s_i)$ defined as $E(s_i) = \sum_{s_j \in
  \mathbb{S}}{\mathbf{R}(s_i,s_j)}$. A path $\omega$ of CTMC $\mathcal{C}$ is a
non-empty sequence $\omega=s_0,t_0,s_1,t_1\ldots$ where $\mathbf{R}(s_i,s_j)>0$
and $t_i\in \mathbb{R}_{\geq 0}$ is the amount of time spent in the state $s_i$
for all $i \geq 0$.  For all $s\in \mathbb{S}$ we denote by
$Path^{\mathcal{C}}(s)$ the set of all paths of $\mathcal{C}$ starting in
state~$s$. There exists the unique probability measure on
$Path^{\mathcal{C}}(s)$ defined, e.g., in~\cite{Kwiatkowska2007}. Intuitively,
any subset of $Path^{\mathcal{C}}(s)$ has the unique probability that can be
effectively computed. For the CTMC $\mathcal{C}$ the transient state
distribution $\pi^{\mathcal{C},s,t}$ gives for all states $s'\in \mathbb{S}$ the
transient probability $\pi^{\mathcal{C},s,t}(s')$ defined as the probability,
having started in the state \textit{s}, of being in state $s'$ at the finite time \textit{t}.

\subsection{Perturbations} 

In our approach we have focused on the behavioural approach for stochastic
systems and thus we will now define a set of perturbed stochastic systems and
their CTMCs. Let each stochastic rate constant $k_r$ have a value interval
$[k_r^{\bot},k_r^{\top}]$ with minimal and maximal bounds expressing an
\emph{uncertainty range} or \emph{variance} of its value. A \emph{perturbation space}
$\mathbf{P}$ induced by a set of stochastic rate constants $k_r$ is defined as
the Cartesian product of the individual value intervals $\mathbf{P} =
\prod_{r=1}^{M}[k_r^{\bot},k_r^{\top}]$. A single \emph{perturbation point} $p
\in \mathbf{P}$ is an \textit{M}-tuple holding a single value of each rate
constant, i.e., $p = (k_{1_p}, \ldots , k_{M_p})$.

A stochastic system $\mathcal{S}_p$ with its stochastic rate constants set to
the point $p\in \mathbf{P}$ is represented by a CTMC $\mathcal{C}_p =
(\mathbb{S}, s_0, \mathbf{R}_p)$ where transition matrix $\mathbf{R}_p$
is defined as:
 $$\mathbf{R}_p(s_i,s_j) \stackrel{def}{=} \sum_{r \in \mathsf{reac}(s_i,s_j) } f_r(k_{r_p},  s_i)$$
A \emph{set
  of parameterised} CTMCs induced by the perturbation space $\mathbf{P}$ is defined
as $\mathbf{C} = \{ \mathcal{C}_p \mid p \in \mathbf{P} \}$. 

Additionally, we consider the perturbation of initial conditions of the
stochastic system that are represented by different initial states of the
corresponding CTMC. In this case we extend the perturbation space such that a
single perturbation point $p\in \mathbf{P}^e = \mathbb{I} \times \mathbf{P}$
where $\mathbb{I} \subseteq \mathbb{S}$ is an \textit{M}+1-tuple holding a
single value of an initial state and a single value of each rate constant, i.e., $p =
(s_p,k_{1_p}, \ldots , k_{M_p})$ and CTMC $\mathcal{C}_p =(\mathbb{S}, s_p,
\mathbf{R}_p)$.

\subsection{Functionality} 
\label{sec:functionality}

To be able to automatically analyse a system's function $\mathcal{A}$ under
scrutiny there must be a formal way of expressing a function of a system. A
function of a system in the biological sense is any intuitively understandable
behaviour such as response, homoeostasis, reproduction, respiration or growth. It
can be a high level concept such as chemotaxis as well as a low level one
e.g. reaching of a state with a given number of molecules of a specific
species.

The inspected function can usually be described by a \textit{property} that is
understood as an abstraction of a system's behaviour expressed in some temporal
logic and given as a formula of that logic. Unlike the intuitive concept of a
biological function mentioned above, a property may be formally verified over a
formal model of a system and proven to hold or to be violated. Since the concept
of robustness builds on the notion of a function that can be measured,
we focus on a quantitative logic for stochastic systems. We use \emph{continuous
  stochastic logic} (CSL) \cite{Aziz1996csl,Baier2003mcctmc} extended with
\emph{reward} operators \cite{Kwiatkowska2006rewards}. Reward operators allow us
to further broaden the scope of possibly captured behaviour. They enable to
express properties such as the probability of a system being in the specified
set of states over a time interval or the probability that a particular reaction
has occurred.

Full CSL with rewards can express properties concerning a system in near future
as well as the infinite steady state situation. In this paper we focus only on
the \emph{bounded time fragment of CSL}. This fragment allows us to speak only
about behaviour within a finite time horizon. For most cases of biochemical
stochastic systems, such as intracellular reaction cascades or multi-cellular
signalling, the bounded time restriction is adequate since a typical behaviour
is recognisable within finite time intervals~\cite{KwiatkowskaBIO}.

As we show later on, there exist several biologically relevant properties that
cannot be directly expressed by CSL with rewards. Therefore, we employ a class of
post-processing functions to specify and analyse robustness of stochastic
systems with respect to such properties. The key idea of these functions is to
process and aggregate the transient state distribution at the given finite
time.

Let $\mathcal{C} =\left(\mathbb{S}, s_0, \mathbf{R}, L \right)$ be a labelled
CTMC such that \textit{L} is a labelling function which assigns to each state $s \in
\mathbb{S}$ the set \textit{L(s)} of atomic propositions that are valid in state
$s$.  We consider the specification of the inspected property using the bounded time fragment of CSL with rewards and post-processing functions. The syntax of this logic is defined in the following way. A state formula $\Phi$ is given as 
$$\Phi::= \mathsf{true} \mid a \mid \neg\Phi \mid \Phi \wedge \Phi \mid
\mathsf{P}_{\sim p}[\phi] \mid \mathsf{R}_{\sim r}[\mathsf{C}^{\leq t}] \mid
\mathsf{R}_{\sim r}[\mathsf{I}^{=t}] \mid \mathsf{E}_{\sim r}[\mathsf{I}^{=t}]$$
where $\phi$ is a path formula given as $\phi::= \mathsf{X}\mbox{~}\Phi \mid$
$\Phi\mbox{~}\mathsf{U}^{I}\mbox{~}\Phi$, \emph{a} is an atomic proposition,
$\sim \in \left\{ <, \leq, \geq, > \right\}$, $p \in [0,1]$ is a probability, $r
\in \mathbb{R}_{\geq 0}$ is an expected reward and $I = [a,b]$ is a bounded time
interval such that $a,b \in \mathbb{R}_{\geq 0} \wedge a \leq b$.  Path
operators $\mathsf{G}$ (always) and $\mathsf{F}$ (eventually) are derived in the
standard way using the operator $\mathsf{U}$.  In order to specify properties
containing rewards ($\mathsf{R}_{\sim r}[\mathsf{C}^{\leq t}]$ is the
\textit{cumulative reward} acquired up to time \textit{t}, $\mathsf{R}_{\sim
  r}[\mathsf{I}^{= t}]$ is the \textit{instantaneous reward} in time \textit{t})
the CTMC~$\mathcal{C}$ is enhanced with reward (cost) structures. Two types of
reward structures can be used, \emph{a state reward} and \emph{a transition
  reward}. For sake of simplicity, we consider in this paper only state rewards,
however, the proposed methods can be easily extended to transition rewards as
well. The state reward $\rho(s)$ defines the rate with which a reward is
acquired in state $s\in \mathbb{S}$. A~reward of $t \cdot \rho(s)$ is acquired
if $\mathcal{C}$ remains in state \textit{s} for \textit{t} time units.

Since the function $\rho$ has to be defined before the actual analysis of the
CTMC, the rewards for particular states have to be known prior to the
specification of the property. This fact limits the class of properties that can
be expressed using such structures. For example, noise expressed by a
\textit{mean quadratic deviation} (\textit{mqd}) of the population probability
distribution of a species at a given time cannot be specified using CSL with
rewards. To compute the \textit{mqd} we need to know the mean of the
distribution to be able to obtain the corresponding coefficients and encode them
into state rewards.

To overcome this problem we introduce the abstract state operator
$\mathsf{E}_{\sim r}[\mathsf{I}^{=t}]$ which \textit{evaluates} the state
distribution $\pi^{\mathcal{C},s_0,t}$ at the given time instant \textit{t} by a
user provided real-valued \textit{post-processing} function
$Post(\pi^{\mathcal{C},s_0,t})$ and compares it to $\sim r$. At the end of this
section we show how to define $Post$ in order to specify biologically relevant
properties such as noise using the \textit{mqd}. The \textit{mqd} is also used 
in the second case study to analyse a noise in different variants of signalling 
pathways.  

The formal semantics of the bounded fragment of CSL with rewards and
post-processing functions is defined similarly as the semantics of full CSL and
thus we refer the readers to original papers. The key part of the semantics is
given by the definition of the satisfaction relation $\vDash$. It specifies when
a state $s$ satisfies the state formula $\Phi$ (denoted as $s \vDash \Phi$) and
when a path $\omega$ satisfies the path formula $\phi$ (denoted as $\omega
\vDash \phi$). The informal definition of $\vDash$ is as follows:
\begin{itemize}
\item  $s \vDash \mathsf{E}_{\sim r}[\mathsf{I}^{=t}]$ iff $Post(\pi^{\mathcal{C},s,t})$ satisfies $\sim r$.
\item $s \vDash \mathsf{P}_{\sim p}[\phi]$ iff the probability of all paths
  $\omega \in Path^{\mathcal{C}}(s) $ that satisfy the path formula $\phi$
  (denoted as $Prob^{\mathcal{C}}(s,\phi)$) satisfies $\sim p$, where
\begin{itemize}
\item $\omega$ satisfies $\mathsf{X}\mbox{~}\Phi$ iff the second state on
  $\omega$ satisfies $\Phi$
\item $\omega$ satisfies $\Phi\mbox{~}\mathsf{U}^{I}\mbox{~}\Psi$ iff there
  exists time instant $t\in I$ such that the state on $\omega$ occupied at $t$
  satisfies $\Psi$ and all states on $\omega$ occupied before $t'\in [0,t)$
  satisfy $\Phi$
\end{itemize}
\item  $s \vDash \mathsf{R}_{\sim r}[\mathsf{C}^{\leq t}]$ iff the
sum of expected rewards over $Path^{\mathcal{C}}(s)$ \emph{cumulated} until $t$
time units (denoted as $Exp^{\mathcal{C}}(s,\mathsf{X}_{\mathsf{C}^{\leq t}})$) satisfies $\sim r$
\item $s \vDash \mathsf{R}_{\sim r}[\mathsf{I}^{= t}]$ iff the sum of expected rewards
over all paths $\omega \in Path^{\mathcal{C}}(s)$ at time \textit{t} (denoted as
$Exp^{\mathcal{C}}(s,\mathsf{X}_ {\mathsf{I}^{= t}})$) satisfies $\sim\!r$.  
\end{itemize}
A set $Sat_{\mathcal{C}}(\Phi)=\{s \in \mathbb{S} \mid s \vDash \Phi \}$ denotes
the set of states that satisfy $\Phi$.

Note that the syntax and semantics can be easily extended with
``quantitative'' formulae in the form $\Phi ::= \mathsf{P}_{= ?}[\phi] \mid
\mathsf{R}_{= ?}[\mathsf{C}^{\leq t}] \mid \mathsf{R}_{= ?}[\mathsf{I}^{= t}]
\mid \mathsf{E}_{=?}[\mathsf{I}^{=t}]$, i.e., the topmost operator of the
formula~$\Phi$ returns a quantitative result, as used, e.g., in
PRISM~\cite{KNP11}.  In this case the result of a decision procedure is not in
the form of a boolean yes/no answer but the actual numerical value of the
probability $Prob^{\mathcal{C}}(s,\phi)$, the expected reward
$Exp^{\mathcal{C}}(s,\mathsf{X})$ for $\mathsf{X} \in \{ \mathsf{X}_
{\mathsf{I}^{= t}}, \mathsf{X}_ {\mathsf{C}^{\leq t}}\}$ or the value of
$Post^{\mathcal{C}}(s,t)$.  The computation of a numerical value is of the same
complexity class as the computation of a result to be compared leading to a
boolean answer, although in some cases the comparison may be carried out on less
precise or preliminary results. As we will show the quantitative result is much
more suitable for robustness analysis.

To demonstrate that the bounded time fragment of CSL with rewards and
post-processing functions can adequately capture relevant biological behaviours
and thus be successfully used in the robustness analysis of stochastic
biochemical systems, we list several formalisations of such behaviours.
\begin{itemize}
\item \emph{stochastic reachability} -  $\mathsf{P}_{\geq 0.8}[ \mathsf{F}^{[5,10]} ( A \geq 3)]$ expresses the property ``The probability that the population of \textit{A} exceeds 3 between 5 and 10 time units is at least $80\%$''.

\item \emph{stochastic stability} -  $\mathsf{P}_{=?}[ \mathsf{G}^{[0,5]} ( A \geq 1 \wedge A \leq 3 )]$ represents the quantitative property ``What is the probability that the population of \textit{A} remains between 1 and 3 during the first 5 time units?''

\item \emph{stochastic temporal ordering of events} - $\mathsf{P}_{<0.2}[(A \le 2)\ \mathsf{U}^{[2,3]}\ \mathsf{P}_{\geq 0.95} [ ( 2 < A 
\leq 5)\ \mathsf{U}^{[0,10]} (A > 5)]]$ expresses the stochastic version of the following temporal pattern:  ``Species A is initially kept below 2 until it reaches 5 and finally exceeds 5.'' The formula quantifies both the time constrains of the events and the probability that the events occur.  It expresses that ``The probability that the system has following probabilistic temporal pattern is less that $20\%$: the population of \textit{A} is initially kept below 2 until  the system between 2 and 3 times units reaches the states satisfying the subformula $\mathsf{P}_{\geq 0.95} [ ( 2 < A \leq 5)\ \mathsf{U}^{[0,10]} (A > 5)]]$." The subformula specifies the states where  ``The probability that the population of \textit{A} remains greater than 2 and less or equal 5 until it exceeds 5 within 10 time units, is greater than $95\%$." 

\item \emph{cumulative reward property} - $\mathsf{R}_{<5}[\mathsf{C}^{\leq 100}]$, where $\forall s\in \mathbb{S}\ \rho(s) = 1$ if $0 \leq A \leq 3$ in \textit{s}, captures the property that ``The overall time spent in states with population of \textit{A} between 0 and 3 within the first 100 time units, is less than 5 time units'', which can also be understood as ``The probability of the system being in a state with population of \textit{A} between 0 and 3 within the first 100 time units is less then 5\%''.

\item \emph{noise as mean quadratic deviation} - $\mathsf{E}_{<10}[\mathsf{I}^{=100}]$, where the post-processing function is defined as 
$Post(\pi) = \sum_{s \in \mathbb{S}}{\lvert s(A) - mean(\pi,A)\rvert^2 \cdot \pi(s)}$, $s(A)$ gives the population of \textit{A} in state \textit{s} and $mean(\pi,A)$ is the mean of the distribution $\pi$ defined as $mean(\pi,A) = \sum_{s \in \mathbb{S}}{s(A) \cdot \pi(s)}$. This property states that ``The mean quadratic deviation of the distribution of species \textit{A} at time instant $t=100$ must be less then 10''.

\end{itemize}
The $\mathsf{E}$ operator could in principle be extended to allow for intervals and be interpreted as an integral of a user-provided post-processing function over the given time interval. This could lead e.g. to the noise over time interval which is more natural then an instantaneous noise, however the computation complexity of such an operator would be very large. 

\begin{figure} 
\includegraphics[width=\textwidth]{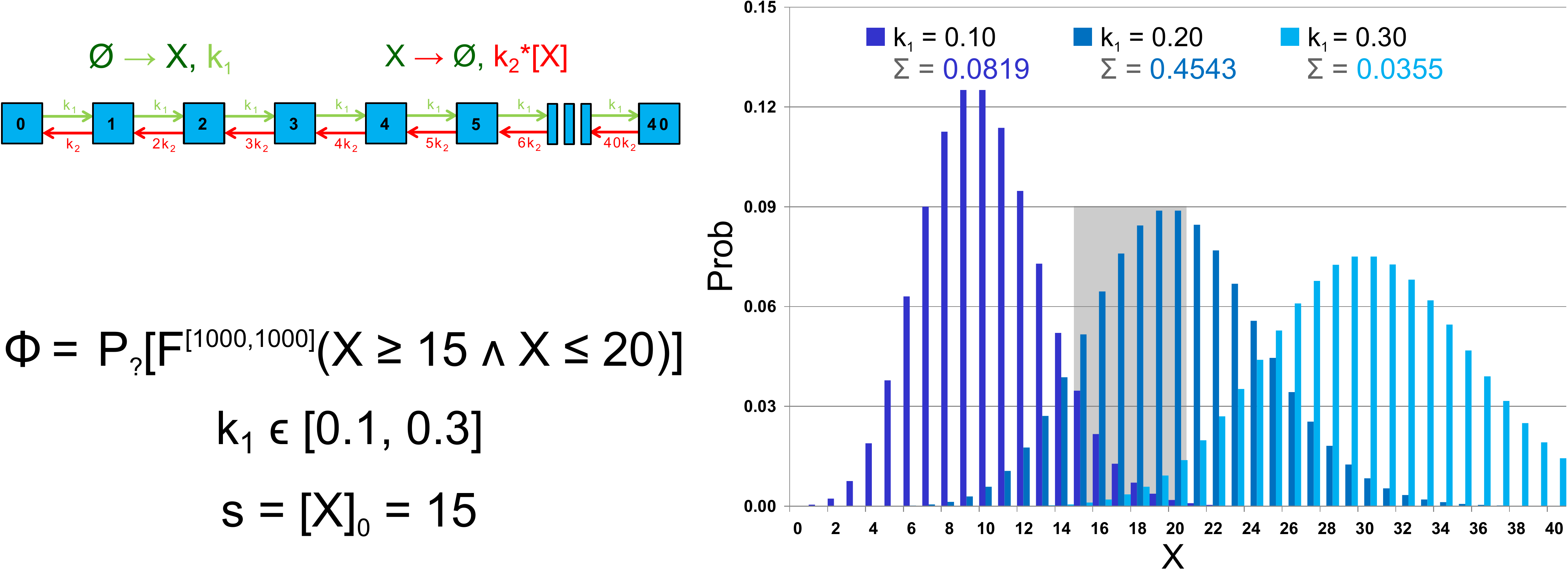}
\caption{{\bf Running example.} 
The example model contains one species \textit{X} with the population bounded to 40, two reactions: production of \textit{X} ($\emptyset \rightarrow X$ with rate $k_1$), degradation of \textit{X} ($X \rightarrow \emptyset$ with rate $k_2 \cdot [X]$, $k_2 = 0.01$) and initial population of \textit{X} is 15.  The corresponding CTMC has 41 states (initial  state $s_0$ corresponds to state with initial population). The inspected formula $\Phi$ represents the quantitative property ``What is the probability that the population of $X$ is between 15 and 20 at time 1000?'' The perturbation space~$\mathbf{P}$ is given by the interval of the rate $k_1 \in [0.1,0.3]$. On the right, there are depicted three transient distributions at time 1000 for three different values of $k_1$ and the resulting probability for the formula $\Phi$ obtained as the sum of probabilities in states with populations from 15 to 20.}
\label{fig:decomposition}
\end{figure}

\subsection{Robustness} 
Let us recap the general definition of Kitano~\cite{Kitano2007} to show how it can be interpreted and how we propose to use it in the context of stochastic systems.

$$
R^{\mathcal{S}}_{\mathcal{A},\mathbf{P}} = \int_{\mathbf{P}}{\psi(p)D^{\mathcal{S}}_{\mathcal{A}}(p)dp} \hspace{1cm} D^{\mathcal{S}}_{\mathcal{A}}(p) = \left\{ 
\begin{array}{cl}
 0             & p \in \mathbf{B} \subset \mathbf{P} \\
 f_{\mathcal{A}}(p)/f_{\mathcal{A}}(0) &  p \in \mathbf{P} \setminus \mathbf{B}
\end{array} \right. 
$$

\subsubsection{Functionality evaluation}

Kitano proposed that the evaluation function $D^{\mathcal{S}}_{\mathcal{A}}(p)$ stating how much the functionality $\mathcal{A}$ is preserved in perturbation \textit{p} should be defined using a subspace $\mathbf{B}$ of all perturbations where the system's function is completely missing and the rest $\mathbf{P} \setminus \mathbf{B}$ where the functions' viability is somehow altered. This definition is meaningful e.g. in cases where the perturbation would lead to a system not having the function at all (speed of reproduction of a dead cell) or in cases where a plain measurement would provide a function's value, however, in reality the system would lack the function altogether (inside temperature during homoeostasis experiment in conditions when an organism loses thermal control and has temperature of environment). These examples have in common that the information about a system lacking its function is provided from outside because if it could be deducible from the system's state alone it could be incorporated into the evaluation function $D^{\mathcal{S}}_{\mathcal{A}}(p)$ itself.

For perturbations $p \in \mathbf{P} \setminus \mathbf{B}$  where the system maintains its function at least partially, Kitano proposes to express the evaluation function $D^{\mathcal{S}}_{\mathcal{A}}(p) = f_{\mathcal{A}}(p)/f_{\mathcal{A}}(0)$ relatively to the ground unperturbed state $f_{\mathcal{A}}(0)$. This is meaningful e.g. for naturally living systems where the ground state is measurable and is considered as an optimal performance state. Such a definition could then enable the comparison of a common property of different species. For example, a reproduction rate for a mouse and a sequoia tree with respect to perturbations of their environment. If a mouse has 20 offsprings per year in base temperature and 22 offsprings for a 2 degrees Kelvin rise then the evaluation function $D^{\mathcal{S}_M}_{\mathcal{A}}(+2 K) = 22/20 = 1.1$. While if a sequoia has 1000 seedlings in ground temperature and 1200 for a 2 degrees Kelvin rise then $D^{\mathcal{S}_S}_{\mathcal{A}}(+2 K) = 1200/1000 = 1.2$. 

We can see that the relativistic nature of Kitano's definition enables comparison of otherwise incomparable organisms and their robustness to perturbations. In our example, the sequoia is more robust to the single perturbation of temperature by $+2K$ than the considered species of mice. However, in cases when no ground state is given the absolute value can be more adequate. The next subsection shows that robustness in stochastic systems can be defined in several different ways providing both the absolute and relative interpretations.

\subsubsection{Robustness in Stochastic systems}

Let $\mathcal{S}$ be a stochastic system with CTMC $\mathcal{C} = \left(\mathbb{S}, s_0, \mathbf{R}, L \right)$, let $\mathbf{P}$ be a space of perturbations to the stochastic kinetic constants of $\mathcal{C}$ and let $\Phi$ be a formula of the bounded time fragment of CSL with rewards and evaluation functions formalising the system's function $\mathcal{A}$. Since the evaluation of $\Phi$ is inherently dependent on the initial conditions of the system that are encoded using the initial state $s_0$, we consider the evaluation function in the form $D^{\mathcal{C},s_0}_{\Phi}$. 

In cases where the set of perturbed stochastic kinetic constants $\mathbf{P}$ is actually extended by initial conditions to $\mathbf{P}^e$, then for a single perturbation point $p = (s_p, k_{1_p}, \ldots , k_{M_p}) \in \mathbf{P}^e$ we consider the initial state $s_0$ of $\mathcal{C}$ to be substituted by $s_p$ in all subsequent expressions, otherwise it remains the original $s_0$.

Let us first define an auxiliary \textit{Eval} function which is then used in the definition of $D^{\mathcal{C},s_0}_{\Phi}$:
\begin{equation}
  Eval_{\Phi}^{\mathcal{C}}( s_0) = \left\{
  \begin{array}{cl}    
    Prob^{\mathcal{C}}(s_0,\phi)                                & \mbox{if~} \Phi \equiv \mathsf{P}_{\star}[\phi] \\[0.3em]
    Exp^{\mathcal{C}}(s_0,\mathsf{X}_{\mathsf{C}^{\leq t}})     & \mbox{if~} \Phi \equiv \mathsf{R}_{\star}[\mathsf{C}^{\leq t}] \\[0.3em]
    Exp^{\mathcal{C}}(s_0,\mathsf{X}_{\mathsf{I}^{= t}})        & \mbox{if~} \Phi \equiv \mathsf{R}_{\star}[\mathsf{I}^{= t}] \\[0.3em]
    Post(\pi^{\mathcal{C},s_0,t}) & \mbox{if~} \Phi \equiv \mathsf{E}_{\star}[\mathsf{I}^{= t}]
  \end{array} \right. 
  \label{eg:eval}
\end{equation}
where $\star \in \{=?, \sim\! r\}$. Given these specifications the evaluation function $D^{\mathcal{C},s_0}_{\Phi}$ can be restated in several different ways: 
\begin{subequations}
\begin{align}
  D^{\mathcal{C},s_0}_{\Phi}(p) &= \left\{ 
  \begin{array}{clcl}
    \hspace{1cm} 0 \hspace{1cm} & p \in \mathbf{B} \subset  \mathbf{P}   & \vee   &  Eval_{\Phi}^{\mathcal{C}_p}(s_0)  \nsim r \\[0.5mm]
    1                           &  p \in \mathbf{P} \setminus \mathbf{B} & \wedge &  Eval_{\Phi}^{\mathcal{C}_p}(s_0)  \sim r
  \end{array} \right.   
  \label{eq:sem_boolean}\\  
  D^{\mathcal{C},s_0}_{\Phi}(p) &= \left\{ 
  \begin{array}{cl}
    \hspace{1cm} 0 \hspace{1cm}                  & p \in \mathbf{B} \subset  \mathbf{P} \\[0.8mm]
     \frac{Eval_{\Phi}^{\mathcal{C}_p}(s_0)}{r}  & \mbox{else if~} \sim \in \{ \geq, > \} \\[0.8mm]
     \frac{r}{ Eval_{\Phi}^{\mathcal{C}_p}(s_0)} & \mbox{else if~} \sim \in \{ \leq, < \}
  \end{array} \right. 
  \label{eq:sem_relative}\\
    D^{\mathcal{C},s_0}_{\Phi}(p) &= \left\{ 
  \begin{array}{cl}
    \hspace{1cm} 0 \hspace{1cm}         & p \in \mathbf{B} \subset  \mathbf{P} \\[0.5mm]
     Eval_{\Phi}^{\mathcal{C}_p}(s_0)   & \mbox{else}
  \end{array} \right. 
  \label{eq:sem_absolute}\\
  D^{\mathcal{C},s_0}_{\Phi}(p) &= \left\{ 
   \begin{array}{cl}
    \hspace{1cm} 0 \hspace{1cm}                            & p \in \mathbf{B} \subset  \mathbf{P} \\[0.8mm]
    \lvert Eval_{\Phi}^{\mathcal{C}_p}(s_0) - X \rvert^2   & \mbox{else, ~} X = agr \{Eval_{\Phi}^{\mathcal{C}_p}(s_0) \mid \mathcal{C}_p \in \mathbf{P} \} \wedge agr\in \{min,max,avg \}
  \end{array} \right. 
  \label{eq:sem_relative2} 
\end{align}
\end{subequations} 

The first definition of the evaluation function~(\ref{eq:sem_boolean}) is possible for the specification where the topmost operator of the formula $\Phi$ includes the threshold $r$ (i.e. $\star = \sim \! r$). Because $D^{\mathcal{C},s_0}_{\Phi}(p)$ returns a qualitative result robustness  $R^{\mathcal{C}}_{\Phi,\mathbf{P}}$ specifies the measure of all perturbations in $\mathbf{P}$ for which the property holds in a strictly boolean sense -- it is the fraction of $\mathbf{P}$ where the property is valid. This definition can be used, e.g., in the property $\Phi_A = \mathsf{P}_{\geq 0.8}[ \mathsf{F}^{[0,5]} (X > 300)]$ which specifies that in $80\%$ of cases  the population of X is larger than 300 within 5 seconds. For this property and a model with a parameter $k \in [0,10]$ the robustness gives us the fraction of the parametric interval $[0,10]$ for which the model satisfies $\Phi_A$. 

In the second definition~(\ref{eq:sem_relative}) $D^{\mathcal{C},s_0}_{\Phi}(p)$ returns the quantitative value that is relative to the threshold~\textit{r}. Therefore, robustness can be interpreted as the average relative validity of the property over $\mathbf{P}$. If \textit{r}~corresponds to the validity of $\Phi$ in conditions considered natural for the inspected system $\mathcal{S}$ (i.e, to the unperturbed state) then this interpretation complies with the original definition of Kitano. Let us consider the same property $\Phi_A$ and the same parametric space $k \in [0,10]$. If in $60\%$ of model behaviours the population of X is larger than 300 within 5 seconds than the robustness is 0.6/0.8 = 0.75. If the probability is different in each \textit{k} then the robustness gives us the average value that meets our expectations.

The third definition~(\ref{eq:sem_absolute}) is possible for specifications using the quantitative semantics of formula $\Phi$ (i.e.~$\star = ?$). The robustness gives the mean validity over all $\mathbf{P}$ regardless of any probability threshold~\textit{r}. This interpretation is convenient when there are no \textit{a priori} assumptions about the system expected behaviour.

Finally, to express the fact that the system behaviour remains the same (with respect to the evaluation function) across the space of perturbations we introduce the fourth definition~(\ref{eq:sem_relative2}). It uses an aggregation function to compute a mean value and then express the variance from the mean.
This definition enables us to compare models which have same numerical values of robustness in the sense of definition~(\ref{eq:sem_absolute}) but which achieve the average value with very different landscapes of evaluation function.

While the last three definitions require the precise computation of the probability value in every $p \in \mathbf{P}$, the first definition is amenable to approximate solutions. In this case it suffices to ensure that the probability is larger or smaller then \textit{r}. In many cases it can be achieved without computing the precise value and thus statistical model checking techniques can be efficiently used. In both case studies we use definition~(\ref{eq:sem_absolute}), since we do not consider any ground unperturbed state. We assume $\mathbf{B}$ to be an empty set and expect all the lack of functionality $\mathcal{A}$ to be fully expressible in terms of the property $\Phi$.

\subsection{Robustness computation} 

Now we look how robustness $R^{\mathcal{C}}_{\Phi,\mathbf{P}^e}$ can be efficiently computed by using the evaluation function $D^{\mathcal{C},s_0}_{\Phi}$. Let us first consider the case where the space of perturbations $\mathbf{P}$ does not contain different initial states.

As will be shown in the next section the computation of $Eval^{\mathcal{C}_p}_{\Phi}(s_0)$ even for a single perturbation point $p$ is rather complex, therefore a computation of the integral over the whole space of perturbations is not possible in an explicit sense. Instead a way to approximate the upper and lower bounds $R_{\Phi,\mathbf{P},\top}^{\mathcal{C}}$ and $R_{\Phi,\mathbf{P},\bot}^{\mathcal{C}}$ is introduced enabling the approximation of the value of the integral as
$$
\begin{array}{rcl}
R_{\Phi,\mathbf{P}}^{\mathcal{C}} &\stackrel{def}{=}& \displaystyle\int_{\mathbf{P}}{\psi(p)D^{\mathcal{C}}_{\Phi}(p)dp} \\[0.6em]
R_{\Phi,\mathbf{P}}^{\mathcal{C}} &\simeq& \displaystyle\frac{1}{2}\left( 
  R_{\Phi,\mathbf{P},\top}^{\mathcal{C}} + R_{\Phi,\mathbf{P},\bot}^{\mathcal{C}} \right) \pm Err_{\Phi,\mathbf{P}}^{\mathcal{C}} \hspace{1cm} 
Err_{\Phi,\mathbf{P}}^{\mathcal{C}} = \displaystyle\frac{1}{2}\left( R_{\Phi,\mathbf{P},\top}^{\mathcal{C}} - R_{\Phi,\mathbf{P},\bot}^{\mathcal{C}} \right)	
\end{array}
$$

The computation of $R_{\Phi,\mathbf{P},\top}^{\mathcal{C}}$ and $R_{\Phi,\mathbf{P},\bot}^{\mathcal{C}}$ is due to the approximation of the upper $D^{C}_{\Phi,\mathbf{P},\top}$ and lower $D^{\mathcal{C}}_{\Phi,\mathbf{P},\bot}$ bounds for values of the evaluation function $D^{\mathcal{C}}_{\Phi}(p)$ over $\mathbf{P}$
$$
D^{\mathcal{C}}_{\Phi,\mathbf{P},\top} \geq max \left\{ D^{\mathcal{C}}_{\Phi}(p) \mid p \in \mathbf{P} \right\}\hspace{1cm} D^{\mathcal{C}}_{\Phi,\mathbf{P},\bot} \leq min \left\{ D^{\mathcal{C}}_{\Phi}(p) \mid p \in \mathbf{P} \right\}
$$ 
Because such an approximation would be too course for most cases a finite decomposition of the perturbation space $\mathbf{P}$  into perturbation subspaces $\mathbf{P} = \mathbf{P}_1 \cup \ldots \cup \mathbf{P}_n$ is used which then under the assumption of equal probability of all perturbations gives better robustness bounds. Hence we get that:
\begin{equation}
R_{\Phi,\mathbf{P},\top}^{\mathcal{C}} = \sum_{i = 1}^{n}{\frac{|\mathbf{P}_i|}{|\mathbf{P}|} \cdot D_{\Phi,\mathbf{P}_i,\top}^{\mathcal{C}}} \hspace{1cm}
R_{\Phi,\mathbf{P},\bot}^{\mathcal{C}} = \sum_{i = 1}^{n}{\frac{|\mathbf{P}_i|}{|\mathbf{P}|} \cdot D_{\Phi,\mathbf{P}_i,\bot}^{\mathcal{C}}} 
\label{eq:robustness_sum}
\end{equation}

Let us now consider the case in which the space of perturbations is extended with initial states $\mathbf{P}^e = \mathbb{I} \times \mathbf{P}$ where $\mathbb{I} \subseteq \mathbb{S}$ and $\mathbf{P}$ is non-singular, for this case the integral defining robustness is actually a finite sum of integrals:
$$
R_{\Phi,\mathbf{P}^e}^{\mathcal{C}} \stackrel{def}{=} \displaystyle\sum_{s \in \mathbb{I}}{\frac{1}{|\mathbb{I}|} \displaystyle\int_{p \in \mathbf{P}}{\psi(p)D^{\mathcal{C}}_{\Phi}(p)dp}} = \frac{1}{|\mathbb{I}|} \displaystyle\sum_{s \in \mathbb{I}}{ R_{\Phi,\mathbf{P}}^{\mathcal{C}} }
$$
where $\psi(p)$ gives the probability of perturbation \textit{p} with respect to $\mathbf{P}$. This expression is valid for uniform distributions of the initial states over the whole space of perturbations $\mathbf{P}^e$, however, it can be straightforwardly modified for non-uniform distributions. Using the expression the robustness computation for perturbations containing a single initial state can be easily extended to perturbations containing different initial states. Moreover, in Section~\ref{sec:GMC}, we show that  for most properties the model checking procedure (utilised in the robustness computation) returns results for an arbitrary set of initial states $\mathbb{I} \subseteq \mathbb{S}$ with the same time complexity as for a single state.

\begin{figure} 
\begin{center}
\includegraphics[width=\textwidth]{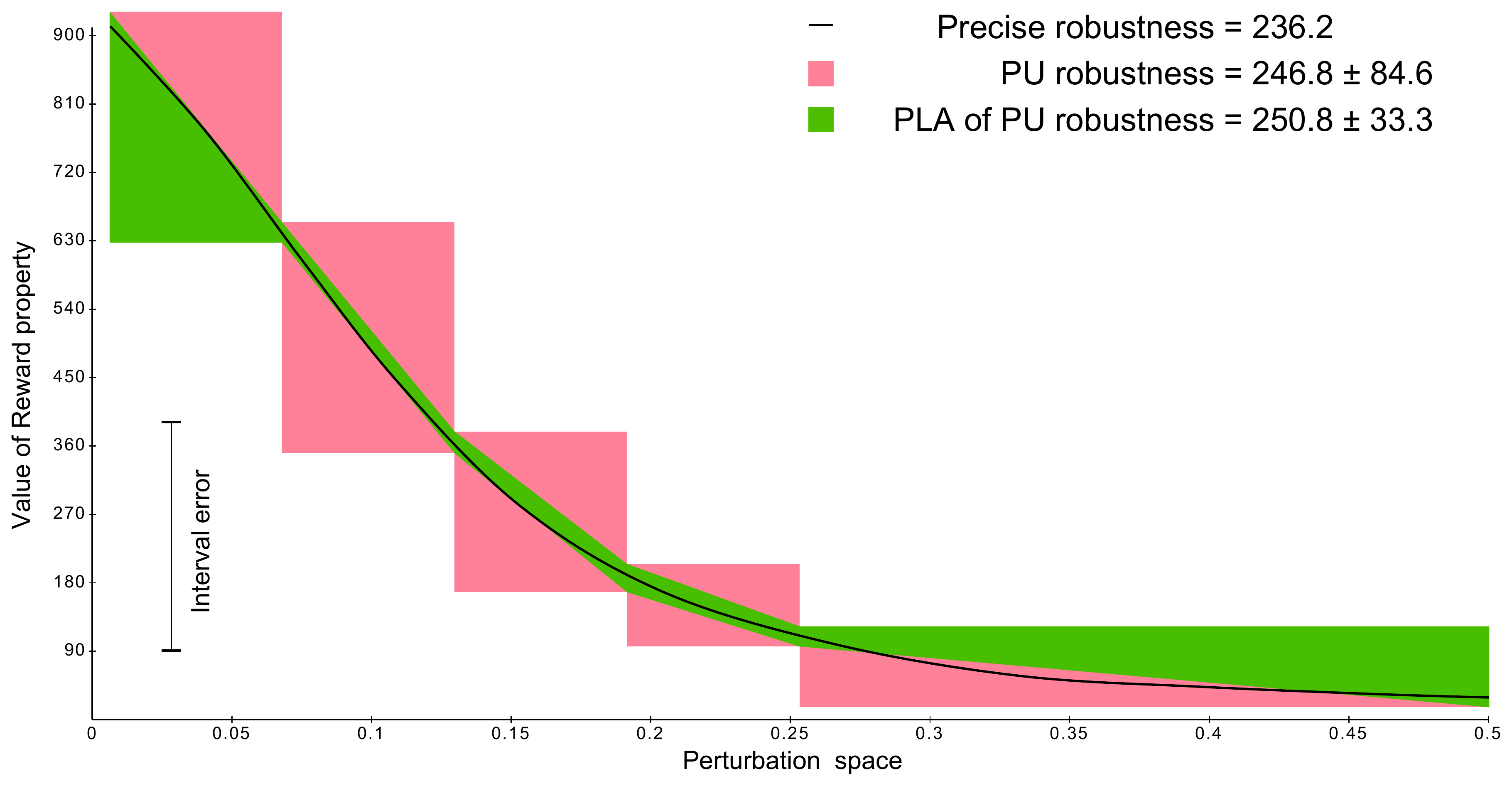}
\end{center}
\caption{{\bf Piecewise linear approximation of robustness.} 
An improved approximation is shown in dark green, it is computed by linearly interpolating grid points in which the upper and lower bounds of a property may be computed more precisely as the minimum resp. maximum of the values from all parameter subintervals sharing boundary grid points. The obtained result is more precise that the original robustness (in light pink) albeit without the conservative guarantee on bounds.}
\label{fig:pla_approx}
\end{figure}

The accuracy of the approximation can be further improved using the \textit{piece-wise linear approximation} of robustness. This concept is illustrated in Figure~\ref{fig:pla_approx}. Since the spaces $\mathbf{P}_i$ and $\mathbf{P}_{i+1}$ have a common point \textit{p} (in a general \textit{n} dimensional perturbation space $2^n$ subspaces intersect in a single point \textit{p}), we can use this to obtain a more precise range of values for the value of the property $\Phi$ in \textit{p} as 
$$D_{\Phi,p,\top}^{\mathcal{C}} = min\left\{ D_{\Phi,\mathbf{P}_i,\top}^{\mathcal{C}} \mid p \in \mathbf{P}_i \right\} \text{ and } D_{\Phi,p,\bot}^{\mathcal{C}} = max\left\{ D_{\Phi,\mathbf{P}_i,\bot}^{\mathcal{C}} \mid p \in \mathbf{P}_i \right\}.$$ 

Under the assumption that the value of a property does not change rapidly over  sufficiently small subspaces $\mathbf{P}_i$ the resulting upper and lower bound of robustness can then be computed from linear interpolation of grid points \textit{p}. The decision in which cases such an assumption is acceptable is up to user since there is in general no efficient way of resolving this situation. In such a case the overall piecewise linear approximation of robustness will usually have a higher precision albeit without the guarantee of strict upper and lower bounds.

To understand how $D_{\Phi,\mathbf{P}_i,\top}^{\mathcal{C}}$ and $D_{\Phi,\mathbf{P}_i,\bot}^{\mathcal{C}}$ can be efficiently computed first the methods for transient analysis and global CSL model checking based on \emph{uniformisation} are revisited~\cite{Baier2000mcviatran,Kwiatkowska2007}. Afterwards we present the \emph{min-max approximation}~\cite{CAV2013} that allows us to approximate the quantitative model checking result for continuous sets of parameterised CTMCs. The key idea is to employ a method called \emph{parameterised uniformization} -- a modification of the standard uniformization technique presented in~\cite{CAV2013}. Finally, we show how to control the \emph{approximation error} in order to obtain the required error bound.

\subsubsection{Transient analysis}
The aim of transient analysis is to compute a transient probability distribution. Given an initial distribution $\pi^{\mathcal{C}, s_0, 0}$ (i.e. $\pi^{\mathcal{C}, s_0, 0}(s) = 1$ if $s_0 = s$, and 0, otherwise) at time 0 of a CTMC $\mathcal{C} = \left( \mathbb{S}, s_0, \mathbf{R}\right)$ what will the transient state distribution $\pi^{\mathcal{C},s_0,t}$ look like in some future yet finite time $t \in \mathbb{R}_{\geq 0}$. 
 
Transient analysis of a CTMC may be efficiently carried out by a standard technique called \emph{uniformization}~\cite{Kwiatkowska2007}. The transient probability in time \textit{t} is obtained as a sum of expressions giving the state distributions after \textit{i} discrete reaction steps of the respective \emph{uniformized} discrete time Markov chain (DTMC) weighted by the \textit{i}th probability of the Poisson process. It is the probability of \textit{i} such steps occurring in time~\textit{t}, assuming the delays between steps of the CTMC $\mathcal{C}$ are exponentially distributed with \emph{rate} $q$. Formally, for the rate $q$ satisfying $q \geq max\{E^{\mathcal{C}}(s) \mid s\in \mathbb{S}\}$ (\textit{E} is the exit rate of state \textit{s}) the uniformized DTMC $\mathsf{unif}(\mathcal{C})$ is defined as $\mathsf{unif}(\mathcal{C}) = \left( \mathbb{S}, s_0, \mathbf{Q}^{\mathsf{unif}(\mathcal{C})} \right)$ where 
$$ \mathbf{Q}^{\mathsf{unif}(\mathcal{C})}(s,s') = \left\{
\begin{array}{cl}
    \frac{\mathbf{R}(s,s')}{q}                       & \mbox{if~~} s \neq s' \\
    1-\sum_{s'' \neq s}{\frac{\mathbf{R}(s,s'')}{q}}   & \mbox{otherwise.}
\end{array} \right.
$$
and the $i$th Poisson probability in time $t$ is given as $\gamma_{i,q \cdot t}= e^{-q \cdot t} \cdot \frac{\left(q\cdot t\right)^i}{i!}$. The transient probability can be computed as follows:
$$\pi^{\mathcal{C}, s_0, t} = \sum^{\infty}_{i=0}{\gamma_{i,q
    \cdot t} \cdot \pi^{\mathcal{C},s_0,0} \cdot (
  \mathbf{Q}^{\mathsf{unif}(\mathcal{C})})^i} \approx
\sum^{R_{\epsilon}}_{i=L_{\epsilon}}{\gamma_{i,q \cdot t} \cdot
  \pi^{\mathcal{C},s_0,0} \cdot ( \mathbf{Q}^{\mathsf{unif}(\mathcal{C})})^i}.$$
Although the sum is in
general infinite, for a given precision~$\epsilon$ the lower and upper bounds
$L_{\epsilon}, R_{\epsilon}$ can be estimated by using techniques such as of Fox
and Glynn~\cite{FoxGlynn1988} which also allow for efficient solutions of the
Poisson process. In order to make the computation of uniformization 
feasible the matrix-matrix multiplication is reduced to a vector-matrix 
multiplication, i.e., $$ \pi^{\mathcal{C}, s_0, 0} \cdot
(\mathbf{Q}^{\mathsf{unif}(\mathcal{C})})^i = ( \pi^{\mathcal{C}, s_0, 0} \cdot
(\mathbf{Q}^{\mathsf{unif}(\mathcal{C})} )^{i-1} ) \cdot
\mathbf{Q}^{\mathsf{unif}(\mathcal{C})}.$$

Standard uniformization can be intractable when the system under study is too complex, i.e., contains more than in order of
$10^7$ states and the upper estimate $R_{\epsilon}$, denoting the number of vector-matrix multiplications as iterations, is high (more than in order of $10^6$). Therefore, many approximation techniques have been studied in order to reduce the state space and to lower the number of iterations $R_{\epsilon}$.  State space reductions are based on the observation that in many cases (especially in biochemical systems) a significant amount of the probability mass in a given time is localized in a manageable set of states. Thus neglecting states with insignificant probability can dramatically reduce the state space while the resulting approximation of the transient probability is still sufficient. Methods allowing efficiently state-space reduction are based on finite projection techniques~\cite{munsky2006finite,Henzinger2009} and dynamic state space truncation~\cite{Didier}.

Since the number of iterations $R_{\epsilon}$ inherently depends on the uniformization rate \textit{q} that has to be greater then the maximal exit rate of all the states of the system, a variant of standard uniformization, so-called \emph{adaptive} uniformization~\cite{Moorsel94}, has been proposed. It uses a uniformization rate that adapts depending on the set of states the system can occupy at a given time, i.e, after a particular number of reactions. In many cases, a significantly smaller rate \textit{q} can be used and thus the number of iterations $R_{\epsilon}$ can be significantly reduced during some parts of the computation. Moreover, adaptive uniformization can be successfully combined with reduction techniques mentioned above~\cite{Didier}. The downside of adaptive uniformization is that the Poisson process has to be replaced with a general \emph{birth} process which is more expensive to solve. See, e.g~\cite{Moorsel94}, for more details. 

For sake of simplicity, we present our methods for the computation of $D_{\Phi,\mathbf{P}_i,\top}^{\mathcal{C}}$ and $D_{\Phi,\mathbf{P}_i,\bot}^{\mathcal{C}}$ using standard uniformization. However, our method can be successfully combined with the aforementioned techniques.

\subsubsection{Global CSL Model Checking}
\label{sec:GMC}

The aim of the global model checking technique is to efficiently compute for any CSL formula $\Phi$ the values  $Eval_{\Phi}^{\mathcal{C}}(s)$ for all states $s \in \mathbb{S}$. On the other hand, the goal of local model checking technique is to compute $Eval_{\Phi}^{\mathcal{C}}(s)$ for a single state $s \in \mathbb{S}$. The crucial advantage of the global approach is the fact that it has the same asymptotic and also practical complexity as the local approach. Therefore, the global model checking technique is much more suitable for robustness analysis over perturbations of initial conditions that are encoded as the initial state of the corresponding CTMC. 

Global model checking returns the vector of size $|\mathbb{S}|$ such that the \textit{i}th position contains the model checking result provided that $s_i$ is the initial state. Let $\mathcal{C} =\left(\mathbb{S}, \mathbf{R}, L \right)$ be a labelled CTMC where the initial state is not specified. The crucial part of this method is to compute the vector of probabilities $\overline{Prob}^{\mathcal{C}, \phi}$ for any path formula $\phi$ and the vector of  expected rewards $\overline{Exp}^{\mathcal{C}, \mathsf{X}}$ for $\mathsf{X} \in \{ \mathsf{X}_ {\mathsf{I}^{= t}}, \mathsf{X}_ {\mathsf{C}^{\leq t}}\}$ such that for all $s\in \mathbb{S}$ the following holds:
$$ \overline{Prob}^{\mathcal{C}, \phi}(s) = Prob^{\mathcal{C}}(s,\phi) \wedge \overline{Exp}^{\mathcal{C}, \mathsf{X}}(s) = Exp^{\mathcal{C}}(s,\mathsf{X}) $$ 

In local model checking the computation of $Prob^{\mathcal{C}}(s,\phi)$ and $Exp^{\mathcal{C}}(s,\mathsf{X})$
is reduced to the computation of the transient probability distribution $\pi^{\mathcal{C},s,t}$, see~\cite{Baier2000mcviatran,Kwiatkowska2007} for more details. Thus, for different initial states $s$ we have to compute the corresponding transient probability distributions separately. The key idea of the global model checking method is to use \emph{backward transient analysis}. The result of backward transient analysis is the vector $\tau^{\mathcal{C},\mathbb{A},t}$ such that for arbitrary set of states~$\mathbb{A}$, the value $\tau^{\mathcal{C},\mathbb{A},t}(s)$ is the probability that  $\mathbb{A}$ is reached from $s$ at the time \textit{t}. Without going into details the vector $\tau^{\mathcal{C},\mathbb{A},t}$ can be computed in a very similar way using the uniformized DTMC $\mathsf{unif}(\mathcal{C})$ as in the case of vector $\pi^{\mathcal{C},s,t}$. Only vector-matrix multiplications is replaced by matrix-transposed-vector multiplication and $\tau^{\mathcal{C},\mathbb{A},0}(s) = 1$ if $s\in\mathbb{A}$, and $0$, otherwise.

The global model checking technique can not be used if $\Phi$ includes the operator $\mathsf{E}_{\sim r}[I^{=t}]$. In such a case we have to compute the value $Post(\pi^{\mathcal{C},s,t})$. Hence the local model checking technique has to be employed, i.e., we first compute the vector $\pi^{\mathcal{C},s,t}$ and then apply the user specified function \textit{Post}.

Now we briefly show how the vector $\overline{Prob}^{\mathcal{C}, \phi}$ is computed using backward transient analysis. Since the definition of next operator $\mathsf{X}\ \Phi$ does not rely on any real time aspects of CTMCs, its evaluation stems from the probability of the next reaction that can be easily obtained from the transition matrix $\mathbf{R}$. 
The evaluation of the until operator $\Phi_1 \mathsf{U}^{I} \Phi_2$ depends on the form of the interval $I$ and is separately solved for the cases of $I=[0,t_1]$ and $I=[t_1, t_2]$ where $t_1, t_2 \in \mathbb{R}_{\geq 0}$. It is based on a modification of the uniformized infinitesimal generator matrix $\mathbf{Q}^{\mathsf{unif}}$ where certain states are made absorbing. This means that all outgoing transitions are ignored in dependence on the validity of $\Phi_1$ and $\Phi_2$ in these states. 

For any CSL formula $\Phi$, let $\mathcal{C}[\Phi] = \left( \mathbb{S}, s_0, \mathbf{R}[\Phi], L \right)$, where $\mathbf{R}[\Phi](s,s') = \mathbf{R}(s,s')$, if $s \vDash \Phi$, and $0$, otherwise. The formula $\phi = \Phi_1 \ \mathsf{U}^{[0,t]} \ \Phi_2$ can be evaluated using the vector $\tau^{\mathcal{C},\mathbb{A},t}$ in the following way:
$$
\overline{Prob}^{\mathcal{C}, \phi}  =  \tau^{\mathcal{C}[\neg \Phi_1 \wedge \Phi_2], \mathbb{A}, t} \mbox{~where~} s \in \mathbb{A} \mbox{~iff~} s \vDash \Phi_2.
$$
For the formula $\phi = \Phi_1 \ \mathsf{U}^{[t_1,t_2]} \ \Phi_2$ the evaluation is split into two parts: staying in states satisfying $\Phi_1$ until time $t_1$ and reaching a state satisfying $\Phi_2$, while remaining in states satisfying $\Phi_1$, within time $t_2-t_1$. The formula $\phi$ can  be evaluated using the vector  $\tau^{\mathcal{C},\overline{v},t}$ that takes a vector $\overline{v}$ instead of a set $\mathbb{A}$ (i.e., $\tau^{\mathcal{C}, \overline{v}, 0} = \overline{v}$) in the following way:    
$$\overline{Prob}^{\mathcal{C}, \phi}  = \tau^{\mathcal{C}[\neg \Phi_1], \overline{v}, t_1} \mbox{~where~} \overline{v} =   \tau^{\mathcal{C}[\neg \Phi_1 \wedge \Phi_2], \mathbb{A}, t_2-t_1} \mbox{~and~} s \in \mathbb{A} \mbox{~iff~} s \vDash \Phi_2.
$$

The backward transient analysis can be also used in the case of reward computation. Since operator $\mathsf{R}_{\sim p}[\mathsf{I}^{= t}]$ expresses the expected reward at time $t$, the vector $\overline{Exp}^{\mathcal{C}, \mathsf{X}_{\mathsf{I}^{= t}}}$  can be computed as follows:
$$ \overline{Exp}^{\mathcal{C}, \mathsf{X}_{\mathsf{I}^{= t}}} =  \tau^{\mathcal{C}, \overline{v}, t} \mbox{~where~} \overline{v} = \rho \mbox{~such that~} \rho \mbox{~is the given state reward structure.}$$
For evaluation of the operator $\mathsf{R}_{\sim p}[\mathsf{C}^{\leq t}]$ we have to use \emph{mixed Poisson probabilities} (see, e.g., \cite{Kwiatkowska2006rewards,Kwiatkowska2007}) in the backward transient analysis. It means that during the uniformization the Poisson probabilities $\gamma_{i,q \cdot t}$ are replaced by the mixed Poisson probabilities $\bar{\gamma}_{i,q \cdot t}$ that can be computed as: 
$$ 
\bar{\gamma}_{i,q \cdot t} = \frac{1}{q} \cdot \left( 1-  \sum_{j = 1}^{i}  \gamma_{j,q \cdot t} \right) \mbox{.}
$$
Using the given state reward structure $\rho$ we can compute the vector $\overline{Exp}^{\mathcal{C}, \mathsf{X}_ {\mathsf{C}^{\leq t}}}$ in the following way:

$$ \overline{Exp}^{\mathcal{C}, \mathsf{X}_ {\mathsf{C}^{\leq t}}} =  \tau^{\mathcal{C}, \overline{v}, t} \mbox{~where~} \overline{v} = \rho \mbox{~and the mixed Poisson probabilities~} \bar{\gamma}_{i,q \cdot t} \mbox{~are used.}$$

To recap the overall method of stochastic model checking of CTMCs over CSL formulae we present the methods from an abstract perspective. The evaluation of a structured formula $\Phi$ proceeds by bottom-up evaluation of a set of atomic propositions, probabilistic or expected reward inequalities and their boolean combinations. This evaluation gives us a discrete set of states that are further used in the following computation. The process continues up the formula until the root is reached. The final verdict is reported either in the form of a boolean yes/no answer or as the actual numerical value of the probability or the expected reward. This process can be easily extended for the operator  $\mathsf{E}_{\sim r}[I^{=t}]$, however, the local model checking method has to be used.

\subsubsection{Min-max approximation}

The key idea of min-max approximation is to approximate the \emph{largest set of states satisfying}~$\Phi$, and the \emph{smallest set of states satisfying} $\Phi$ with respect to the space of perturbations $\mathbf{P}$. Let $\mathbf{C}$ be a set of parameterised CTMCs induced by the space of perturbations $\mathbf{P}$ in the system $\mathcal{S}$. We compute the approximation $Sat^{\top}_{\mathbf{C}}(\Phi)$ and $Sat^{\bot}_{\mathbf{C}}(\Phi)$ such that
$$Sat^{\top}_{\mathbf{C}}(\Phi) \supseteq \bigcup_{\mathcal{C}_p \in \mathbf{C}} Sat_{\mathcal{C}_p}(\Phi) \ \wedge \  Sat^{\bot}_{\mathbf{C}}(\Phi) \subseteq \bigcap_{\mathcal{C}_p \in \mathbf{C}} Sat_{\mathcal{C}_p}(\Phi)$$ 
where $s\in Sat_{\mathcal{C}_p}(\Phi)$ iff $s\vDash \Phi$ in CTMC $\mathcal{C}_p$.
To obtain such approximations we extended the satisfaction relation $\vDash$ and showed that it is sufficient for an arbitrary path formula $\phi$, and $\mathsf{X} \in
\{\mathsf{X}_ {\mathsf{C}^{\leq t}}, \mathsf{X}_ {\mathsf{I}^{= t}}\}$
to compute the vectors $\overline{Prob}^{\mathbf{C}, \phi}_{\top},\ \overline{Prob}^{\mathbf{C}, \phi}_{\bot},\  \overline{Exp}^{\mathbf{C}, \mathsf{X}}_{\top}$ and $\overline{Exp}^{\mathbf{C}, \mathsf{X}}_{\bot}$ such that for each $s\in \mathbb{S}$ the following holds:
\begin{equation}
\label{eq:glob}
\begin{array}{rl}
\overline{Prob}^{\mathbf{C}, \phi}_{\top}(s) & \geq max\{ \overline{Prob}^{\mathcal{C}_p, \phi}(s) \mid \mathcal{C}_p \in \mathbf{C} \}  \\ 
\overline{Prob}^{\mathbf{C}, \phi}_{\bot}(s)  &  \leq  min \{  \overline{Prob}^{\mathcal{C}_p, \phi}(s) \mid \mathcal{C}_p \in \mathbf{C} \}  \\
\overline{Exp}^{\mathbf{C}, \mathsf{X}}_{\top}(s) & \geq   max\{ \overline{Exp}^{\mathcal{C}_p, \mathsf{X}}(s) \mid \mathcal{C}_p \in \mathbf{C} \}  \mbox{~for~} \mathsf{X} \in \{ \mathsf{X}_ {\mathsf{I}^{= t}}, \mathsf{X}_ {\mathsf{C}^{\leq t}} \}\\
\overline{Exp}^{\mathbf{C}, \mathsf{X}}_{\bot}(s) & \leq   min\{ \overline{Exp}^{\mathcal{C}_p, \mathsf{X}}(s) \mid \mathcal{C}_p \in \mathbf{C} \}  \mbox{~for~}\mathsf{X} \in \{ \mathsf{X}_ {\mathsf{I}^{= t}}, \mathsf{X}_ {\mathsf{C}^{\leq t}} \}.\\
\end{array}
\end{equation}

The min-max approximation can be easily extended to the operator $\mathsf{E}_{\sim r}[I^{=t}]$. For the given state $s \in \mathbb{S}$ and the time \textit{t} it is sufficient to compute the values $Post^{\mathbf{C}}_{\top}(s,t)$ and $Post^{\mathbf{C}}_{\top}(s,t)$ such that the following holds:
\begin{equation}
\label{eq:loc}
\begin{array}{rl}
Post^{\mathbf{C}}_{\top}(s,t) & \geq  max\{ Post(\pi^{\mathcal{C}_p,s,t}) \mid \mathcal{C}_p \in \mathbf{C} \}  \\ 
Post^{\mathbf{C}}_{\bot}(s,t) & \leq  min\{ Post(\pi^{\mathcal{C}_p,s,t}) \mid \mathcal{C}_p \in \mathbf{C} \}.  \\
\end{array}
\end{equation}

The approximated sets ${Sat}^{\top}_{\mathbf{C}}(\Phi)$ and
${Sat}^{\bot}_{\mathbf{C}}(\Phi)$ are further used in the computation of
$D_{\Phi,\mathbf{P},\top}^{\mathcal{C},s}$ and $D_{\Phi,\mathbf{P},\bot}^{\mathcal{C},s}$. If the topmost operator of the formula $\Phi$
is $\mathsf{P}_{= ?}[\phi]$ then 
$$ D_{\Phi,\mathbf{P},\bot}^{\mathcal{C},s} =
\overline{Prob}^{\mathbf{C}, \phi}_{\bot}(s) \wedge D_{\Phi,\mathbf{P},\top}^{\mathcal{C}, s} =
\overline{Prob}^{\mathbf{C}, \phi}_{\top}(s). $$
If the topmost operator of the formula $\Phi$ is $\mathsf{R}_{= ?}[\mathsf{C}^{\leq t}]$ and $\mathsf{R}_{= ?}[\mathsf{I}^{= t}]$ then
$$D_{\Phi,\mathbf{P},\bot}^{\mathcal{C},s} =
\overline{Exp}^{\mathbf{C}, \mathsf{X}}_{\bot}(s) \wedge D_{\Phi,\mathbf{P},\top}^{\mathcal{C},s} =
\overline{Exp}^{\mathbf{C}, \mathsf{X}}_{\top}(s) \mbox{~for~} \mathsf{X} =\mathsf{X}_
{\mathsf{C}^{\leq t}} \mbox{~and~} \mathsf{X} = \mathsf{X}_ {\mathsf{I}^{= t}} \mbox{, respectively}.$$
Similarly, if the topmost operator of the formula $\Phi$ is $\mathsf{E}_{= ?}[\mathsf{I}^{= t}]$ then
$$D_{\Phi,\mathbf{P},\bot}^{\mathcal{C}, s} = Post^{\mathbf{C}}_{\bot}(s,t)
 \wedge D_{\Phi,\mathbf{P},\top}^{\mathcal{C}, s} =
Post^{\mathbf{C}}_{\top}(s,t).$$

\subsubsection{Parameterised uniformisation}
\label{sec:parameteriseduniformisation}

Recall that the most crucial part of the robustness computation is given by the fact that the space of perturbations of stochastic rate constants $\mathbf{P}$ is dense and thus the set $\mathbf{C}$ is infinite. Therefore, it is not possible to employ the standard model checking techniques to compute the result for each CTMC $\mathcal{C}_p \in$~$\mathbf{C}$ individually.

In order to overcome this problem we employ parameterised uniformisation introduced in~\cite{CAV2013}. It is a modification of the standard uniformisation technique that allows us to compute strict approximations of the minimal and
maximal transient probability with respect to the set $\mathbf{C}$, moreover, the modification preserves the asymptotic time complexity of standard uniformisation.  For the given state $s \in \mathbb{S}$ and time $t \in \mathbb{R}_{\geq 0}$
the parameterised uniformisation returns vectors $\pi^{\mathbf{C}, s, t}_{\top}$ and $\pi^{\mathbf{C}, s,
  t}_{\bot}$ such that for each state $s' \in \mathbb{S}$ the following holds:
$$
\pi^{\mathbf{C}, s, t}_{\top}(s') \geq max\{\pi^{\mathcal{C}_p, s, t}(s') \mid \mathcal{C}_p \in \mathbf{C}\}  \ \wedge \  \pi^{\mathbf{C}, s, t}_{\bot}(s') \leq min\{\pi^{\mathcal{C}_p, s, t}(s') \mid \mathcal{C}_p \in \mathbf{C}\}
$$

The modification is based on the computation of the local maximum (minimum) of $\pi^{\mathcal{C}_p, s, t}(s')$ 
over all $\mathcal{C}_p \in \mathbf{C}$ for each state \textit{s'} and in each iteration $i$ of standard uniformisation.  
It means that in the \textit{i}th iteration of the computation for a state \textit{s'} we consider only the maximal
(minimal) values in the relevant states in the iteration \textit{i-1}, i.e., the states
that affect $\pi^{\mathcal{C}_p, s, t}(s')$. 

In~\cite{CAV2013} we have defined the function $\sigma(s)$ (formally $\sigma(p,s,\pi)$) which for each state $s \in \mathbb{S}$, perturbation point $p\in \mathbf{P}$ and probability distribution $\pi$ (or pseudo-distribution with the sum smaller or larger than $1$) returns the difference of probability mass inflow and outflow to/from state \textit{s}. If all reactions are described by mass action kinetics the resulting $\sigma$ functions are monotonic with respect to any single perturbed stochastic rate constant $k_r$. This allows us to efficiently compute for each state $s'$ the local maximum (minimum) of $\pi^{\mathcal{C}_p, s, t}(s')$ over all $\mathcal{C}_p \in \mathbf{C}$ corresponding to $\mathbf{P}$.

However, in the case of more complex rate functions than those resulting from mass action kinetics, the corresponding $\sigma(s)$ function does not have to be in general monotonic over $k_r \in [k_r^{\bot},k_r^{\top}]$ for all states \textit{s}. This makes the computation of local extremes with respect to $k_r$ more complex however still tractable. In the following let us assume the space of perturbations $\mathbf{P} = [k_r^{\bot},k_r^{\top}] \times \mathbf{P}'$ will be decomposed along the $k_r$ axis. 

The key idea is for each state \textit{s} to be able to efficiently decompose $\mathbf{P}$ into subspaces $\mathbf{P} = \mathbf{P}_1 \cup \ldots \cup \mathbf{P}_n$, such that for each $\mathbf{P}_i$ the function $\sigma(s)$ over $\mathbf{P}_i$ is monotonic and then use the original method. The problem is a computation of such a strict decomposition into monotonic subspaces is computationally demanding. Therefore we use a simplification, by off-line functional analysis we identify properties of $\sigma$ functions for a given class of reaction kinetics and then obtain a partial decomposition of $\mathbf{P}$ based on function derivations into subspaces where monotonicity is guaranteed. For the remaining subspaces $\mathbf{P}_j$ where monotonicity of $\sigma$ is not guaranteed we employ a less accurate approximation.

We decompose the function $\sigma(s)$ over ${\mathbf{P}_j}$ into functions $\alpha^{s,\mathbf{P}_j}_k$ and  $\beta^{s,\mathbf{P}_j}_l$
such that:
$$
\sigma(s) = \sum_{k=1}^K \alpha^{s,\mathbf{P}_j}_k - \sum_{l=1}^L \beta^{s,\mathbf{P}_j}_l
$$
and each $\alpha^{s,\mathbf{P}_j}_k$ and  $\beta^{s,\mathbf{P}_j}_l$ is monotonic. This allows us to use the original method to compute the maximum and minimum of the functions $\alpha^{s,\mathbf{P}_j}_k$ and $\beta^{s,\mathcal{P}_j}_l$ over the interval $\mathbf{P}_j$, denoted as  $max(\alpha^{s,\mathbf{P}_j}_k)$, $min(\alpha^{s,\mathbf{P}_j}_k)$, $max(\beta^{s,\mathbf{P}_j}_k)$ and $min(\beta^{s,\mathbf{P}_j}_k)$, respectively. Note that, this decomposition can be easily obtained from the definition of the rate function $f_r$. Now the maximum and minimum of $\sigma(p,s)$ over $\mathbf{P_j}$ can be approximated in the following way:
$$
max\{\sigma(p,s) \mid p \in \mathbf{P}_j\} \leq \sum_{k=1}^K max(\alpha^{s,\mathbf{P}_j}_k) - \sum_{l=1}^L min(\beta^{s,\mathbf{P}_j}_l)
$$
$$
min\{\sigma(p,s) \mid p \in \mathbf{P}_i\} \geq \sum_{k=1}^K min(\alpha^{s,\mathbf{P}_i}_k) - \sum_{l=1}^L max(\beta^{s,\mathbf{P}_j}_l).
$$
This approximation increases the inaccuracy of parameterised uniformisation, however, the subspaces~$\mathbf{P}_j$ where  
the monotonicity of $\sigma(s)$ is not guaranteed are usually small and together with perturbation space decomposition introduced in the following section keep on getting smaller. Hence, the additional inaccuracy of the presented extension is manageable. Despite the fact that the time demands of this approximation are orders of magnitudes lower than other numerical methods computing maximum/minimum of $\sigma(s)$ over $\mathbf{P}_i$, they still significantly slow down the computation of parameterised uniformisation.

The aforementioned parameterised uniformisation can be straightforwardly employed also for backward transient analysis. It means that we can efficiently compute the vectors $\tau^{\mathbf{C},\mathbb{A},t}_{\top}$ and $\tau^{\mathbf{C},\mathbb{A},t}_{\bot}$ such that for the given set of states $\mathbb{A}$ and each state $s \in \mathbb{S}$ the following holds: 
$$
\tau^{\mathbf{C}, \mathbb{A}, t}_{\top}(s) \geq max\{\tau^{\mathcal{C}_p, \mathbb{A}, t}(s) \mid \mathcal{C}_p \in \mathbf{C}\}  \ \wedge \  \tau^{\mathbf{C}, \mathbb{A}, t}_{\bot}(s) \leq min\{\tau^{\mathcal{C}_p, \mathbb{A}, t}(s) \mid \mathcal{C}_p \in \mathbf{C}\}
$$

Once we know how to compute the vectors $\tau^{\mathbf{C}, \mathbb{A}, t}_{\top}$ and $\tau^{\mathbf{C}, \mathbb{A}, t}_{\bot}$ the global model checking technique for non-parameterised CTMCs can be directly employed. To obtain the vectors  $\overline{Prob}^{\mathbf{C}, \phi}_{\top},\ \overline{Prob}^{\mathbf{C}, \phi}_{\bot},\  \overline{Exp}^{\mathbf{C}, \mathsf{X}}_{\top}$ and $\overline{Exp}^{\mathbf{C}, \mathsf{X}}_{\bot}$ satisfying Equation~\ref{eq:glob}, it is sufficient to replace the backward transient distribution $\tau^{\mathcal{C}, \mathbb{A},t}_{\bot}$ by the vectors $\tau^{\mathbf{C}, \mathbb{A}, t}_{\top}$ and $\tau^{\mathbf{C}, \mathbb{A}, t}_{\top}$. However for a general class of user-defined post-processing functions \textit{Post}, the vectors $\pi^{\mathbf{C}, s, t}_{\top}$ and $\pi^{\mathbf{C}, s,t}_{\bot}$ cannot be directly used to compute values of $Post^{\mathbf{C}}_{\top}(s,t) = Post(\pi^{\mathbf{C}, s, t}_{\top})$ nor $Post^{\mathbf{C}}_{\bot}(s,t) = Post(\pi^{\mathbf{C}, s, t}_{\bot})$ that would satisfy Equation~\ref{eq:loc} since there is no guarantee about the projective properties of the function \textit{Post}.

Now we show the main idea how to compute $Post^{\mathbf{C}}_{\top}(s,t)$ and $Post^{\mathbf{C}}_{\bot}(s,t)$ for the post-processing function \textit{Post} defined as the mean quadratic deviation of a probability distribution. This function allows us to quantify and analyse a noise in different variants of signalling pathways that are studied in the second case study.   The post-processing function is defined as 
$Post(\pi) = \sum_{s \in \mathbb{S}}{\lvert s(A) - mean(\pi,A)\rvert^2 \cdot \pi(s)}$, where $s(A)$ gives the population of \textit{A} in state \textit{s} and $mean(\pi,A)$ is the mean of the distribution $\pi$ defined as $mean(\pi,A) = \sum_{s \in \mathbb{S}}{s(A) \cdot \pi(s)}$.

Let us suppose we have an upper and lower bound on the probability distribution $\pi_{\top}, \pi_{\bot}$ obtained by the parameterised uniformisation. It means that $\forall \mathcal{C}_p \in \mathbf{C}$ and  $\forall s \in \mathbb{S}. \ \pi_{\bot}(s) \leq \pi^{\mathcal{C}_p}(s) \leq \pi_{\top}(s)$. To find the maximal value $max \left\{ Post(\pi^{\mathcal{C}_p}) \mid \mathcal{C}_p \in \mathbb{C} \right\}$ means to find the distribution $\pi^{max}$ such that  $\sum_{s \in \mathbb{S}}{\pi^{max}(s)} = 1$, $\forall s\in \mathbb{S}. \ \pi_{\bot}(s) \leq \pi^{max}(s) \leq \pi_{\top}(s)$ and the probability mass in $\pi^{max}$ is distributed with the farest distance from the mean. Clearly, such a distribution has a maximal mean quadratic deviation. Note that the number of distributions satisfying the first two conditions is uncountable. Thus we cannot employ direct searching strategy.  

Our searching strategy builds on the observation that only distributions that localise most of the mass as far as possible from the mean (i.e., maximizing the mean quadratic deviation and still meeting the bounds $\pi_{\bot}, \pi_{\top}$), have to be considered. These distributions can be linearly ordered with respect to the sum of mass \textit{x} localised at the low populated part of the state space. It can be shown that the function that evaluates $Post$ on all these distributions is piece-wise quadratic with respect to \textit{x} and has $O(|\mathbb{S}|)$ segments. Therefore, $O(|\mathbb{S}|)$ many steps are sufficient to compute $max \left\{ Post(\pi^{\mathcal{C}_p}) \mid \mathcal{C}_p \in \mathbb{C} \right\}$.

To compute the minimal value $min \left\{ Post(\pi^{\mathcal{C}_p}) \mid \mathcal{C}_p \in \mathbb{C} \right\}$ we proceed analogously, i.e., only the distributions that localise most of the mass as close as possible to the mean are considered.  This leads again to a piece-wise quadratic function. It is also important to note that the perturbation space decomposition presented in the next section allows us to obtain the values $Post^{\mathbf{C}}_{\top}(s,t)$ and $Post^{\mathbf{C}}_{\bot}(s,t)$ with the desired precision. 

\begin{figure} 
\includegraphics[width=\textwidth]{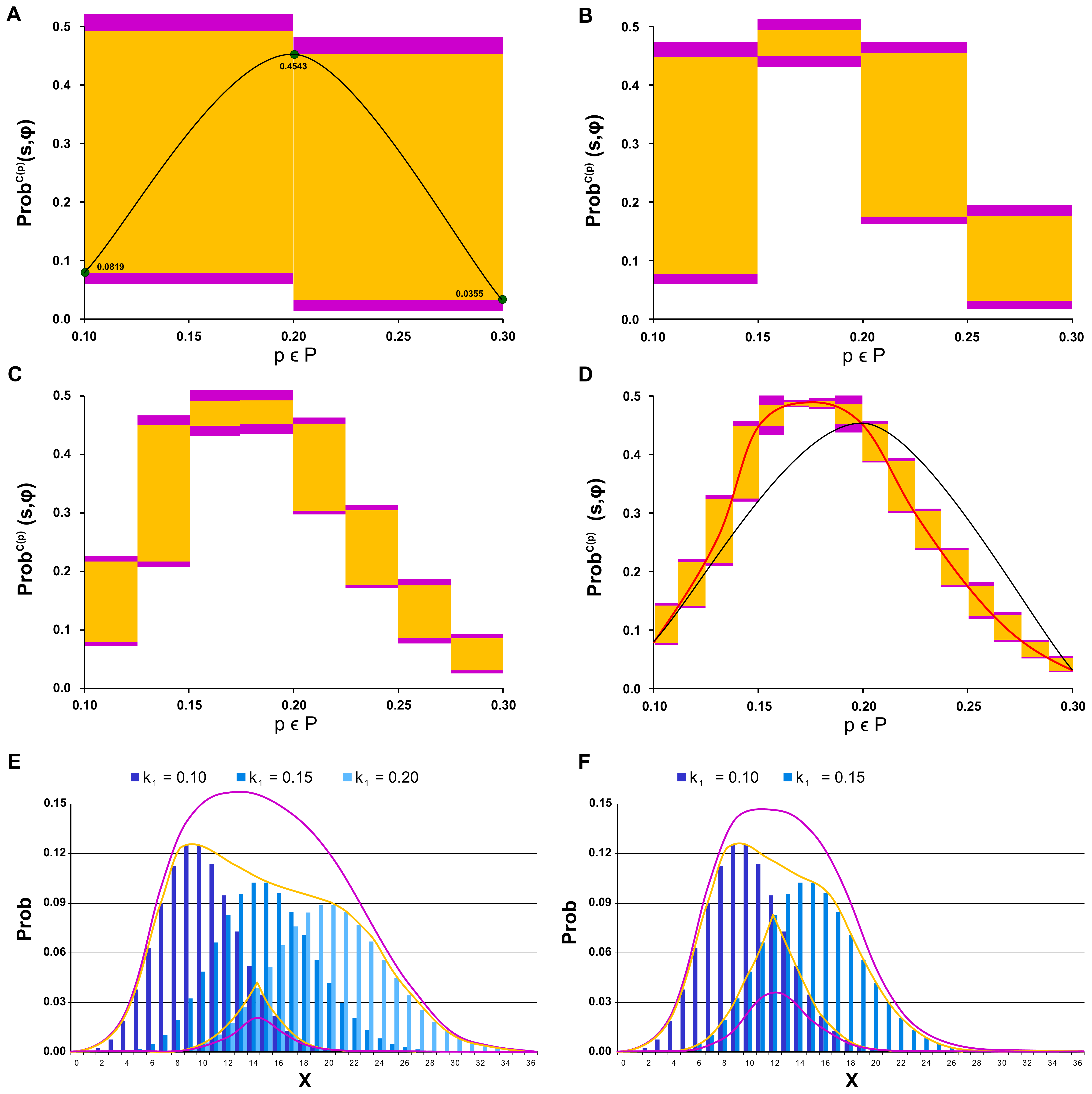}
\caption{{\bf Perturbation space refinement.}
Part (A) depicts three resulting probabilities (green dots) of the formula $\Phi$ for three values of the rate $k_1$ corresponding to three perturbation points $p \in \mathbf{P}$ from Figure~\ref{fig:decomposition} for the initial state $s_0$ denoted as $Prob^{\mathcal{C}_p}(s_0,\Phi)$. The shape of $Prob^{\mathcal{C}_p}(s_0,\Phi)$ for all $p\in \mathbf{P}$ is estimated upon these three points by polynomial interpolation and shown as a black curve. The top four parts (A), (B), (C) and (D) illustrate the min-max approximation of $Prob^{\mathcal{C}_p}(s_0,\Phi)$ for all $p\in \mathbf{P}$ using the decomposition of $\mathbf{P}$ into 2, 4, 8 and 16 subspaces. The exact shape of the probability function for $\Phi$ is visualised as the red thick curve in the (D) and is compared to the initial estimate.
Two types of errors are illustrated: the approximation error is depicted as yellow rectangles and the uniformisation error as the pink rectangles. As can be seen a more refined decompositions reduces both types of errors in each further refined subspace. The bottom parts (E) and (F) depict how the errors arise and how they can be reduce using perturbation space decomposition.}
\label{fig:refinement}
\end{figure}
 
\subsubsection{Perturbation space decomposition}

As we already mentioned, a finite decomposition $\mathbf{P} = \mathbf{P}_1 \cup \ldots \cup \mathbf{P}_n$ into perturbation subspaces is used in order to obtain more accurate approximation of the evaluation function $D^{\mathcal{C}}_{\Phi}$ over the perturbation space $\mathbf{P}$. Before we describe perturbation space decomposition we briefly discuss the key characteristics of parameterised uniformisation that helps us to understand the source of the inaccuracy.   The most important fact is that parameterised uniformisation for the
set $\mathbf{C}$ in general does not correspond to standard uniformisation for
any CTMC $\mathcal{C}_p \in \mathbf{C}$. The reason is that we consider a
behaviour of a parameterised CTMC that has no equivalent counterpart in any particular $\mathcal{C}_p$. First, the
parameter $k_r$ (minimizing/maximizing the inspected value) is determined locally for each
state. Therefore, in a single iteration there can exist two different states such that in one state the parameterised uniformisation selects $k_r=k_r^{\top}$ while in another state it selects $k_r=k_r^{\bot}$. Second, the
parameter is determined individually for each iteration and thus for a state
$s_i$ the parameter $k_r$ can be chosen differently in individual iterations.

Inaccuracy of the proposed  min-max approximation related to the computation of 
 parameter\-ised uniformisation, called \emph{unification error}, is given
as:
$$ (D^{\mathcal{C}}_{\Phi,\mathbf{P},\top} - max \{ D^{\mathcal{C}}_{\Phi}(p) \mid p \in \mathbf{P} \}) +
(min \{ D^{\mathcal{C}}_{\Phi}(p) \mid p \in \mathbf{P} \} - D^{\mathcal{C}}_{\Phi,\mathbf{P},\bot}).$$
Apart from the unification error our approach introduces an inaccuracy related to
approximation of the evaluation function, called
\emph{approximation error}, given as:
$$  max \{ D^{\mathcal{C}}_{\Phi}(p) \mid p \in \mathbf{P} \} - min \{ D^{\mathcal{C}}_{\Phi}(p) \mid p \in \mathbf{P} \}.$$
Finally, the \emph{overall error} of the min-max approximation, denoted as
$\mathsf{Err}^{\mathcal{C}}_{\Phi, \mathbf{P}}$, is defined as a sum of both errors,
i.e., $D^{\mathcal{C}}_{\Phi,\mathbf{P},\top} - D^{\mathcal{C}}_{\Phi,\mathbf{P},\bot}$.  
Figure~\ref{fig:refinement} illustrates
 both types of errors. The approximation error is depicted as yellow rectangles
and the unification error is depicted as the purple rectangles.

We are not able to effectively distinguish the proportion of the approximation
error and the unification error nor to reduce the unification error as
such. Therefore, we design a method based on the perturbation space decomposition
that allows us to effectively reduce the overall error of the min-max
approximation to a user specified \emph{absolute error bound}, denoted as
$\textsc{Err}$.    

In order to ensure that the min-max approximation meets the given absolute error
bound $\textsc{Err}$, we iteratively decompose the perturbation space $\mathbf{P}$
into finitely many subspaces such that $\mathbf{P} = \mathbf{P}_1 \cup \ldots \cup
\mathbf{P}_n$ and each partial result satisfies the overall error bound, i.e., $\forall \ j: 1 \leq j \leq n: \mathsf{Err}^{\mathcal{C}}_{\Phi, \mathbf{P}_j} \leq \textsc{Err}$. Therefore, the overall error equals to 
$$\mathsf{Err}^{\mathcal{C}}_{\Phi, \mathbf{P}}= \sum_{j =
  1}^{n}{\frac{|\mathbf{P}_j|}{|\mathbf{P}|} \left(D^{\mathcal{C}}_{\Phi,\mathbf{P}_j,\top} - D^{\mathcal{C}}_{\Phi,\mathbf{P}_j,\bot} \right)} \leq \sum_{j =
  1}^{n}{\frac{|\mathbf{P}_j|}{|\mathbf{P}|} \textsc{Err}} = \textsc{Err}.$$
Figure~\ref{fig:refinement} illustrates such a decomposition and demonstrates
convergence of $\mathsf{Err}^{\mathcal{C}}_{\Phi, \mathbf{P}_j}$ to 0 provided that the evaluation function $D^{\mathcal{C}}_{\Phi}$
over $\mathbf{P}$ is continuous.

For sake of simplicity, we present the parametric decomposition only on the computation of
$\pi^{\mathbf{C}, s, t}_{\star}$ and $\tau^{\mathbf{C}, \mathbb{A}, t}_{\star}$ for  $\star \in \{\top,\bot\}$ and $\mathbf{P}$ since it can be easily extended to the computation of $D_{\Phi,\mathbf{P},\star}^{\mathcal{C}}$ for any formula~$\Phi$. The key part of the parametric decomposition is to decide when the inspected subspace should be further decomposed. The condition for the decomposition is different for $\pi^{\mathbf{C}, s, t}_{\star}$ and $\tau^{\mathbf{C}, \mathbb{A}, t}_{\star}$. Since the vector $\pi^{\mathbf{C}, s, t}_{\star}$ gives us the transient probability distribution from the state $s$ that is further used to compute $D^{\mathcal{C},s}_{\Phi,\mathbf{P},\star}$, we consider the following condition. The space $\mathbf{P}$ (represented by the CTMC $\mathbf{C}$) is decomposed if during the computation of parameterised uniformisation in an iteration $i$ it holds that:
$$  \sum_{k =1}^{|\mathbb{S}|} {\pi^{\mathbf{C}, s, i}_{\top}(s_k)} - \sum_{k =1}^{|\mathbb{S}|} {\pi^{\mathbf{C}, s, i}_{\bot}(s_k)} > \textsc{Err}$$
where $ \pi^{\mathbf{C}, s, i}_{\star}$  denotes the corresponding approximation of $\pi^{\mathbf{C}, s, 0} \cdot (\mathbf{Q}^{\mathsf{unif}(\mathbf{C})})^i$. 

In contrast to $\pi^{\mathbf{C}, s, t}_{\star}$,  the value $\tau^{\mathbf{C}, \mathbb{A}, t}_{\top}(s)$  for each state $s\in \mathbb{S}$ is further used to  $D^{\mathcal{C},s}_{\Phi,\mathbf{P},\star}$ and thus we consider the different condition. The space $\mathbf{P}$  is decomposed if during the computation of parameterised uniformisation in an iteration $i$ for any state $s$ it holds that:
$$   \tau^{\mathbf{C}, \mathbb{A}, i}_{\top}(s) - \tau^{\mathbf{C}, \mathbb{A}, i}_{\bot}(s) > \textsc{Err}.$$

If the decomposition takes place we cancel the current computation and decompose the perturbation space  $\mathbf{P}$ to $n$ subspaces such that $\mathbf{P} = \mathbf{P}_1 \cup
\ldots \cup \mathbf{P}_n$. Each subspace $\mathbf{P}_j$ defines a new set
of CTMCs $\mathbf{C}_j = \{\mathcal{C}_j \mid j \in \mathbf{P}_j \}$ that is
independently processed in a new computation branch. Note that we could reuse
the previous computation and continue from the iteration~$i-1$. However, the
most significant part of the error is usually cumulated during the previous iterations and thus the decomposition would have only a negligible impact on error reduction.

\emph{A minimal decomposition with respect to the perturbation space  $\mathbf{P}$} defines a
minimal number of subspaces \textit{m} such that $\mathbf{P} = \mathbf{P}_1 \cup \ldots
\cup \mathbf{P}_{m}$ and for each subspace $\mathbf{P}_j$ where $1 \leq j \leq m$
holds that $\mathsf{Err}^{\mathcal{C},s}_{\Phi, \mathbf{P}_j} \leq \textsc{Err}$ where $\mathsf{Err}^{\mathcal{C},s}_{\Phi, \mathbf{P}_j} = D^{\mathcal{C}, s}_{\Phi,\mathbf{P}_j,\top} - D^{\mathcal{C},s}_{\Phi,\mathbf{P}_j,\bot}$. Note that
the existence of such decomposition is guaranteed only if the evaluation function
$D^{\mathcal{C},s}_{\Phi}$
over $\mathbf{P}$ is continuous. If the evaluation function is continuous there can exist more than one minimal decomposition.
However, it can not be straightforwardly found. To overcome this problem
we have considered and implemented several heuristics allowing to iteratively
compute a decomposition satisfying the following: (1) it ensures the required
error bound whenever $D^{\mathcal{C},s}_{\Phi}$
over $\mathbf{P}$ is continuous, (2) it guarantees the refinement termination in the situation where$D^{\mathcal{C},s}_{\Phi}$
over $\mathbf{P}$ is not continuous and the discontinuity causes that $\textsc{Err}$ can not be achieved. To ensure the termination an additional parameter has to be introduced as a lower bound on the subspace size. Hence this parameter provides a supplementary termination criterion.

\subsubsection{Implementation}

We delivered a prototype implementation of the framework for the robustness analysis on top of the tool PRISM 4.0~\cite{KNP11}. This tool provides the appropriate modelling and specification language. Our implementation builds on sparse engine that uses data structures based on the sparse matrices. They provides suitable representation of  models for the time efficient numerical computation.  

In the case that large number of perturbation subspaces is required to obtain the desired accuracy of the approximation the sequential computation can be extremely time consuming. However, our framework allows very efficient parallelization since the the computation of particular subspaces is independent and thus can be executed in parallel. Our implementation enables the parallel computation and thus the robustness analysis can be significantly accelerated using high performance parallel hardware architectures.

\section{Results}

\subsection{Gene Regulation of Mammalian Cell Cycle}

\begin{figure} 
\begin{center}
\includegraphics[width=0.5\textwidth]{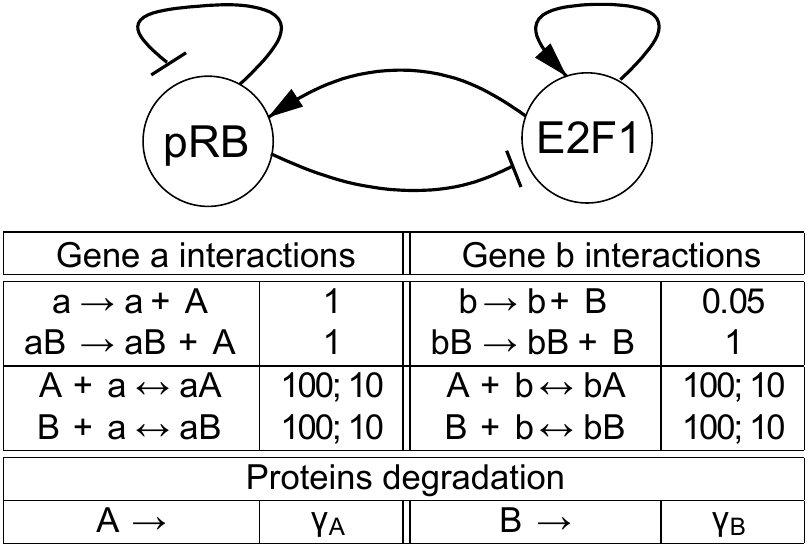}
\end{center}
\caption{{\bf Model of regulation of the mammalian cell cycle.}
The core gene regulatory module controlling the $G_1/S$-phase transition in the cell cycle of mammalian cells~\cite{Keletal00} is depicted in the upper part. The retinoblastoma protein \textit{pRB} (A) [HumanCyc:HS06650] interacts with the retinoblastoma-binding transcription factor $E_2F_1$ (B) [HumanCyc:HS02261]. In high concentration levels, the $E_2F_1$ protein activates the $G_1/S$ transition mechanism. On the other hand, a low concentration of $E_2F_1$ prevents committing to $S$-phase. Positive autoregulation of $E_2F_1$ causes bi-stability.\\[2mm]
Stochastic mass action reformulation of the $G_1/S$ regulatory circuit is shown in the table below. The gene regulation is modelled by means of a set of second-order reactions simplifying the elementary processes behind transcription. In particular, the model includes the interactions among transcription factors (\textit{A}, \textit{B} stand for \textit{pRB} and $E_2F_1$, respectively) and respective genes and protein production/degradation reactions. The interactions are represented by reversible TF-gene binding reactions in the second row of the table (genes are denoted by small letters). Individual protein production reactions controlled by these interactions are represented by the irreversible gene expression reactions in the first row of the table. Protein degradation is modelled as spontaneous by means of first-order reactions. Kinetic coefficients are set only approximately provided that they are considered equal for all instances of a particular process (binding, dissociation, promoted protein production). The only exception is the spontaneous (basal) expression of \textit{b} which is set to a low rate. This mimics the fact that $E_2F_1$ is only rapidly produced under the circumstances of self-activation~\cite{Swatetal04}. Degradation parameters are left unspecified.}
\label{fig:modelg1s}
\end{figure}

We have applied the robustness analysis  to the gene regulation model published
in~\cite{Keletal00}, the regulatory network is shown in
Fig.~\ref{fig:modelg1s} (left). The model explains regulation of a transition between
early phases of the mammalian cell cycle. In particular, it targets the
transition from the control $G_1$-phase to \textit{S}-phase (the synthesis
phase). $G_1$-phase makes an important checkpoint controlled by a \emph{bistable
  regulatory circuit} based on an interplay of the retinoblastoma protein \textit{pRB},
denoted by \textit{A} (the so-called tumour suppressor, HumanCyc:HS06650) and the
retinoblastoma-binding transcription factor $E_2F_1$, denoted by \textit{B} (a central
regulator of a large set of human genes, HumanCyc:HS02261). In high
concentration levels, the $E_2F_1$ protein activates the $G_1$/$S$ transition
mechanism. On the other hand, a low concentration of $E_2F_1$ prevents committing
to \textit{S}-phase.

Positive autoregulation of \textit{B} causes bi-stability of its concentration depending
on the parameters. Especially, of specific interest is the degradation rate of
\textit{A}, $\gamma_{A}$. In~\cite{Swatetal04} it is shown that for increasing
$\gamma_{A}$ the low stable mode of \textit{B} switches to the high stable
mode.  
When mitogenic stimulation increases under conditions of active growth, rapid phosphorylation of
\textit{A} starts and makes the degradation of unphosphorylated \textit{A} stronger (the
degradation rate $\gamma_{A}$ increases). This causes \textit{B} to lock in the high
stable mode implying the cell cycle commits to \textit{S}-phase. Since mitogenic stimulation influences the degradation rate of \textit{A}, our goal is to study the population distribution around the low and high steady state and to explore the effect of $\gamma_{A}$ by means of the evaluation function.

It is necessary to note that the original ODE model in~\cite{Swatetal04} has been formalised by means of Hill kinetics representing the cooperative action of transcription factor molecules.
Since Hill kinetics cannot be directly transferred to stochastic modelling~\cite{Garaietal12,Sanftetal11}, we have reformulated the model in the framework of stochastic mass action
kinetics~\cite{Gillespie1977}. The resulting reactions are shown in
Fig.~\ref{fig:modelg1s} (right). Since the detailed knowledge of elementary chemical
reactions occurring in the process of transcription and translation is incomplete, we
use the simplified form as suggested in~\cite{Gillespieetal05}. In the minimalist
setting, the reformulation requires addition of rate parameters describing the
transcription factor--gene promoter interaction while neglecting cooperativeness
of transcription factors activity. Our parameterisation is based on time-scale
orders known for the individual processes~\cite{Yang2003} (parameters considered in $s^{-1}$). Moreover, we assume
the numbers of \textit{A} and \textit{B} are bounded by 10 molecules. 
Correctness of the upper bounds for \textit{A} and \textit{B} was validated by observing thousand independent stochastic simulations. We consider minimal population number distinguishing the two stable modes.
All other species are bounded by the initial number of DNA molecules
(genes \textit{a} and \textit{b}) which is conserved and set to~1. The corresponding CTMC 
has 1078 states and 5919 transitions.

We consider two hypotheses: (1) stabilisation in the low mode where $B<3$, 
(2)~stabilisation in the high mode where $B>7$.  Both hypotheses are expressed within time horizon 1000 seconds reflecting the
time scale of gene regulation response.  According
to~\cite{Swatetal04}, we consider the perturbation space $\gamma_{A}\in
[0.005,0.5]$. For both hypothesis we consider three different settings of $\gamma_{B}$: $\gamma_B=0.05$, $\gamma_B=0.10$, and $\gamma_B=0.15$.

\begin{figure} 
\centering
\includegraphics[width=\textwidth]{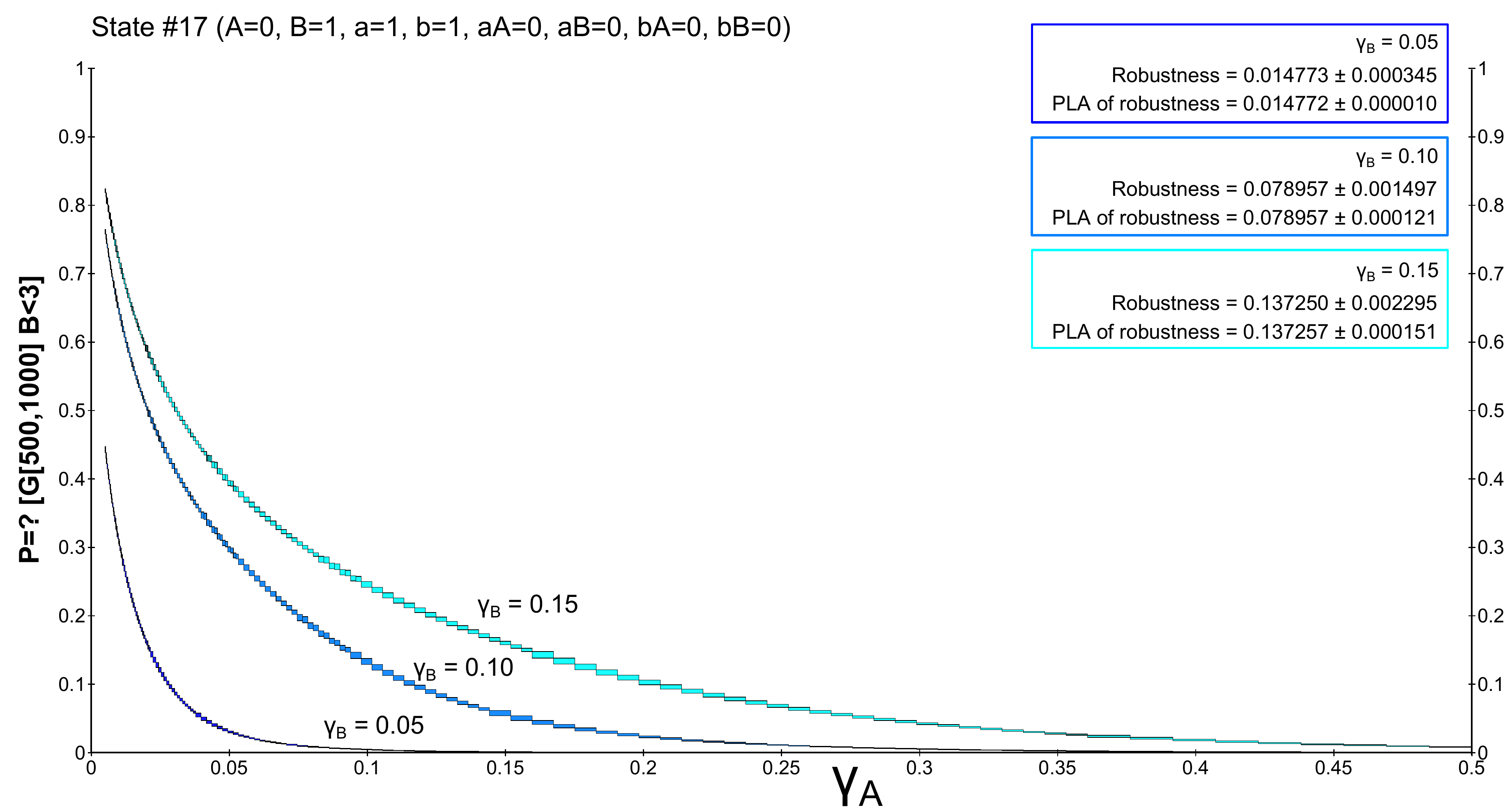}
\caption{{\bf Results of robustness analysis for hypothesis (1) using a until operator.} 
Hypothesis (1) requires stabilisation of $E_2F_1$ in the low concentration mode ($B<3$). A CSL formula with the until operator is used in this case. Each of the curves represents the evaluation function over $\gamma_A$ degradation obtained for a particular setting of $\gamma_B$. More precisely, the horizontal axis shows the perturbation of \textit{pRB} degradation rate and the vertical axis shows the probability of the hypothesis to be satisfied. In the upper left corner, robustness values are shown for each of the curves. The values are displayed with the absolute error quantifying the precision of the approximate method. For comparison, the values are computed also on piece-wise affine approximations of the evaluation function. It can be seen that the robustness values are small which is due to the fact that fluctuations of molecular numbers cause frequent exceeding of the required bound in the considered time horizon.}
\label{fig:cs1_until_0-2}
\end{figure}

 We employ two
alternative CSL formulations to express the hypothesis (1).
First, we express the property of being inside the given bound during the time interval $I=[500,1000]$ using globally operator: $\mathsf{P}_{=?}[\mathsf{G}^I\, (B< 3)]$. The interval starts from 500 seconds in order to bridge the initial fluctuation region and let the system stabilise. The resulting landscape visualisation is depicted in Figure~\ref{fig:cs1_until_0-2} together with the robustness values computed for individual cases.
Since the stochastic noise causes molecules to repeatedly escape the requested bound, the resulting probability is
significantly lower than expected. Namely, in the case $\gamma_B=0.05$ the resulting probability is close to 0 for almost all considered parameter values implying very small robustness. Increasing of the \textit{B} degradation rate causes an observable increase in robustness.

\begin{figure} 
\centering
\includegraphics[width=\textwidth]{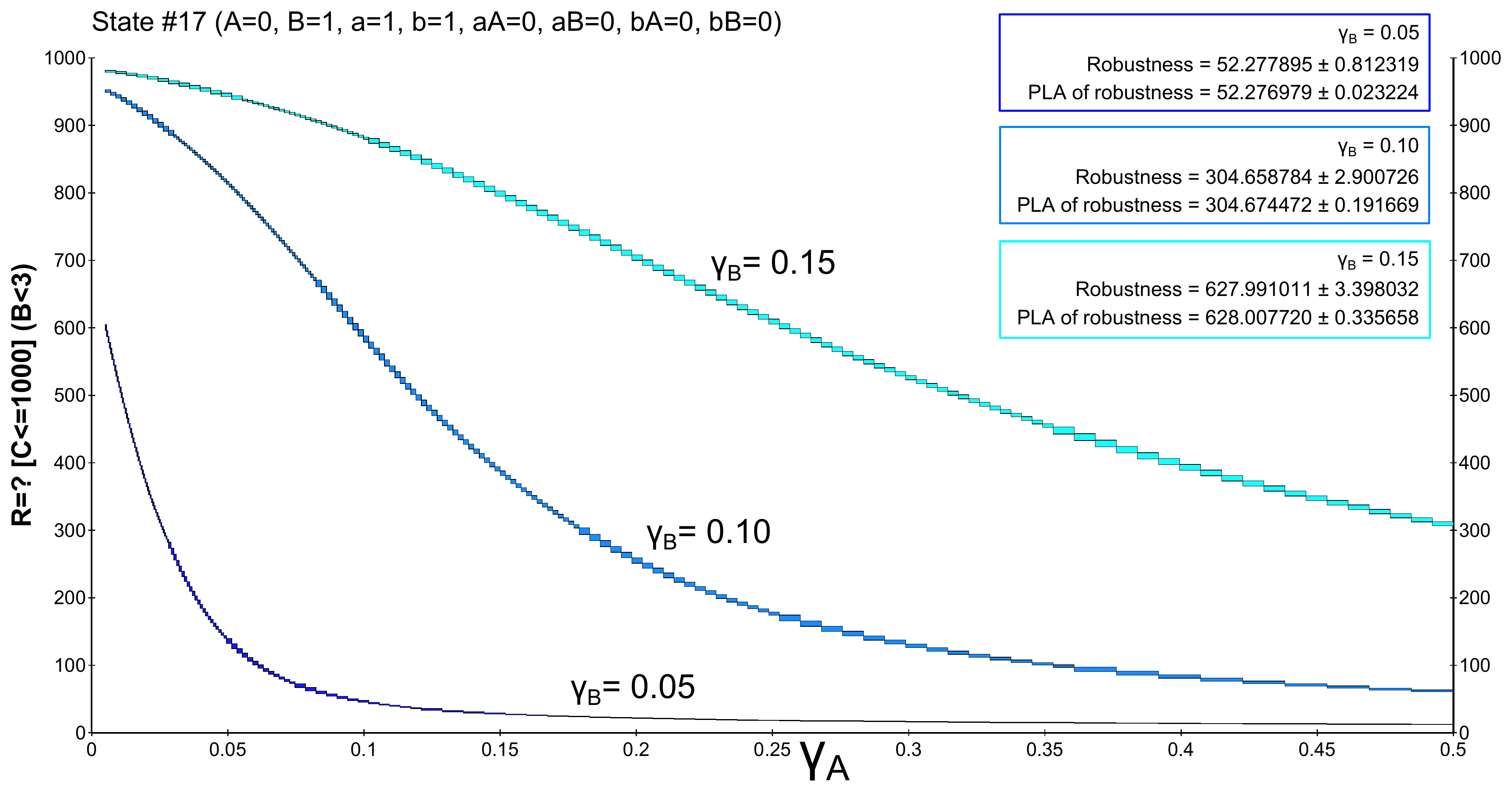}
\caption{{\bf Results of robustness analysis for hypothesis (1) using a reward operator.} 
Hypothesis (1) requires stabilisation of $E_2F_1$ in the low concentration mode ($B<3$). A CSL formula with cumulative reward operator is used in this case. Each of the curves represents the evaluation function over $\gamma_A$ degradation obtained for a particular setting of $\gamma_B$. More precisely, the horizontal axis shows the perturbation of \textit{pRB} degradation rate and the vertical axis shows the probability of the hypothesis to be satisfied. In the upper left corner, robustness values are shown for each of the curves. The values are displayed with the absolute error quantifying the precision of the approximate method. For comparison, the values are computed also on piece-wise affine approximations of the evaluation function. It can be seen that the robustness values change rapidly with different settings of $\gamma_B$. This observation goes with the fact that with faster degradation of $E_2F_1$ there is a higher probability that the positively self-regulated protein is locked in the stable mode of no production. The decrease of the value with increasing $\gamma_A$ is due to the weakening effect of inhibition by \textit{pRB}.}
\label{fig:cs1_reward_0-2}
\end{figure}

In order to avoid fluctuations of affecting the result, we use a cumulative reward property to capture the fraction of the time the system has the required number of molecules within the time interval $[0,1000]$: 
$\mathsf{R}_{=?}[\mathsf{C}^{\leq t}](B<3)$ where $t = 1000$ and $\mathsf{R}_{=?}[\mathsf{C}^{\leq t}](B \sim X)$ denotes that state reward $\rho$ is defined such that $\forall s\in \mathbb{S}.\rho(s) = 1$ iff $B \sim X$ in $s$. The resulting landscape visualisation is shown in Figure~\ref{fig:cs1_reward_0-2}. Here the effect of increase of robustness value with respect to increasing $\gamma_B$ is significantly stronger.

After normalising the robustness values, we can observe that the model is significantly more robust with respect to the cumulative reward-based formulation of the hypothesis. This goes with the fact that the reward property neglects the frequent fluctuations in the given time horizon. 

When focusing on the phenomenon of bistability, we can conclude that the most significant variance in the molecule population with respect to the two stable modes is observed in the range $\gamma_A=[0.15,0.3]$ with $\gamma_B=0.10$. Here the distribution of the behaviour targeting the low and high mode is diversified nearly uniformly (especially for $\gamma_A=0.2$). Note that in this case there is a significant amount of behaviour (around $40\%$) not converging to either of the two modes.

\begin{figure} 
\centering
\includegraphics[width=\textwidth]{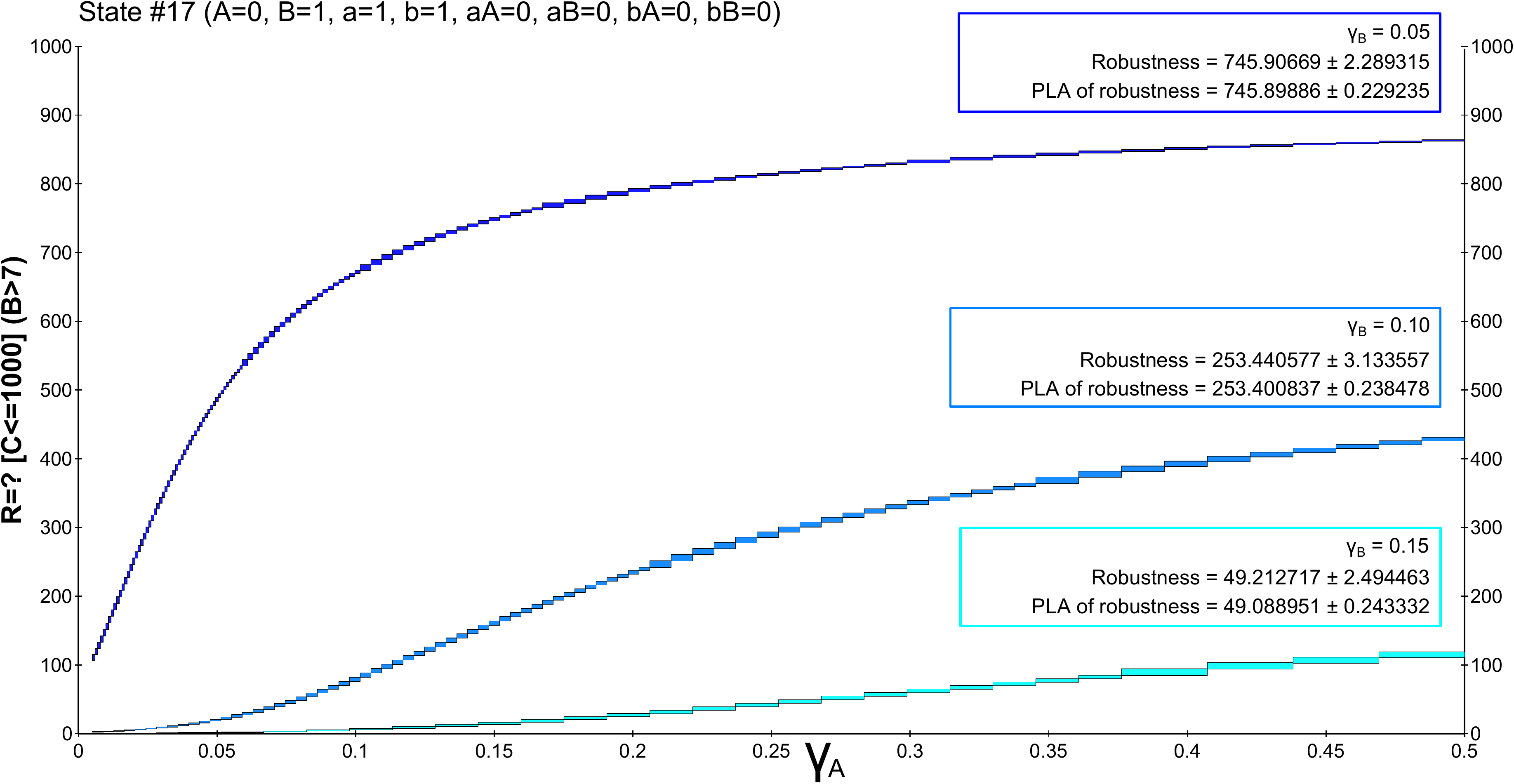}
\caption{{\bf Results of robustness analysis for hypothesis $2$.} 
Hypothesis (2) requires stabilisation of $E_2F_1$ in the high concentration mode ($B>7$). A CSL formula with cumulative reward operator is employed. Each of the curves represents the evaluation function over $\gamma_A$ degradation obtained for a particular setting of $\gamma_B$. The horizontal axis shows the perturbation of \textit{pRB} degradation rate and the vertical axis shows the probability of the hypothesis to be satisfied. In the upper left corner, robustness values are shown for each of the curves. The values are displayed with the absolute error quantifying the precision of the approximate method. For comparison, the values are computed also on piece-wise affine approximations of the evaluation function. It can be seen that the robustness values change rapidly with different settings of $\gamma_B$. This observation goes with the fact that with faster degradation of $E_2F_1$ there is a lower probability that the positively self-regulated protein is locked in the stable mode of no production. In particular, the high stable mode is preferred for lower values of $\gamma_B$. The increase of the value with increasing $\gamma_A$ is due to the weakening effect of inhibition by \textit{pRB}.}
\label{fig:cs1_reward_8-10}
\end{figure}

To encode the hypothesis $2$ we employ the reward-based formulation: $\mathsf{R}_{=?}[\mathsf{C}^{\leq t}](B>7)$.  The time interval is set to be the same as in the previous case ($t=1000$). The resulting landscape visualisations for individual settings of $\gamma_B$ are depicted in Figure~\ref{fig:cs1_reward_8-10}. It can be observed that the effect of $\gamma_B$ is now inverse which goes with the fact that higher rate of $E_2F_1$ degradation causes the rapid dynamics of the protein and decreases the amenability of the cell to commit to \textit{S}-phase (by making the hypothesis $1$ more robust than hypothesis $2$). 

\begin{figure} 
\centering
\includegraphics[width=\textwidth]{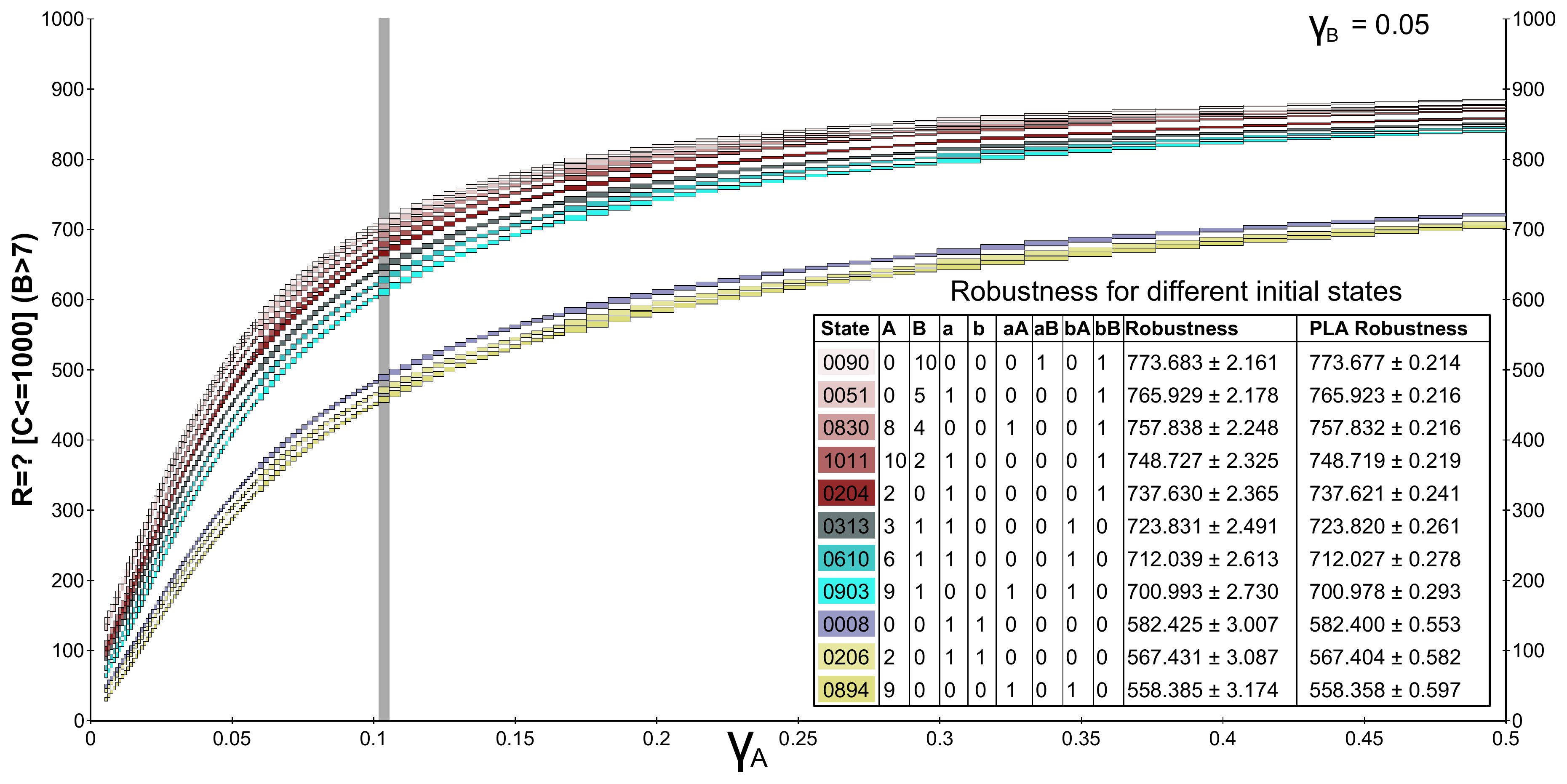}
\caption{{\bf Landscape visualisation for hypothesis (2) and several selected initial states.} 
The landscape visualisation of hypothesis (2) (stabilisation of $E_2F_1$ in the high concentration mode $B>7$) is shown for several selected initial states of the whole state space. A CSL formula with cumulative reward operator is employed. Each of the curves represents the evaluation function over $\gamma_A$ degradation obtained for a particular initial state and $\gamma_B$ set to 0.05. The legend shows the amount of individual species in particular initial states and the robustness of the hypothesis is given together with the absolute error. The results obtained by piece-wise affine approximation are also shown. It can be seen that the hypothesis is only negligibly sensitive to initial conditions. Especially, only states with zero initial concentration of $E_2F_1$ cause $E_2F_1$ to attain low molecular numbers thus lowering the robustness of the hypothesis. The grey vertical line shows the small perturbation in $\gamma_A$ which is further explored in detail in Figure~\ref{fig:cs1_selected_interval_statemap}.}
\label{fig:cs1_init_states_combine}
\end{figure}

An interesting observation coming out of the analysis is that the selection of an initial state has only a negligible impact on the result. This is exploited in Figure~\ref{fig:cs1_init_states_combine} where we have selected 11 states uniformly distributed throughout the state space. Although low initial numbers of \textit{B} slightly decrease robustness of hypothesis (2), the difference is not very big. 

\begin{figure} 
\centering
\includegraphics[width=.8\textwidth]{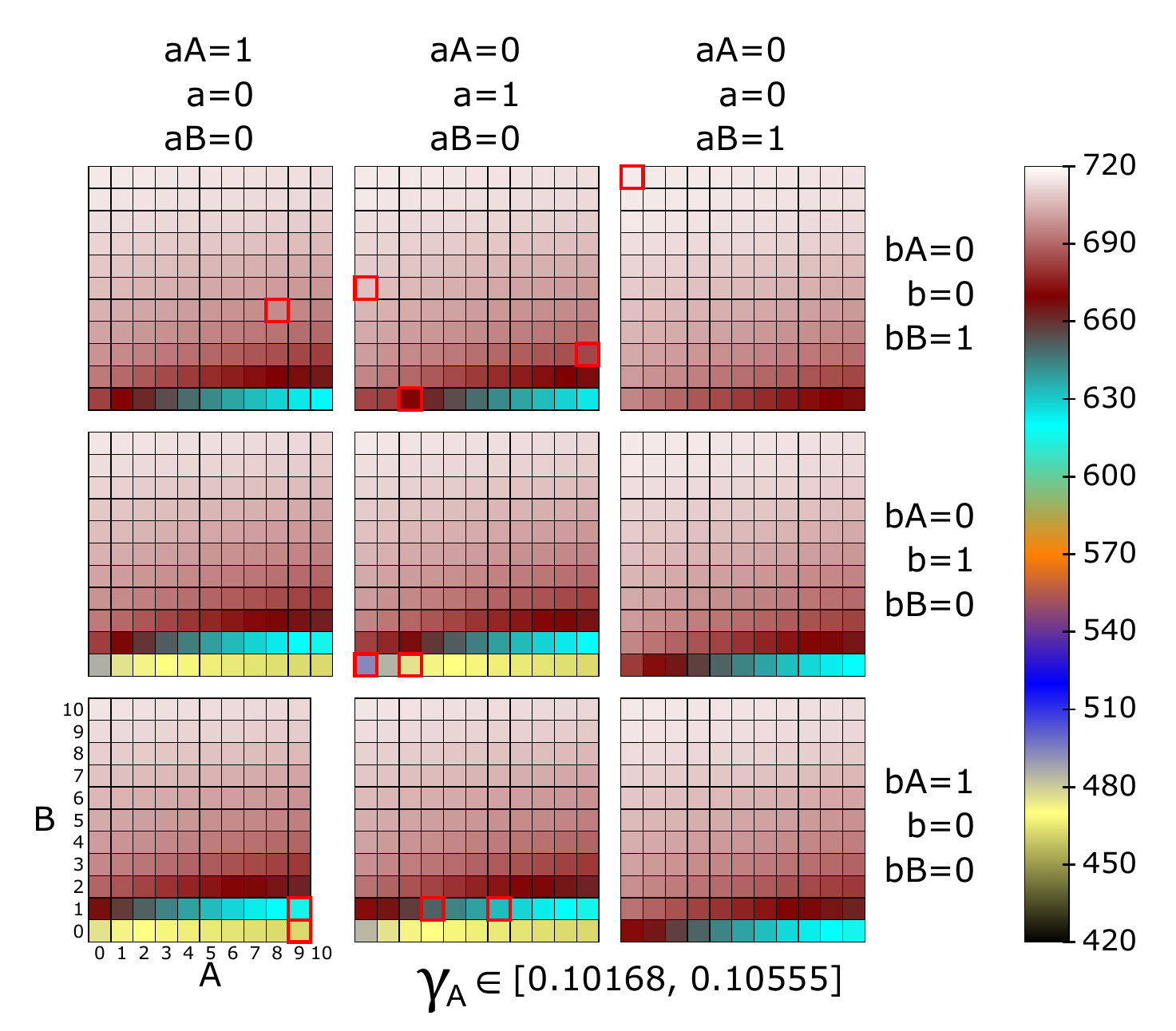}
\caption{{\bf Analysis of hypothesis (2) for all initial states.} 
Hypothesis (2) (stabilisation of $E_2F_1$ in the high concentration mode $B>7$) is computed and visualised for all initial states in the considered perturbation space $(\gamma_A,\gamma_B) \in [0.10168, 0.10555] \times [0.05]$. Because we assume at most a single molecule of DNA in the system, state variables denoting genes and gene-protein complexes have a binary domain. There are only two variables having a larger domain (0-10), in particular, these are the proteins \textit{pRB} and $E_2F_1$. Therefore each of the (binary) combinations is visualised for the entire domain of \textit{A} and \textit{B} in a separate box. The colour intensity of each box in the grid shows the upper bound of the cumulative reward evaluated for the respective initial state. It can be seen that the hypothesis is mostly insensitive to selection of initial states. Only the initial zero level of $E_2F_1$ (\textit{B, bB, aB}) causes a decrease of the resulting value. States selected in Figure~\ref{fig:cs1_init_states_combine} are highlighted in red.}
\label{fig:cs1_selected_interval_statemap}
\end{figure}

More detailed insight can be inferred from Figure~\ref{fig:cs1_selected_interval_statemap} where hypothesis (2) evaluation is exploited for a small perturbation of $\gamma_A$ with respect to the entire initial state space. The considered perturbation is highlighted in Figure~\ref{fig:cs1_init_states_combine} by the grey vertical line. The colour intensity of the grid shows the upper bound of the cumulative reward evaluated for the respective initial state. It can be seen that the hypothesis is really insensitive to selection of initial states. Only the initial zero level of \textit{B} causes a decrease of the resulting value. Moreover, this happens (naturally) just in two kinds of states: (\textit{i}) no molecule of \textit{B} is bound to any of the genes, i.e., self-activation of \textit{b} is inactive and the expression of \textit{b} occurs in the spontaneous mode having a low rate 0.05; (\textit{ii}) a molecule of \textit{A} is bound to \textit{b} thus imposing the inhibition of \textit{b} and causing the same scenario.

\subsection{Robustness of two-component signalling systems response}

Signalling pathways make the main interface between cells and their environment. Their main role is to sense biochemical conditions outside the cell and to transfer this information into the internal logical circuits (gene regulation) of the cell. Since signal processing is realised by several dedicated protein complexes (signalling components), it is naturally amenable to intrinsic noise in these protein populations caused by stochasticity of transcription/translation processes. Robust input-output signal mapping is crucial for cell functionality. Many models and experimental studies have been conducted attempting to explain mechanisms of robust signal processing in procaryotic cells, e.g.,~\cite{Batchelor21012003,Shinar11122007}.

\begin{figure} 
\begin{center}
\centering
\includegraphics[width=0.5\textwidth]{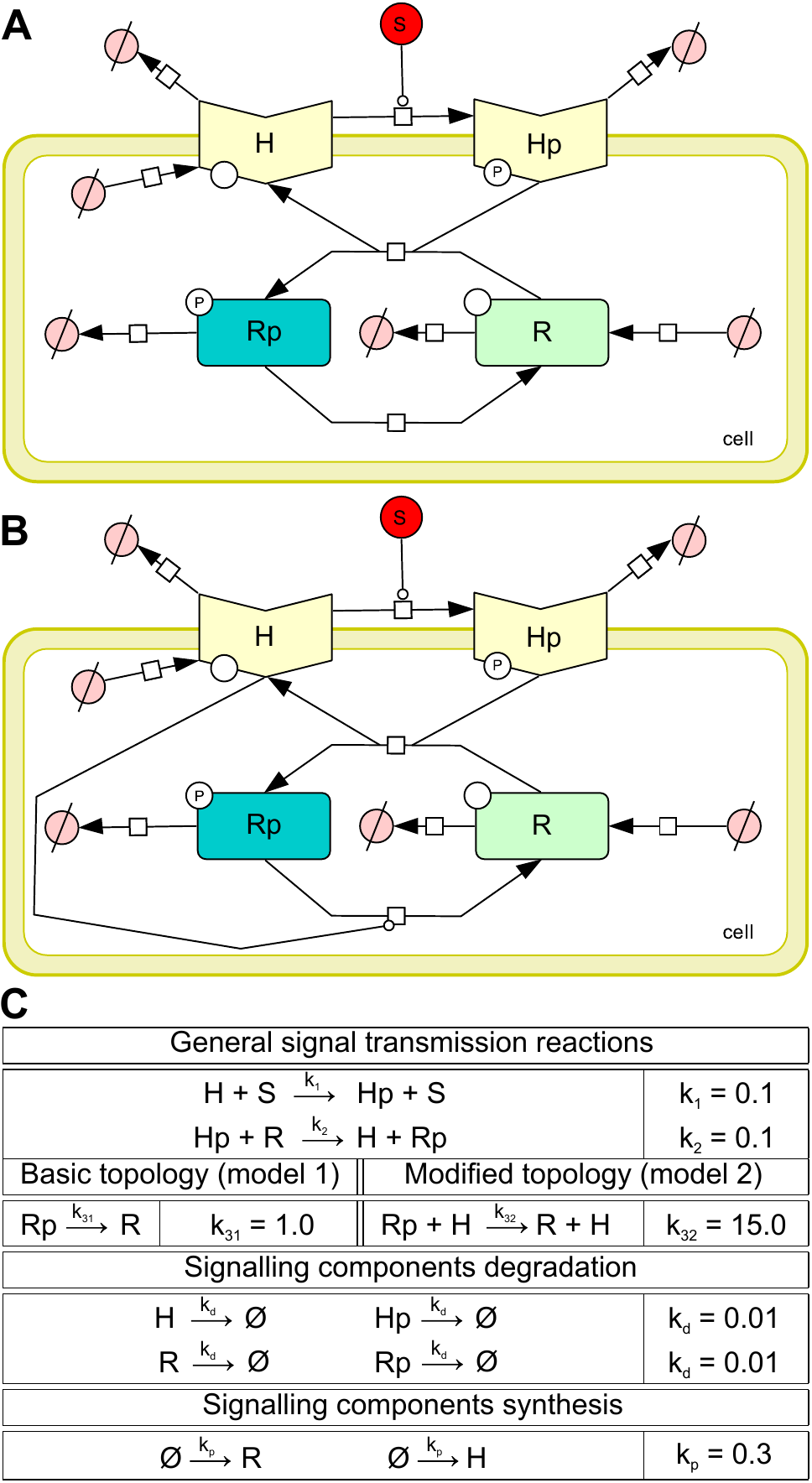}
\end{center}
\caption{{\bf Model of a two-component signalling pathway.}
(A) Basic topology of the two-component signalling pathway. (B) Modified topology of the two-component signalling pathway, additionally, histidine kinase \textit{H} catalyses dephosporylation of the response regulator \textit{R}. (C) Reactions specifying the biochemical model of the two considered topologies of the two-component signalling pathway. Phosphorylation of the first component \textit{H} catalysed by the input signal \textit{S} and phosporylation of the second component \textit{R} are shared by both topologies, the only difference is in the second component dephophorylation. Additionally, we consider unregulated proteosynthesis/degradation reactions for both topology variants. Reaction topology in (A) and (B) was created using CellDesigner~\cite{Celldesigner}.
}
\label{fig:sigtopol}
\end{figure}

In order to construct robust signalling circuits in synthetically modified procaryotic cells, Steuer et al.~\cite{steuer2011robust} has suggested and analysed a modification of a well-studied two-component signalling pathway that is insensitive to signalling component concentration fluctuations. The study has been performed by using a simplified model consisting of the two signalling components each considered in both phosphorylated and unphosporylated forms.
The first component, the histidine kinase \textit{H}, is a membrane-bound receptor phosphorylated by an external signalling ligand \textit{S}. In its phosphorylated form \textit{Hp}, the histidine kinase transfers the phospho-group onto the second component -- the response regulator \textit{R}. That way it activates the response regulator by transforming it into the phosphorylated form \textit{Rp} which is diffusible and functions as the internal signal for the cell. The basic topology of the pathway is depicted in Figure~\ref{fig:sigtopol}A. The modification suggested by Steuer et al. is depicted in Figure~\ref{fig:sigtopol}B. The difference is in the addition of catalytic activation of \textit{Rp} dephosporylation by the unphosphoshorylated histidine kinase \textit{H}. In~\cite{steuer2011robust} it has been rigorously proven that under the deterministic setting this modification leads to globally robust steady-state response of the signalling pathway that is not achievable with the basic topology. 

We reformulate the model in the stochastic setting and employ our method to provide detailed analysis of the input-output signal response under fluctuations in population of both signalling components. In contrast to~\cite{steuer2011robust} where average steady-state population is analysed with respect to fluctuations in signalling components, our analysis refines the steady population in terms of distributions. That way we obtain for a stable input signal a detailed view of distribution of the output response. In particular, instead of studying the effect of perturbations on the average population, we see how perturbations affect the distribution, i.e., the variance (fluctuation) in the output response. That way the stochastic framework gives a more detailed insight into the input-output signal response mechanism.

The biochemical model of both topology variants is given in Figure~\ref{fig:sigtopol}C.
The input signal \textit{S} is considered to be fixed and therefore it makes a constant parameter of the model. The signalling components in both phosporylated and unphosporylated forms make the model variables \textit{H, Hp, R}, and~\textit{Rp}. 

Depending on which topology is chosen the original deterministic model~\cite{steuer2011robust} exhibits different relationships between the steady-state concentrations of the input signal \textit{S} and the output signal~\textit{Rp}: 

\begin{align*}           
Rp \text{ steady-state in model 1} \hspace{2cm}& Rp \text{ steady-state in model 2}\\
[Rp]=\frac{k_1}{k_{31}}[S][H]      \hspace{3cm} &\hspace{1.5cm}   [Rp]=\frac{k_1}{k_{32}}[S]
\end{align*}
In particular, it can be seen that the steady-state concentration of the output signal [\textit{Rp}] in model 1 is affected not only by the input signal \textit{S} but also by the number of unphosphorylated receptors \textit{R}, this can be interpreted in such a way that the concentration of the signalling components should be kept stable in order to obtain the robust output. This is, however, not an issue in model 2 where \textit{Rp} depends only on~\textit{S}.
Since the steady-state analysis has been carried out under the deterministic setting additionally imposing assumptions of conserved total amounts of \textit{H + Hp} and \textit{R + Rp}, it is appropriate only for high molecular populations.

The question we want to answer is ``Is there a difference in the way the two models handle noise (fluctuations) for low molecular numbers of signalling components?'' In such conditions, populations of \textit{H + Hp} and \textit{R + Rp} can not be considered conserved since the proteins are subject to degradation and production. Production of proteins from genes as well as degradation is inherently noisy as it has been demonstrated in the previous case study. Different levels of noise can be affected by, e.g., regulatory feedback loops or varying numbers of gene copies. Even for a noiseless output signal \textit{S} these internal fluctuations of protein concentrations transfer noise to \textit{Rp}. We formalise our question in terms of the CSL property $\mathsf{E}_{=?}[\mathsf{I}^{=t}]$ which asks for the value of a post-processing function in a future time \textit{t}, where the post-processing function is defined as the \textit{mean quadratic deviation} of the distribution of \textit{Rp}.

For the model to have low numbers of molecules to exhibit stochastic fluctuations and enable responses to varying levels of \textit{S} we have chosen $k_{p} = 0.3$ \textit{molecules}$\cdot s^{-1}$ and $k_{d} = 0.01 \mbox{~} s^{-1}$  which leads to an average total population of 30 molecules for both $H + Hp$ and $R + Rp$. To make the analysis straightforward we assume same speed of degradation of phosphorylated and unphosphorylated variants of each protein.

\begin{figure} 
\centering
\includegraphics[width=\textwidth]{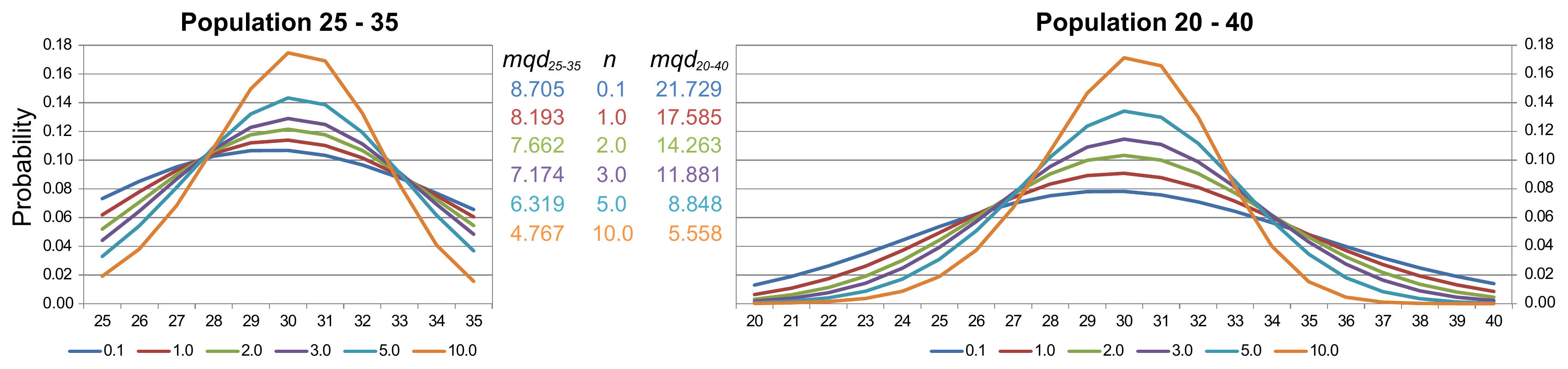}
\caption{{\bf Influence of state space truncation to mean quadratic deviation of a distribution.}
A simple birth death model is considered to show the influence of different settings of the state space truncation on the measured noise evaluated in the form of a \textit{mean quadratic deviation} (\textit{mqd}) of the state space distribution. The model has a single species \textit{X} and two reactions $\emptyset \stackrel{sig(0.3,n)}{\longrightarrow} X, X \stackrel{0.01}{\longrightarrow} \emptyset$ which stabilise the population around an average of 30. For different values of the sigmoid coefficient \textit{n} we can see different \textit{mqd} values, the larger the \textit{n} the smaller the noise. If \textit{X} is restricted to $25\leq X\leq 35$ the overall noise is smaller since the probability mass can not spread to states placed further from the mean. In a less restricted version with populations between 20 and 40 the noise is about $2.5\times$ larger. If sigmoid regulation is weak and the regulation is strong then the difference in the amount of noise is less then 20\%.}
\label{fig:cs2_truncation_noise}
\end{figure}

To reduce the size of the state space we have truncated total populations to $25 \leq H + Hp \leq 35$ and $25 \leq R + Rp \leq 35$ which leads to $116281$ states in total. The initial state is considered with populations $s_0 = (H = 30, Hp = 0, R = 30, Rp = 0)$. The state space reduction has a significant impact on the measured absolute values of noise but conserves general trends as is shown in Figure~\ref{fig:cs2_truncation_noise}.

In order to control fluctuations in protein production we extend our model with two populations of genes, one for \textit{H} and one for \textit{R}, respectively, and for each of the genes we introduce an autoregulatory negative feedback loop via binding of the proteins to their corresponding genes. That way we restrict the protein production. By modifying the number of gene copies in the cell and the rate of protein-gene binding we are able to regulate the overall noise in the transcription. This approach however leads to rapid increase in state space size because of the necessary introduction of new variables representing genes and protein-gene complexes thus making the analysis inefficient. To this end, we decided to abstract from details of the underlying autoregulatory mechanism and to model it using a sigmoid production function which mimics the desired behaviour accordingly. By numerical analysis, we have verified that such an approximation can be employed in the stochastic framework. 
The function is defined in the following way:
$$
\emptyset \stackrel{sig(k_p,n)}{\longrightarrow} X \hspace{4em} sig(k_p,n) = \frac{2}{1 + \left( \frac{X}{30} \right)^n} \cdot k_p
$$ where $n$ is the so-called Hill coefficient controlling the steepness of the sigmoid (caused by cooperativity of transcription factors in protein-gene interactions) and $k_P$ is the maximal production rate.
We use this approach for modelling the production of both species \textit{H} and \textit{R} by sigmoid coefficients denoted $n_H$ and $n_R$, respectively. The sigmoid function regulates the population by enabling production when it is below average and represses it when the population is above the average. The larger \textit{n} is the more steep the sigmoid function is leading to stronger regulation and lower noise. The case \textit{n=0} corresponds to an unregulated model and when increased to \textit{n=20} it corresponds to over 10 copies of each gene in the fully modelled feedback loop mechanism. The effect of different levels of sigmoid regulation to noise can be seen in a simplified birth death model in Figure~\ref{fig:cs2_truncation_noise}.

To see long term effects of intrinsic noise we decided to examine the system in the situation when the output response is stabilised. Since the min-max approximation method cannot be employed with steady-state computation, transient analysis in a suitable time horizon has been performed instead. To estimate the closest time \textit{t} when the system behaviour can be observed stable, we have computed values of output response noise for the unregulated variant of the model (\textit{n = 0}) using standard numerical steady state numerical analysis (we employed the tool PRISM~\cite{KNP11}) and compare it to probability distributions obtained by transient analysis in $t=20$, $t=50$ and $t=100$ seconds. Consequently, we have compared the probability distribution in the steady state with the probability distribution in $t=100$ seconds. The results clearly show that that the difference in distributions is negligible and the transient distribution can be considered stable after $t=100$.

To further speed up the computation, we have precomputed the distribution of \textit{H} and \textit{R} in the time horizon $t=100$ without enabling phosphorylation reactions. This has lead to a significant reduction to $121$ states. Starting with the achieved probability distribution, we have subsequently computed the transient analysis with enabled phosporylation reactions in next $5$ seconds. The rationale behind is that the protein production and degradation are two orders of magnitude slower than phosphorylation. Therefore total populations of \textit{H} and \textit{R} dictate the time at which the system is nearly stable and thus the next $5$ seconds are sufficient for the fast-scale phosporylation to stabilise the fractions $\frac{H}{Hp}$ and $\frac{R}{Rp}$.

To compute the noise (variance) in \textit{Rp} we employ the \textit{mean quadratic deviation} post-processing function for state space distributions. 
Our goal is to compare the levels of \textit{Rp} noise in both models for different levels of the output signal \textit{S} and for different values of intrinsic noise appearing in protein production (controlled by sigmoid coefficients $n_H$ and $n_R$). After computing lower and upper bounds of the state space distributions, we have computed the lower and upper bounds of the post-processing function using the algorithm informally introduced in Section~\ref{sec:parameteriseduniformisation}. Consequently, we obtain robustness values for the output response $R_p$ over the respective perturbation subspaces in the form \textit{average} $\pm$ \textit{error}. Finally, we define the perturbation space of the interest. In particular, for the signal we choose the value interval $S \in [2.0, 20.0]$ and for sigmoid coefficients $n_H, n_R \in [0.1, 10.0]$.

\begin{figure} 
\centering
\includegraphics[width=\textwidth]{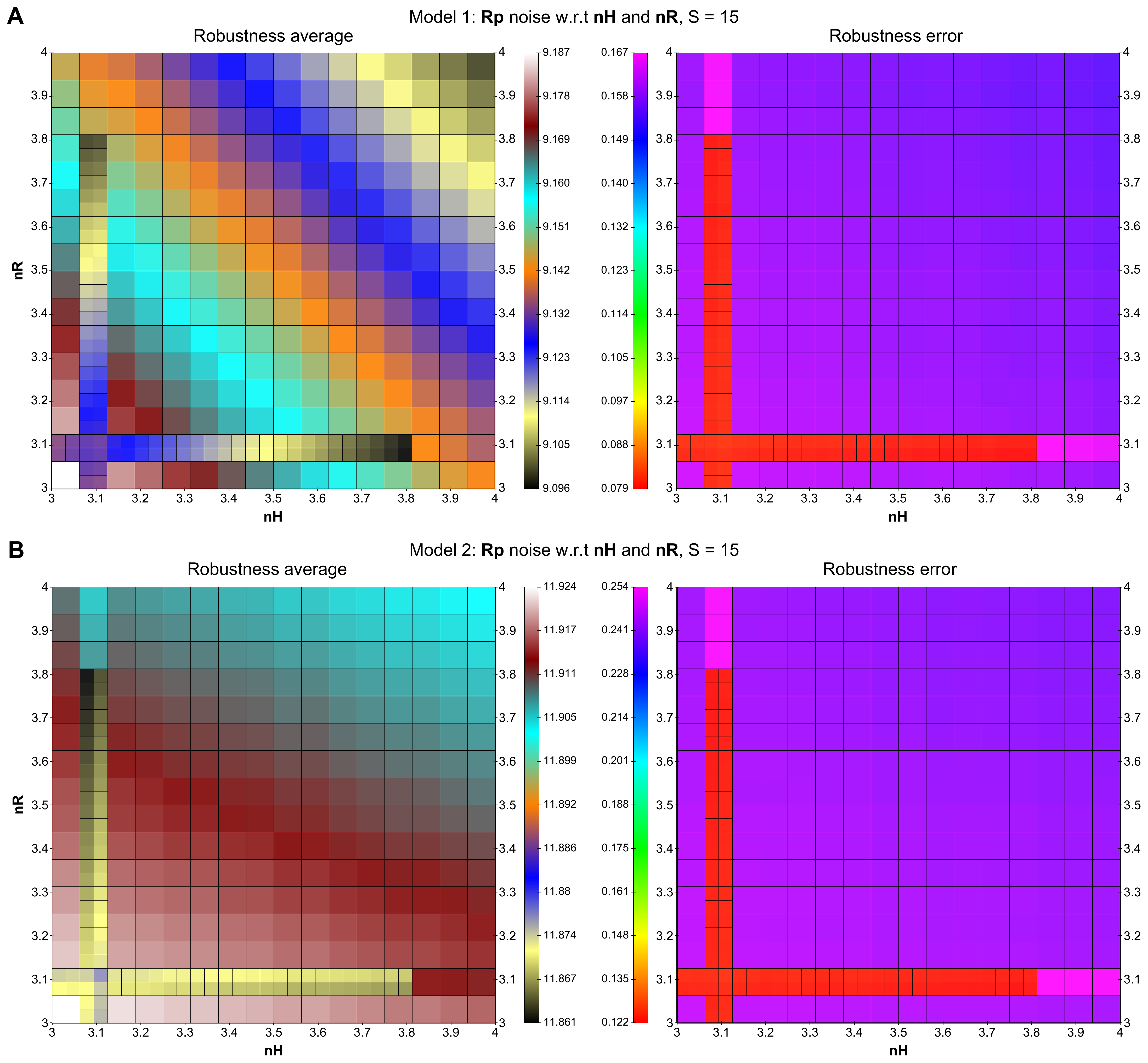}
\caption{{\bf Influence of genetic regulation on noise in model 1 and 2.}
In the upper part two schemes noise of \textit{Rp} in model 1 is computed over perturbations of both sigmoid production constants $n_H$ and $n_R$ in $[3.0,4.0] \times [3.0,4.0]$. The upper and lower bounds on noise (mean quadratic deviation of the resulting probability distribution projected onto populations of \textit{Rp}) are recomputed into the form \textit{average} $\pm$ \textit{error}, the average values are shown on the left and errors are shown on the right. The densely subdivided subspaces around the value $3.1$ are due to conservative over/under approximations in the computation of the probability distribution in states where inflow and outflow of the probability mass is not strictly a monotonous function over the given perturbation interval, thus the error is locally increased and the subspaces must be further divided to obtain the required precision.
The lower two schemes show the same results for model 2.
By comparing both results we can see that model 1 has an overall lower noise and also computation error given the same level of refinement then model 2, in model 1 the results are  symmetrical with respect to perturbations in $n_H$ and $n_R$ with $n_H$ having a slightly larger influence. In model 2 $n_R$ has a larger influence, however we considered the difference negligible and combined both parameters into a single sigmoid production constant \textit{n}.}
\label{fig:cs2_models_nhnr}
\end{figure}

Since the full computation over the 3-dimensional perturbation space has turned out to be intractable, we have to find a way how to reduce its dimension. To this end, we focus on a subspace $S = 15.0, (n_H,n_R) \in [3.0,4.0] \times [3.0,4.0]$ where both models have symmetric sensitivity to both sigmoid production coefficients $n_H, n_R$. This symmetry allows us to merge $n_H, n_R$ into a single coefficient \textit{n}.
Results are visualised in Figure~\ref{fig:cs2_models_nhnr} where it can be seen that in Model 1 the influence of $n_H$ and $n_R$ is almost perfectly symmetrical with $n_H$ being slightly more influential. In Model 2 the influence is evidently stronger in $n_R$ but the response seems to be symmetrical enough to justify the sigmoid coefficients merging. An interesting property of parameterised uniformization and the perturbation space decomposition algorithm can be seen in Figure~\ref{fig:cs2_models_nhnr} where the decomposition of the perturbation spaces around both sigmoid coefficients set to~$3.1$ is very dense. This is due to the non-linearity of the sigmoid production functions which leads to non-monotonicity of probability inflow/outflow differences in states during parameterised uniformization (see Section~\ref{sec:methods}). In order to preserve conservativeness of estimates we have to locally over/under approximate these inflow/outflow rates thus leading to increase of error. To obtain the desired level of accuracy, we dynamically refine all those subspaces where this has occurred.

\begin{figure} 
\centering
\includegraphics[width=\textwidth]{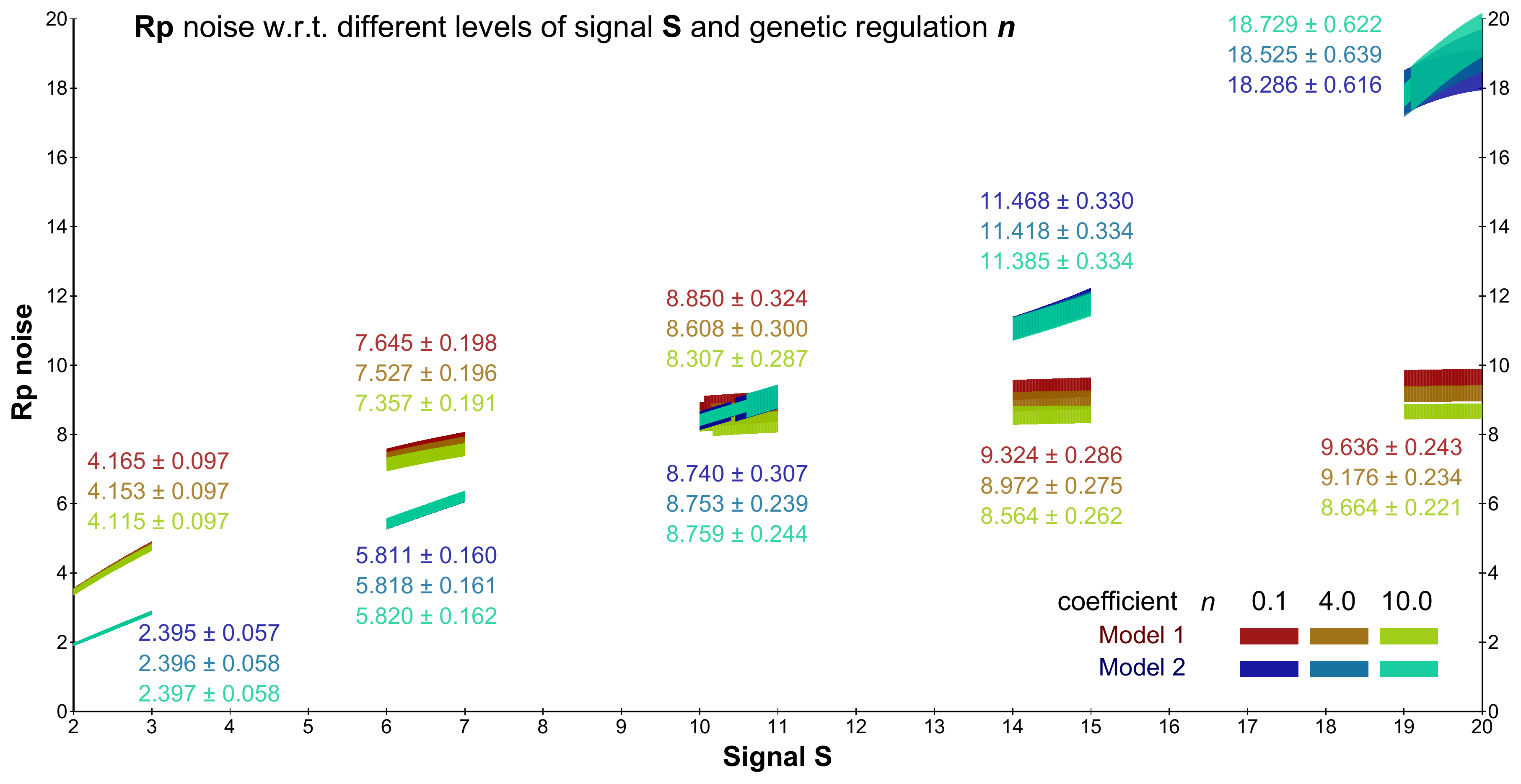}
\caption{{\bf Comparison of models by \textit{Rp} noise robustness.}
Robustness \textit{Rp} noise in both models has been computed with respect to perturbations of signal \textit{S} over five selected intervals of the input signal $S \in [2,3] \cup [6,7] \cup [10,11] \cup [14,15] \cup [19,20]$ and for three distinct levels of the intrinsic noise in signalling component dynamics represented by sigmoid coefficient $n \in \left\{0.1, 4.0, 10.0\right\}$. Perturbations were not computed over the whole interval $(S,n) \in [2, 20] \times [0.1, 10.0]$ due to very high computational demands. From the computed values of individual refined subspaces as well as the aggregated robustness values for each input signal interval we can see that for lower values of signal \textit{S} (up-to 10) Model 2 embodies lower output response noise then Model 1 (spontaneous dephosphorylation). While output response noise in Model 1 tends to converge to values between 8 and 10, Model 2 exhibits a permanent (almost linear) increase in the output response noise over most of the studied portion of the perturbation space. A super-linear increase of the noise is observed for strong input signals. Another interesting aspect is that while with increasing levels of gene regulation given by sigmoid coefficient \textit{n} the overall noise in \textit{Rp} decreases over the whole interval of signal values for Model 1 and most of it for Model 2, there is an anomaly in Model 2 in the high signal region [19.0, 20.0] where with decreasing noise in \textit{R} and \textit{H} (see Figure~\ref{fig:cs2_model12_hr_noise}) the noise in \textit{Rp} increases. We have not yet explained this phenomenon satisfactorily.}
\label{fig:cs2_model12_rp_noise}
\end{figure}

\begin{figure} 
\centering
\includegraphics[width=\textwidth]{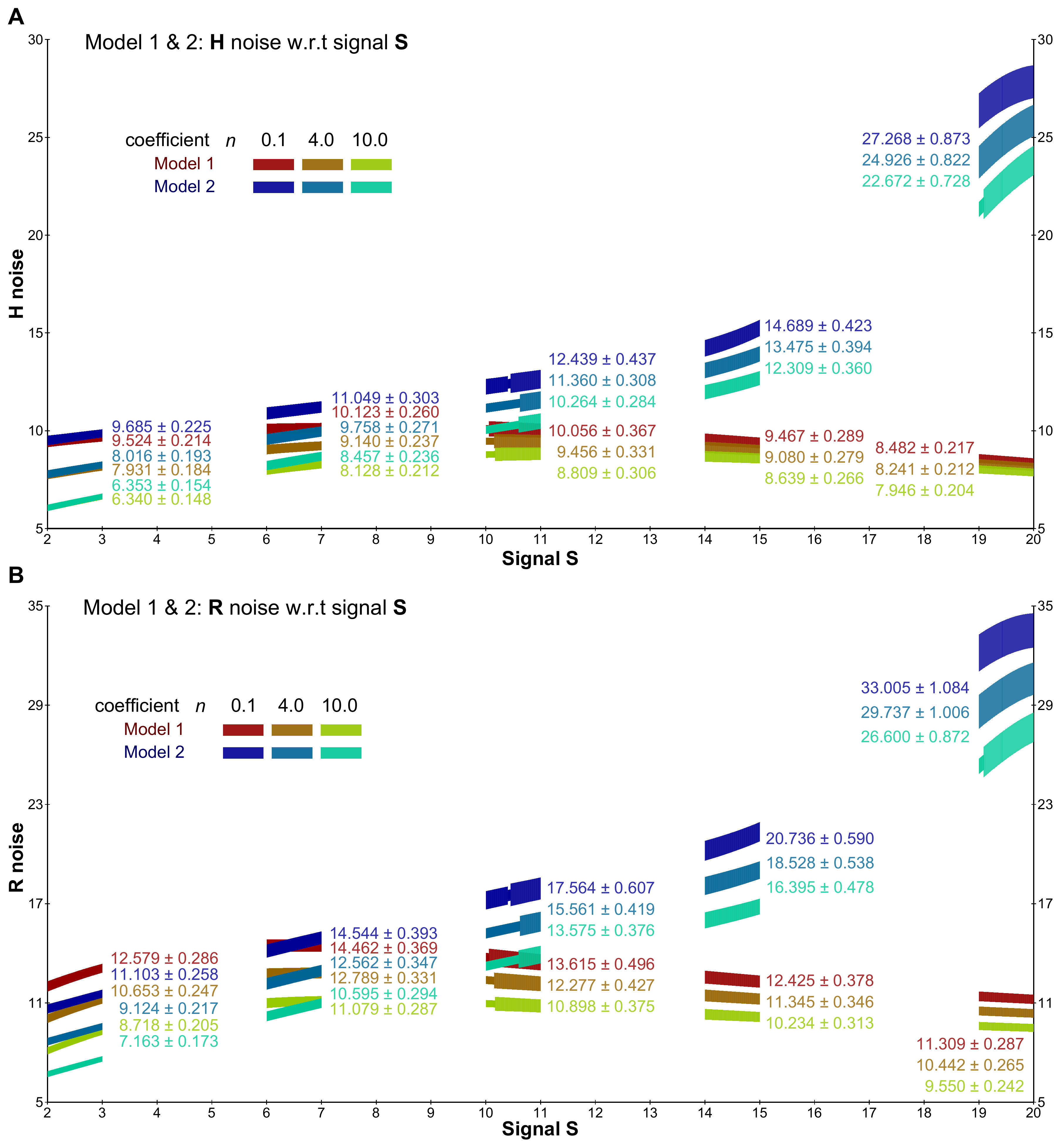}
\caption{{\bf Noise in populations or \textit{H} and \textit{R} in both models.}
Noise in \textit{H} (A) and \textit{R} (B) in both models has been computed with respect to perturbations of signal \textit{S} over five selected intervals $S \in [2,3] \cup [6,7] \cup [10,11] \cup [14,15] \cup [19,20]$ and for three distinct levels of inherent production noise represented by sigmoid coefficient $n \in \left\{0.1, 4.0, 10.0\right\}$. We can see that in all cases with increasing regulation by \textit{n} the intrinsic noise in the dynamics of each of the signalling components decreases.}
\label{fig:cs2_model12_hr_noise}
\end{figure}

\begin{figure} 
\centering
\includegraphics[width=\textwidth]{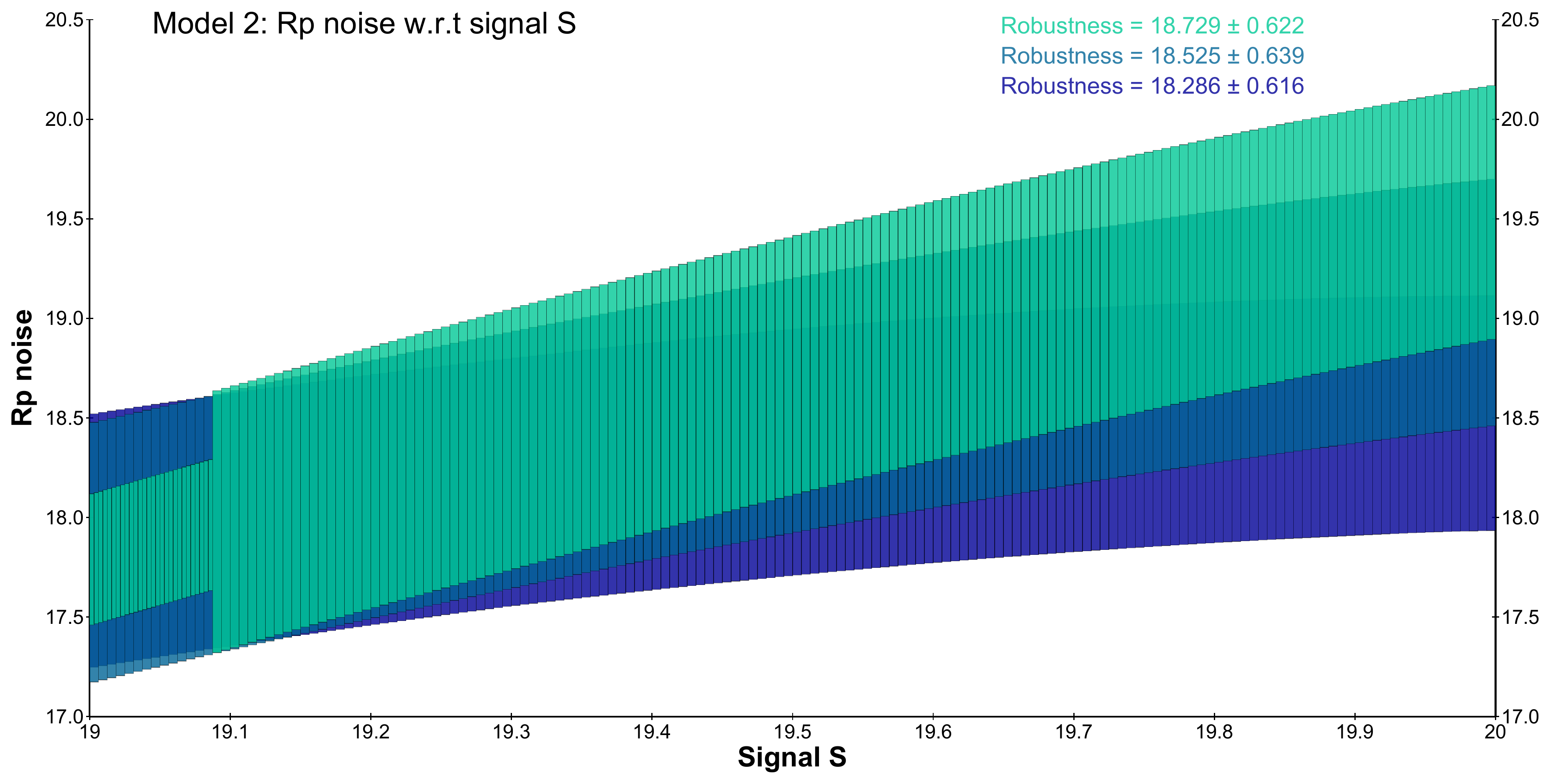}
\caption{{\bf High signal region in model 2.}
A closeup of the high signal region in model 2, where increasing levels of regulation by the sigmoid coefficient \textit{n} leads to a paradoxical increase of output response noise instead of decrease. Even though the inaccuracy is large we consider the trend to be strong and thus real.}
\label{fig:cs2_model2highsignal}
\end{figure}

Finally, we inspect selected subintervals of the perturbation space given by five exclusive intervals of the input signal value domain, $S \in [2,3] \cup [6,7] \cup [10,11] \cup [14,15] \cup [19,20]$, and three distinct levels of production noise represented by sigmoid coefficient $n \in \left\{0.1, 4.0, 10.0\right\}$. The results of this main experiment can be seen in Figure~\ref{fig:cs2_model12_rp_noise} and Figure~\ref{fig:cs2_model12_hr_noise}. The trends that can be seen in Figure~\ref{fig:cs2_model12_rp_noise} are that for lower signals up to \textit{S=10}. Model 2 has encountered lower noise in \textit{Rp} than Model 1 but in the higher signal region it is outperformed by Model 1 which quickly converges to values between $8$ and $10$. However, \textit{Rp} noise produced in Model 2 linearly increases with increasing value of the input signal \textit{S}. For most of the inspected subspaces a stronger regulation of \textit{H} and \textit{R} production by the sigmoid coefficient \textit{n} leads to a reduction of \textit{Rp} noise. An exception to this observation can be seen in Model 2 at the signal interval $[19.0, 20.0]$ where this trend is inverted. To show that this is an emergent behaviour arising from the nontrivial interaction of phosphorylation and dephosphorylation reactions not present in the basic production and degradation of components \textit{H} and \textit{R}, their respective influences are displayed in Figure~\ref{fig:cs2_model12_hr_noise}. There we can see that in Model 1 both \textit{H} and \textit{R} follow an initial increase of noise with increasing \textit{S} but then the noise stabilises. This leads us to a hypothesis that the regulation of noise in signalling components dynamics looses its influence as signal \textit{S} increases. This is however due to the fact that more \textit{S} leads to faster phosphorylation of \textit{H} which effectively reduces the population of \textit{H} thus also reducing its absolute noise. In the case of Model 2 the situation is different since we can observe a permanent increase of noise in both \textit{H} and \textit{R} populations. The inversion of noise with increased regulation seen in~\ref{fig:cs2_model12_rp_noise} and closely shown in Figure~\ref{fig:cs2_model2highsignal} has not yet been explained satisfactorily.

\section{Discussion}

In this paper we proposed a novel framework for robustness analysis of stochastic biochemical systems. It allows us to quantify and analyse how the validity of a hypothesis formulated as a temporal property depends on the perturbations of stochastic kinetic parameters and initial concentrations.  The framework extends the quantitative model checking techniques and numerical methods for CTMCs and adapts them to the needs of stochastic modelling in biology. Therefore, in contrast to statistical methods such as Monte Carlo simulation and parameter sampling our framework is customizable with respect to the required precision of computation. This is obtained by providing the lower and upper bounds of the results.

Case studies have demonstrated that the framework can be successfully applied to the robustness analysis of nontrivial biochemical systems. They have shown how to use CSL to specify properties targeting transient behaviour under fluctuations. From the first case study we can conclude that the reward-based formulation of stability properties is more appropriate to distinguish the individual parameter settings under the requested range of uncertainty. The inspected biological hypothesis in the second case study can not be directly formulated using CSL with rewards. Therefore, we have employed post-processing functions to express and study the mean quadratic deviation of the molecule population distribution of the signal response regulator protein.

The time complexity of our framework in practice depends mainly on the size of the state space, the number of reaction steps that have to be considered, and the number of perturbation sets that have to be analysed to provide the desired precision. The size of the state space is given by the number of species and their populations. The framework is suitable for low populations and is relevant especially in the case of gene regulation. In the first case study we have considered only a single molecule of DNA and thus the state space of resulting CTMC was manageable. In the second case study we had to abstract from the feedback loop mechanism using a sigmoid production function to reduce the state space and to make the analysis feasible. If such an abstraction can not be used, our framework can be effectively combined with general state space reduction methods for CTMCs, e.g., finite projection techniques~\cite{munsky2006finite,Henzinger2009} and dynamic state space truncation~\cite{Didier}. The number of reaction steps can be reduced using separation of fast and slow reactions as demonstrated in the second case study or using adaptive uniformisation~\cite{Moorsel94, Didier}.

In the first case study several hundreds of perturbation subsets had to be analysed and the overall robustness analysis took a few hours. However, in the second case study several thousands of perturbation subsets were required to achieve reasonable precision. In order to speedup the computation we analysed the subsets in parallel using a high performance multi-core workstation were the analysis took several hours. To further improve the accuracy of the robustness analysis without decreasing the performance, we have employed a piecewise linear approximation. It allows us to obtain more precise result without increasing the number of perturbation sets, however, it does not guarantee the conservative error bounds.

The presented method as employed in the first case study gives us a tool for exact analysis of bistability from the global point of view (with respect to all initial conditions, the considered time bound, and the given range of parameters). It can be considered as an analogy to bifurcation analysis known from the ODE world. When comparing our approach with the bifurcation analysis performed in~\cite{Swatetal04}, our approach provides a detailed mesoscopic insight into the analysed phenomenon. Instead of identifying just the points where the population diverges, we obtain the precise knowledge of how the population is distributed around the two stable states. Especially, the method shows that reachability of the cancer-inducing high stable mode of the retinoblastoma-binding transcription factor is almost always possible despite the initial state of the regulatory system. The exhaustive analysis is performed with uncertainty in the degradation parameters of the two most important cell-cycle regulating proteins. However, if the degradation of the tumour suppressor protein is sufficiently high, there is always possibility allowing the population to switch into the safe low stable mode. Moreover, robustness of having the possibility to avoid the cell malfunction is positively affected by increasing the retinoblastoma-binding transcription factor degradation. In contrast to~\cite{Swatetal04}, the switching mechanism is described at the single cell level which allows to quantify the portion of population amenable to mall-function and thus can provide a preliminary guide to further analysis targeting elimination of the undesired behaviour.      

The second case study has shown new insights into the phenomenon of noise in two-component signalling pathways appearing in procaryotic organisms. The previous study~\cite{steuer2011robust} conducted in the framework of deterministic models targeted global robustness of steady concentrations of output signalling components by means of analytically finding invariant perturbation space. The result has shown that a synthetic pathway topology including additional catalysis of signal response regulator by histidine kinase leads to globally robust input-output signal mapping with respect to fluctuations in signalling components concentration. On the contrary, the basic topology without histidine-modulated dephosphorylation does not fulfil global robustness. Since signalling pathways are understood to be amenable to intrinsic noise due to relatively low molecule populations of signalling proteins (typically hundreds of molecules), the respective stochasticity might affect the input-output signal response. To this end, we have reformulated the model in the stochastic framework and instead of studying the effect of perturbations on the average population, we study in detail how perturbations affect the distribution, i.e., the variance (fluctuation) in the output response. Our study has shown that both pathway topologies result with fluctuations in output response, but robustness of input-output mapping varies in both models with increasing the level of the (constant) input signal. For low input signals the synthetic topology gives response with smaller variance in the output whereas for high input signals the output variance rapidly increases. Therefore the basic topology seems to be more suitable for processing of strong signals while the synthetic topology is more appropriate for low level signals. Our study has also shown that both topologies are quite robust with respect to scaling the noise in signalling components dynamics.

\section*{Acknowledgments}
This work has been supported by the Czech Science Foundation
    grant No. GAP202/11/0312. M. \v{C}e\v{s}ka has been supported by Ministry of Education, Youth, and Sport project No. CZ.1.07/2.3.00/30.0009~-~Employment of Newly Graduated Doctors of Science for Scientific Excellence. D. \v{S}afr\'anek has been supported by EC OP project No. CZ.1.07/2.3.00/20.0256.

\bibliography{bibliography}

\begin{thebibliography}{10}

\bibitem{Verena2011}
Aleksandr Andreychenko, Linar Mikeev, David Spieler, and Verena Wolf.
\newblock {Parameter Identification for Markov Models of Biochemical
  Reactions}.
\newblock In {\em Computer Aided Verification}, LNCS, pages 83--98. Springer,
  2011.

\bibitem{Aziz1996csl}
Adnan Aziz, Kumud Sanwal, Vigyan Singhal, and Robert Brayton.
\newblock {Verifying continuous time Markov chains}.
\newblock In {\em Computer Aided Verification}, volume 1102 of {\em LNCS},
  pages 269--276. Springer, 1996.

\bibitem{Baier2003mcctmc}
C.~Baier, B.~Haverkort, H.~Hermanns, and J.P. Katoen.
\newblock {Model-checking algorithms for continuous-time Markov chains}.
\newblock {\em IEEE Transactions on Software Engineering}, 29(6):524--541,
  2003.

\bibitem{Baier2000mcviatran}
Christel Baier, Boudewijn Haverkort, Holger Hermanns, and Joost-Pieter Katoen.
\newblock {Model Checking Continuous-Time Markov Chains by Transient Analysis}.
\newblock In {\em Computer Aided Verification}, volume 1855 of {\em LNCS},
  pages 358--372. Springer, 2000.

\bibitem{Paolo}
Paolo Ballarini, Michele Forlin, Tommaso Mazza, and Davide Prandi.
\newblock {Efficient Parallel Statistical Model Checking of Biochemical
  Networks}.
\newblock In {\em Parallel and Distributed Methods in verifiCation}, volume~14
  of {\em EPTCS}, pages 47--61, 2009.

\bibitem{Barnat12coloredmc}
J.~Barnat, L.~Brim, A.~Krej\v{c}\'{i}, A.~Streck, D.~\v{S}afr\'{a}nek,
  M.~Vejn\'{a}r, and T.~Vejpustek.
\newblock On parameter synthesis by parallel model checking.
\newblock {\em IEEE/ACM Transactions on Computational Biology and
  Bioinformatics}, 9(3):693 --705, may-june 2012.

\bibitem{BBS10}
J.~Barnat, L.~Brim, and D.~\v{S}afr\'{a}nek.
\newblock {High-Performance Analysis of Biological Systems Dynamics with the
  DiVinE Model Checker}.
\newblock {\em Briefings in Bioinformatics}, 11(3):301--312, 2010.

\bibitem{Bartocci2013}
E.~{Bartocci}, L.~{Bortolussi}, L.~{Nenzi}, and G.~{Sanguinetti}.
\newblock {On the Robustness of Temporal Properties for Stochastic Models}.
\newblock {\em ArXiv e-prints}, September 2013.

\bibitem{Batchelor21012003}
Eric Batchelor and Mark Goulian.
\newblock Robustness and the cycle of phosphorylation and dephosphorylation in
  a two-component regulatory system.
\newblock {\em Proceedings of the National Academy of Sciences},
  100(2):691--696, 2003.

\bibitem{Bernardini}
F.~Bernardini, C.~Biggs, J.~Derrick, M.~Gheorghe, M.~Niranjan, and
  G.~Sanguinetti.
\newblock {Parameter Estimation and Model Checking in a Model of Prokaryotic
  Autoregulation}.
\newblock Technical report, University of Sheffield, 2007.

\bibitem{Bortolussi2012}
Luca Bortolussi and Jane Hillston.
\newblock Fluid model checking.
\newblock In Maciej Koutny and Irek Ulidowski, editors, {\em CONCUR 2012 –
  Concurrency Theory}, volume 7454 of {\em Lecture Notes in Computer Science},
  pages 333--347. Springer Berlin Heidelberg, 2012.

\bibitem{CAV2013}
Lubo\v{s} Brim, Milan \v{C}e\v{s}ka, Sven Dra\v{z}an, and David
  \v{S}afr\'{a}nek.
\newblock Exploring parameter space of stochastic biochemical systems using
  quantitative model checking.
\newblock In {\em Computer Aided Verification}, volume 8044 of {\em LNCS},
  pages 107--123. Springer Berlin Heidelberg, 2013.

\bibitem{Petzold2012}
Bernie Daigle, Min Roh, Linda Petzold, and Jarad Niemi.
\newblock {Accelerated Maximum Likelihood Parameter Estimation for Stochastic
  Biochemical Systems}.
\newblock {\em BMC Bioinformatics}, 13(1):68--71, 2012.

\bibitem{Didier}
Frederic Didier, Thomas~A. Henzinger, Maria Mateescu, and Verena Wolf.
\newblock {Fast Adaptive Uniformization of the Chemical Master Equation}.
\newblock In {\em High Performance Computational Systems Biology}, pages
  118--127. IEEE Computer Society, 2009.

\bibitem{Donze2010breach}
Alexandre Donz\'{e}.
\newblock {Breach, A Toolbox for Verification and Parameter Synthesis of Hybrid
  Systems}.
\newblock In {\em Computer Aided Verification}, volume 6174 of {\em LNCS},
  pages 167--170. Springer, 2010.

\bibitem{Donze2010robust}
Alexandre Donz\'{e} and Oded Maler.
\newblock {Robust satisfaction of temporal logic over real-valued signals}.
\newblock In {\em Formal Modeling and Analysis of Timed Systems}, volume 6246
  of {\em LNCS}, pages 92--106. Springer, 2010.

\bibitem{Donze11rob-behaviour}
Alexandre Donzé, Eric Fanchon, Lucie~Martine Gattepaille, Oded Maler, and
  Philippe Tracqui.
\newblock Robustness analysis and behavior discrimination in enzymatic reaction
  networks.
\newblock {\em PLoS ONE}, 6(9):e24246, 09 2011.

\bibitem{Gillespieetal05}
Hana El~Samad, Mustafa Khammash, Linda Petzold, and Dan Gillespie.
\newblock Stochastic modelling of gene regulatory networks.
\newblock {\em International Journal of Robust and Nonlinear Control},
  15(15):691--711, 2005.

\bibitem{Fages2004biocham}
F.~Fages, S.~Soliman, and N.~Chabrier-Rivier.
\newblock Modelling and querying interaction networks in the biochemical
  abstract machine biocham.
\newblock {\em Journal of Biological Physics and Chemistry}, 4(2):64--73, 2004.

\bibitem{Fainekos2009robustness}
G.E. Fainekos and G.J. Pappas.
\newblock Robustness of temporal logic specifications for continuous-time
  signals.
\newblock {\em Theoretical Computer Science}, 410(42):4262--4291, 2009.

\bibitem{FoxGlynn1988}
Bennett~L. Fox and Peter~W. Glynn.
\newblock {Computing Poisson probabilities}.
\newblock {\em Commun. ACM}, 31(4):440--445, April 1988.

\bibitem{Celldesigner}
Akira Funahashi, Mineo Morohashi, Hiroaki Kitano, and Naoki Tanimura.
\newblock Celldesigner: a process diagram editor for gene-regulatory and
  biochemical networks.
\newblock {\em BIOSILICO}, 1(5):159 -- 162, 2003.

\bibitem{Garaietal12}
Ashok Garai, Bartlomiej Waclaw, Hannes Nagel, and Hildegard Meyer-Ortmanns.
\newblock Stochastic description of a bistable frustrated unit.
\newblock {\em Journal of Statistical Mechanics: Theory and Experiment},
  2012(01):P01009, 2012.

\bibitem{Gillespie1977}
Daniel~T. Gillespie.
\newblock {Exact Stochastic Simulation of Coupled Chemical Reactions}.
\newblock {\em Journal of Physical Chemistry}, 81(25):2340--2381, 1977.

\bibitem{Wilkinson}
Andrew Golightly and Darren~J. Wilkinson.
\newblock {Bayesian Parameter Inference for Stochastic Biochemical Network
  Models Using Particle Markov Chain Monte Carlo}.
\newblock {\em Interface Focus}, 1(6):807--820, 2011.

\bibitem{Verena2013}
J.~Hasenauer, V.~Wolf, A.~Kazeroonian, and F.J. Theis.
\newblock Method of conditional moments (mcm) for the chemical master equation.
\newblock {\em Journal of Mathematical Biology}, pages 1--49, 2013.

\bibitem{Henzinger2009}
Thomas~A. Henzinger, Maria Mateescu, and Verena Wolf.
\newblock {Sliding Window Abstraction for Infinite Markov Chains}.
\newblock In {\em Computer Aided Verification}, volume 5643 of {\em LNCS},
  pages 337--352. Springer, 2009.

\bibitem{Hill1910}
Archibald~Vivian Hill.
\newblock The possible effects of the aggregation of the molecules of
  hamoglobin on its dissociation curves.
\newblock {\em The Journal of Physiology}, 40(Suppl):iv--vii, 1910.

\bibitem{Hoops2006copasi}
Stefan Hoops, Sven Sahle, Ralph Gauges, Christine Lee, J\"{u}rgen Pahle,
  Natalia Simus, Mudita Singhal, Liang Xu, Pedro Mendes, and Ursula Kummer.
\newblock {COPASI - a COmplex PAthway SImulator}.
\newblock {\em Bioinformatics}, 22(24):3067--3074, 2006.

\bibitem{Clarke}
Sumit~K. Jha, Edmund~M. Clarke, Christopher~J. Langmead, Axel Legay, Andr{\'e}
  Platzer, and Paolo Zuliani.
\newblock {A Bayesian Approach to Model Checking Biological Systems}.
\newblock In {\em Computational Methods in Systems Biology}, pages 218--234.
  Springer, 2009.

\bibitem{Keletal00}
Alexander~E. Kel, Igor Deineko, Olga~V. Kel-Margoulis, Edgar Wingender, and
  Vadim Ratner.
\newblock Modeling of gene regulatory network of cell cycle control. role of
  e2f feedback loops.
\newblock In {\em German Conference on Bioinformatics'00}, pages 107--114,
  2000.

\bibitem{Kitano2004}
H.~Kitano.
\newblock Biological robustness.
\newblock {\em Nat Rev Genet}, 5:826--837, 2004.

\bibitem{Kitano2007}
Hiroaki Kitano.
\newblock {Towards a theory of biological robustness}.
\newblock {\em Molecular Systems Biology}, 3:137, 2007.

\bibitem{Cago}
Chuan~Hock Koh, Sucheendra Palaniappan, PS~Thiagarajan, and Limsoon Wong.
\newblock {Improved Statistical Model Checking Methods for Pathway Analysis}.
\newblock {\em BMC Bioinformatics}, 13(Suppl 17):S15, 2012.

\bibitem{Kwiatkowska2006rewards}
M.~Kwiatkowska, G.~Norman, and A.~Pacheco.
\newblock {Model Checking Expected Time and Expected Reward Formulae with
  Random Time Bounds}.
\newblock {\em Computers \& Mathematics with Applications}, 51(2):305 -- 316,
  2006.

\bibitem{KNP11}
M.~Kwiatkowska, G.~Norman, and D.~Parker.
\newblock {PRISM} 4.0: Verification of probabilistic real-time systems.
\newblock In {\em Computer Aided Verification}, volume 6806 of {\em LNCS},
  pages 585--591. Springer, 2011.

\bibitem{Kwiatkowska2007}
Marta Kwiatkowska, Gethin Norman, and David Parker.
\newblock {Stochastic model checking}.
\newblock In {\em Formal Methods for Performance Evaluation}, volume 4486 of
  {\em LNCS}, pages 220--270. Springer, 2007.

\bibitem{KwiatkowskaBIO}
Marta~Z. Kwiatkowska, Gethin Norman, and David Parker.
\newblock {Using Probabilistic Model Checking in Systems Biology}.
\newblock {\em SIGMETRICS Performance Evaluation Review}, 35(4):14--21, 2008.

\bibitem{Loew2001vcell}
L.M. Loew and J.C. Schaff.
\newblock {The Virtual Cell: a software environment for computational cell
  biology}.
\newblock {\em Trends in biotechnology}, 19:401--406, 2001.

\bibitem{Madsen2012}
C.~Madsen, C.J. Myers, N.~Roehner, C.~Winstead, and Zhen Zhang.
\newblock Utilizing stochastic model checking to analyze genetic circuits.
\newblock In {\em Computational Intelligence in Bioinformatics and
  Computational Biology}, pages 379--386. IEEE Computer Society, 2012.

\bibitem{Verena2012}
Linar Mikeev, MartinR. Neuh{\"a}u{\ss}er, David Spieler, and Verena Wolf.
\newblock {On-the-fly Verification and Optimization of DTA-properties for Large
  Markov Chains}.
\newblock {\em Form. Method. Syst. Des.}, pages 1--25, 2012.

\bibitem{munsky2006finite}
Brian Munsky and Mustafa Khammash.
\newblock The finite state projection algorithm for the solution of the
  chemical master equation.
\newblock {\em The Journal of chemical physics}, 124:044104, 2006.

\bibitem{Reinker}
S.~Reinker, R.M. Altman, and J.~Timmer.
\newblock {Parameter Estimation in Stochastic Biochemical Reactions}.
\newblock {\em IEEE Proc. Syst. Biol.}, 153(4):168--78, 2006.

\bibitem{Rizketal09}
A.~Rizk, G.~Batt, F.~Fages, and S.~Soliman.
\newblock A general computational method for robustness analysis with
  applications to synthetic gene networks.
\newblock {\em Bioinformatics}, 25:169--178, 2009.

\bibitem{Rizk09rob-property}
Aurélien Rizk, Gregory Batt, François Fages, and Sylvain Soliman.
\newblock A general computational method for robustness analysis with
  applications to synthetic gene networks.
\newblock {\em Bioinformatics}, 25(12):i169--i178, 2009.

\bibitem{Sanftetal11}
K.R. Sanft, D.T. Gillespie, and L.R. Petzold.
\newblock Legitimacy of the stochastic michaelis-menten approximation.
\newblock {\em Systems Biology, IET}, 5(1):58--69, 2011.

\bibitem{Shinar11122007}
Guy Shinar, Ron Milo, María~Rodríguez Martínez, and Uri Alon.
\newblock Input–output robustness in simple bacterial signaling systems.
\newblock {\em Proceedings of the National Academy of Sciences},
  104(50):19931--19935, 2007.

\bibitem{steuer2011robust}
Ralf Steuer, Steffen Waldherr, Victor Sourjik, and Markus Kollmann.
\newblock Robust signal processing in living cells.
\newblock {\em PLoS computational biology}, 7(11):e1002218, 2011.

\bibitem{Stewart2009uniformization}
William~J. Stewart.
\newblock {\em Probability, Markov Chains, Queues, and Simulation: The
  Mathematical Basis of Performance Modeling}.
\newblock Princeton University Press, Princeton, NJ, USA, 2009.

\bibitem{Swatetal04}
Maciej Swat, Alexander Kel, and Hanspeter Herzel.
\newblock Bifurcation analysis of the regulatory modules of the mammalian g1/s
  transition.
\newblock {\em Bioinformatics}, 20(10):1506--1511, 2004.

\bibitem{ueda2007stochastic}
Masahiro Ueda and Tatsuo Shibata.
\newblock Stochastic signal processing and transduction in chemotactic response
  of eukaryotic cells.
\newblock {\em Biophysical journal}, 93(1):11, 2007.

\bibitem{Moorsel94}
A.~P.~A. van Moorsel and W.~H. Sanders.
\newblock Adaptive uniformization.
\newblock {\em ORSA Communications in Statistics: Stochastic Models, vol. 10,
  no. 3}, pages 619--648, 1994.

\bibitem{Yang2003}
Edward Yang, Erik van Nimwegen, Mihaela Zavolan, Nikolaus Rajewsky, Mar~k
  Schroeder, Marcelo Magnasco, and James~E. Darnell.
\newblock {Decay Rates of Human mRNAs: Correlation With Functional
  Characteristics and Sequence Attributes}.
\newblock {\em Genome Research}, 13(8):1863--1872, 2003.

\bibitem{Alon2004}
Alon Zaslaver, Avi~E. Mayo, Revital Rosenberg, Pnina Bashkin, Hila Sberro, Miri
  Tsalyuk, Michael~G. Surette, and Uri Alon.
\newblock {Just-in-time transcription program in metabolic pathways}.
\newblock {\em Nature Genetics}, 36(5):486--491, April 2004.

\end{thebibliography}

\end{document}